\documentclass[12pt]{amsart}
\usepackage{amsmath}
\usepackage{amsxtra}
\usepackage{amscd}
\usepackage{amsthm}
\usepackage{amsfonts}
\usepackage{amssymb}
\usepackage{eucal}
\usepackage{enumerate}
\usepackage{epsfig}
\usepackage{bm}
\usepackage{graphics}
\usepackage{graphicx}
\usepackage{accents}
\usepackage{dynkin-diagrams}
\usepackage{comment}
\usepackage{enumitem}
\usepackage{mathrsfs}
\usepackage{mathtools}
\usepackage{tikz-cd}
\usepackage{mathtools}
\usepackage{pgf}
\usetikzlibrary{arrows, matrix}
\usepackage{ytableau}
\usepackage{here} 
\usepackage{simplewick}
\usepackage[hidelinks]{hyperref}
\usepackage{cancel}

\numberwithin{equation}{section}
\newtheorem{thm}{Theorem}[section]
\newtheorem{prop}[thm]{Proposition}
\newtheorem{lem}[thm]{Lemma}

\AtBeginEnvironment{rem}{%
  \pushQED{\qed}%
}
\AtEndEnvironment{rem}{\popQED\endexample}

\newtheorem{cor}[thm]{Corollary}
\newtheorem{conj}[thm]{Conjecture}
\newtheorem{dfn}[thm]{Definition}

\theoremstyle{definition}
\newtheorem{example}[thm]{Example}
\AtBeginEnvironment{example}{%
  \pushQED{\qed}%
}
\AtEndEnvironment{example}{\popQED\endexample}

\newcommand{\bra}[1]{\langle #1 |}        
\newcommand{\ket}[1]{{| #1 \rangle}}

\newcommand{\C}{{\mathbb C}}
\newcommand{\Z}{{\mathbb Z}}

\newcommand{\E}{{\mathcal E}}

\newcommand{\bs}{\boldsymbol}

\newcommand{\gl}{\mathfrak{gl}}

\newcommand{\sL}{\mathfrak{sl}}

\newcommand{\la}{\lambda}

\newcommand{\ga}{\gamma}
\newcommand{\Res}{\mathop{\rm res}}

\newcommand{\mc}{\mathcal}
\newcommand{\al}{\alpha}

\newcommand\genmathchar[6]{%
    \newcommand*#1{#3{%
        \mathchoice
            {\hbox{#4\char#2}}%
            {\hbox{#4\char#2}}%
            {\hbox{#5\char#2}}%
            {\hbox{#6\char#2}}%
    }}%
}

\font\sixmathb=mathb6
\font\eightmathb=mathb8
\font\twelvemathb=mathb12
\genmathchar {\boxxslash}     {"6D} {\mathbin} {\twelvemathb} {\eightmathb} {\sixmathb}
\genmathchar {\boxbackslash} {"6E} {\mathbin} {\twelvemathb} {\eightmathb} {\sixmathb}

\usepackage{graphicx}

\let\oldboxxslash\boxxslash
\renewcommand{\boxxslash}{\raisebox{0.2ex}{\ensuremath{\oldboxxslash}}}

\usepackage{manfnt}
\let\oldmancube\mancube
\renewcommand{\mancube}{\scalebox{0.75}{\oldmancube}}

\newcommand{\U}{U}

\newcommand{\be}{\begin{equation*}}
\newcommand{\ee}{\end{equation*}}

\usepackage{ifsym}
\usepackage{manfnt}

\def\squarebox at (#1,#2){\draw (#1,#2) rectangle (#1+0.5, #2 +0.5);}

\def\diagonalbox at (#1,#2){\draw (#1,#2) rectangle (#1+0.5, #2 +0.5);
\draw (#1,#2) -- (#1+0.5, #2+0.5);}

\def\thindiagonalbox at (#1,#2){\draw[ultra thin] (#1,#2) rectangle (#1+0.5, #2 +0.5);
\draw[ultra thin] (#1,#2) -- (#1+0.5, #2+0.5);}

\def\ibox at (#1,#2){\draw (#1,#2) rectangle (#1+1, #2 +1);  \node at (#1+0.5, #2 +0.5) {\tiny{$i$}};}
\def\inextbox at (#1,#2){\draw (#1,#2) rectangle (#1+1, #2 +1);  \node at (#1+0.5, #2 +0.5) {\tiny{$i+1$}};}
\def\iprevbox at (#1,#2){\draw (#1,#2) rectangle (#1+1, #2 +1);  \node at (#1+0.5, #2 +0.5) {\tiny{$i-1$}};}

\def\zerobox at (#1,#2){\draw (#1,#2) rectangle (#1+0.5, #2 +0.5);  \node at (#1+0.25, #2 +0.25) {\tiny{$0$}};}

 \def\redbox at (#1,#2){%
  \fill[red] (#1,#2) rectangle (#1+0.5,#2+0.5);
}

\def\rzerobox at (#1,#2){\draw (#1,#2) rectangle (#1+0.5, #2 +0.5);  \node at (#1+0.25, #2 +0.25) {\textcolor{red}{$0$}};}

\def\abox at (#1,#2){\draw (#1,#2) rectangle (#1+0.5, #2 +0.5);  \node at (#1+0.25, #2 +0.25) {$A$};}

\def\onebox at (#1,#2){\draw (#1,#2) rectangle (#1+0.5, #2 +0.5);  \node at (#1+0.25, #2 +0.25) {\tiny{$1$}};}

\def\twobox at (#1,#2){\draw (#1,#2) rectangle (#1+0.5, #2 +0.5);  \node at (#1+0.25, #2 +0.25) {\tiny{$2$}};}

\def\threebox at (#1,#2){\draw (#1,#2) rectangle (#1+0.5, #2 +0.5);  \node at (#1+0.25, #2 +0.25) {\tiny{$3$}};}

\def\dotsbox at (#1,#2){\draw (#1,#2) rectangle (#1+0.5, #2 +0.5);  \node at (#1+0.27, #2 +0.25) {$\dots$};}

\def\ddotsbox at (#1,#2){\draw (#1,#2) rectangle (#1+0.5, #2 +0.5);  \node at (#1+0.23, #2 +0.33) {$\ddots$};}

\def\vdotsbox at (#1,#2){\draw (#1,#2) rectangle (#1+0.5, #2 +0.5);  \node at (#1+0.25, #2 +0.34) {$\vdots$};}

\def\bbox at (#1,#2){\filldraw[fill=gray,-] (#1,#2)--++(0.6,0)--++(-0.4,-0.3)--++(-0.6,0)--(#1,#2);}

\def\fbbox at (#1,#2){\filldraw[fill=gray,- ] (#1,#2)--++(0,0.6)--++(0.6,0)--++(0,-0.6)--(#1,#2);}

\def\fbboxh at (#1,#2){%
  \filldraw[fill=gray] (#1,#2) -- ++(0,0.6) -- ++(0.6,0) -- ++(0,-0.6) -- cycle;
  \draw[line width=2] (#1,#2+0.3) -- ++(0.6,0);
}

\def\rbbox at (#1,#2){\filldraw[fill=gray,-] (#1,#2)--++(0,0.6)--++(0.4,0.3)--++(0,-0.6)--(#1,#2); \draw[line width=2pt]
        (#1,#2+0.6) -- ++(0.4,0.3);}

\def\rbboxh at (#1,#2){%
  \filldraw[fill=gray]
    (#1,#2) -- ++(0,0.6) -- ++(0.4,0.3) -- ++(0,-0.6) -- cycle;
  \draw[line width=2pt]
    (#1,#2+0.3) -- ++(0.4,0.3);
}

\def\wbox at (#1,#2){\draw[fill=white,-] (#1,#2)--++(0.6,0)--++(-0.4,-0.3)--++(-0.6,0)--(#1,#2);}

\def\fwbox at (#1,#2){\draw[fill=white,-] (#1,#2)--++(0,0.6)--++(0.6,0)--++(0,-0.6)--(#1,#2);}

\def\fwboxh at (#1,#2){\draw[fill=white,-] (#1,#2)--++(0,0.6)--++(0.6,0)--++(0,-0.6)--(#1,#2); \draw[line width=2] (#1,#2+0.3) -- ++(0.6,0);}

\def\rwbox at (#1,#2){\draw[fill=white,-] (#1,#2)--++(0,0.6)--++(0.4,0.3)--++(0,-0.6)--(#1,#2);}

\def\rwboxh at (#1,#2){\draw[fill=white,-] (#1,#2)--++(0,0.6)--++(0.4,0.3)--++(0,-0.6)--(#1,#2);\draw[line width=2pt]
    (#1,#2+0.3) -- ++(0.4,0.3);}

\begin{document}

\author{Maico Freitas and Evgeny Mukhin}
\address{MF: Department of Mathematical Sciences,
Indiana University Indianapolis,
402 N. Blackford St., LD 270, 
Indianapolis, IN 46202, USA}
\email{mgfreita@iu.edu}
\address{EM: Department of Mathematical Sciences,
Indiana University Indianapolis,
402 N. Blackford St., LD 270, 
Indianapolis, IN 46202, USA}
\email{emukhin@iu.edu}

\title{Thin modules for Quantum Toroidal $\gl_n$}

\begin{abstract}  We develop tools for a systematic study of irreducible thin modules over the quantum toroidal algebra $\E_n$ of type $\gl_n$. As a result, we give a full classification and explicit construction of irreducible thin $\E_n$-modules whose restriction to every quantum affine $\sL_2$ is a sum of finite-dimensional modules.

We also apply our method to construct and study irreducible thin modules for the shifted quantum toroidal algebras $\E_{n,\bs s}$ of type $\gl_n$. 
\end{abstract}

\maketitle

Keywords:  quantum toroidal algebra, shifted quantum toroidal algebra, thin module, $q$-character, non-intersecting paths.

\tableofcontents

\section{Introduction}
In recent years, the representations of quantum toroidal algebras have become central in mathematics and physics. In almost all cases appearing in applications, these representations have distinguished bases described by various plane partitions with various boundary conditions, see \cite{ORV}, \cite{FJMM1}, \cite{K}, \cite{KN}, \cite{NW}. In these bases, the Cartan subalgebra is diagonalizable and has simple spectrum. Such a representation is called thin.

Thin modules have been studied in many cases. For Yangians of type A, see \cite{NT}; for quantum affine algebras of type A and B, see \cite{MY}, \cite{BrM}; for quantum toroidal algebras of type $D(2,1;\alpha), $ see \cite{DM}. For the quantum toroidal algebras of type $\gl_n,$ the Fock modules and MacMahon modules are in \cite{FJMM1}; the evaluation modules of weight modules of quantum affine $\sL_2$ are in \cite{JM}; the evaluation modules of admissible quantum affine $\sL_n$ are in \cite{BMRT}, and the thin modules over the quantum toroidal algebra of type $\gl_{m|n}$ are  in \cite{BM}.

Here, we offer a systematic treatment of thin modules over the quantum toroidal algebra $\E_n$ of type $\gl_n$. Our approach is based on the algorithm of combinatorial description of $q$-characters suggested in \cite{FM}. Given an $\ell$-weight (a set of eigenfunctions of the Cartan currents), the FM algorithm uses combinatorial expansions to produce other $\ell$-weights and their corresponding multiplicities, providing a candidate for the $q$-character of a module. There is no general theorem specifying when the FM algorithm produces the correct $q$-character. However, in the case of thin modules, the situation is more tractable since there are no multiplicities of $\ell$-weights. In this case, we expect that the FM algorithm gives a correct answer and, in many cases, we can provide a proof of that.
 On the other hand, the FM algorithm was originally formulated for finite-dimensional modules over the quantum affine algebras, so it requires a modification in the context of infinite-dimensional modules over the quantum toroidal algebras. In particular, in addition to the standard FM expansions in the direction of negative simple roots, we also perform expansions corresponding to positive simple roots. We also have to find non-standard proofs of theorems about injectivity of the $q$-character map.

We consider only $\ell$-weights which are $n$-tuples of rational functions of the spectral parameter. For example, all highest weight modules with finite-dimensional weight spaces are of that form. An important case is when all coordinates have the form $q^{\deg P} P(q^{-1}/z)(P(q/z))^{-1},$ where $P(x)$ is a rational function with $P(0)=1$. Such functions appear as $\ell$-weights of finite-dimensional $U_q\hat\sL_2$ modules. We call an $\E_n$-module locally finite if the restriction to each $U_q\hat\sL_2$ is a sum of finite-dimensional type one modules.  Then, all $\ell$-weights of a locally finite module have this form. This is the case we use to describe our techniques and where all thin modules are classified and described with full proofs.

Our strategy is as follows. 
Given an $\ell$-weight of a locally finite thin module, we picture it as a union of paths and formulate the expansions as  simple explicit changes in the path, see Figure \ref{path expansion}. 

The key observation is that, in all locally finite thin modules, the coordinates of the $\ell$-weights cannot have poles of order more than one, see Lemma \ref{thin exponents}. It turns out that only the non-intersecting paths do not produce double poles after expansions, see Theorem \ref{no intersections}. We call a union of non-intersecting paths a finite chain configuration. Each path contains two rays which are invariant under the expansions. The module is then parameterized by the collection of these rays. We explain that the combinatorics of the FM algorithm is easily understood in terms of finite chain configurations.

Then, we explicitly construct the corresponding irreducible $\E_n$-module by finding it inside a tensor product of Fock modules and vector modules, see Section \ref{construction chain}. In particular, it justifies the FM algorithm.

The matrix coefficients of the action of $\E_n$ generators $E_i(z)$ and $F_i(z)$ are given by $\delta$-functions multiplied by constants. The support of the delta function is at the pole of the $\ell$-weight, so, given the $q$-character, the only unknown parts of the action are the constants in front of the delta functions. These constants are called the  supplement.  We give a combinatorial proof that the supplement is unique up to renormalization of the basis of  $\ell$-weight vectors, see Section \ref{injectivity of q character}.  We also give explicit formulas for the supplement using the ideas of \cite{JM}, see Theorem \ref{plane partition module}. That gives an alternative proof of the existence of the modules corresponding to  finite chain configurations.

To show that there are no irreducible modules whose $q$-character is a sum of expansions of several finite chain configurations, we study and classify extensions between the simple modules we already constructed. The answer is given in terms of linking conditions, see Sections \ref{linked}, \ref{linear independece q}. Using that information, we prove Theorem \ref{U is reducible} ruling out other thin irreducible locally-finite modules. Another technical point is the proof of  uniqueness of an irreducible module with a given $q$-character. We prove it using a combinatorial property of the graph of the FM expansion which we call ``square complete", see Section \ref{square complete}. 

That completes the case of locally finite thin modules. The main result is Theorem \ref{irreducible is tensor product} which classifies and explicitly describes all irreducible thin locally finite  $\E_n$-modules up to isomorphism. Theorem \ref{prime thm} describes such prime modules.

\medskip

Then, we show that our method is applicable to all thin modules of $\E_n$ and all thin modules of the shifted quantum toroidal algebra $\E_{n,\bs s}$, see \cite{N}, \cite{FT} and Section \ref{shifted algebras}.

The only difference is that we must treat infinite chain configurations 
where some poles of an $\ell$-weight produce
infinitely many paths and some zeroes of an $\ell$-weight give truncation boxes. The truncation boxes limit the expansions in the FM algorithm. Given such a non-intersecting truncated chain configuration, we construct the corresponding thin module, see Theorem \ref{truncated chain theorem}. We conjecture that all totally prime modules (the shifted version of prime) are of this form, see Conjecture \ref{main conjecture}.  We expect that the proof of this conjecture can be obtained by similar methods as in the locally finite case. We plan to address this problem in a future publication.

In the general case, a parameterization of all irreducible thin modules (not only totally prime) is less clear. In particular, even in the totally prime case,  the same module can be described by different chain configurations, see e.g. Figure \ref{Infinite same monomial}. 

\medskip

In the case of thin modules, the study of representations becomes combinatorial. Moreover, our main claim is that the combinatorics of thin modules of $\E_{n,\bs s}$ presented via truncated chain configurations is rather simple. We made an effort to present sufficiently many details, examples, and pictures to make it accessible to a reader. 

The combinatorics of non-intersecting paths is motivated by \cite{MY} and it is well suited for our study. However, it is well known that it is equivalent to the language of plane partitions used in most papers on this subject. We explain the relation in detail in Section \ref{sec: plane partition}.

As an example, we have a family of thin locally finite modules whose weights are uniformly bounded.
Curiously, every irreducible thin finite-dimensional $U_q\hat\sL_n$ module appears as a submodule of an irreducible,  locally finite, bounded and  thin $\E_n$-module, see Section \ref{main results}.

 For simplicity of exposition, in this paper we treat the case $n\geq 3$. However, the cases of $n=2$ and $n=1$ are the same with only minor adjustments.

\medskip 

The paper is organized as follows. In Section \ref{locally finite section}, we recall the generalities of quantum toroidal algebras $\E_n$ of type $\gl_n$ focusing on thin locally finite modules.  We identify the $\ell$-weights of locally finite thin modules with certain monomials and discuss the expansions of these monomials following a version of the FM algorithm. 
In Section \ref{seed monomials section}, we view the monomials as collections of paths and study the combinatorics of the FM expansions.
In Section \ref{construction chain}, we construct the $\E_n$-modules, which we call chain modules, whose $\ell$-weights are parametrized by all possible expansions of a given finite chain configuration as submodules in tensor products of the vector, covector, Fock, and dual Fock representations. 
In Section \ref{injectivity of q character}, we prove that every irreducible and thin locally finite $\E_n$-module is isomorphic to a unique tensor product of shifted chain modules (up to a permutation of factors). 
In Section \ref{language of plane partitions}, we use the language of plane partitions to view the chain configurations as stacks of generalized Young diagrams. Using this language, we describe the supplements of the action in a chain module explicitly. In particular, this gives us an alternative proof of the existence of the chain modules. 
In Section \ref{sec: shifted}, we recall the  shifted quantum toroidal algebras. Then, we consider truncated chain configurations.  Using the notion of supplements once more, we define an action of the shifted quantum toroidal algebra on the vector space with a basis parametrized by the expansions of a truncated chain configuration. We conjecture that every irreducible, totally prime thin module is isomorphic to a shift of a module obtained from a truncated chain configuration. We finish by discussing various examples of thin modules for the shifted quantum toroidal algebra.

\section{Quantum Toroidal Algebras}  

\label{locally finite section}

\subsection{Type A quantum toroidal algebras} 
Fix $n\in \mathbb Z_{\geq 3}$.  

Let $I = \{1,\dots, n-1\}$ and $\hat I = \{0\}\cup I $ be the sets of Dynkin nodes of types $A_{n-1}$ and A$^{(1)}_{n-1}$, respectively. We often refer to the elements of  $\hat I$ as \textit{colors}. We order the  sets $I$ and $\hat I$ as subsets of $\Z$.

For $a,b\in \mathbb Z$, we write $a\equiv b$ whenever $a,b$ represent the same integer modulo $n$.  When $x\in\mathbb Z$ is interpreted as a color, we understand it modulo $n$.

Fix parameters $q,d\in \mathbb C^\times$, and let
$$q_1 = d/q, \quad q_2=q^2, \quad q_3 = 1/(dq).$$
Then, $q_1q_2q_3=1.$  We assume that $q_1,q_2$ are generic, in the sense that $q_1^aq_2^b=1$ for $a,b \in \mathbb Z$ if and only if $a=b=0$. In particular, $q_1,q_2$ are not roots of unity. 

We fix the square root of $q_2$ to be $q$, that is, $q_2^{1/2}=q$. We fix a logarithm of $q$ and set  $q^c=e^{c\log q}$ for all $c\in \mathbb C.$

The\textit{ quantum toroidal $\mathfrak{gl}_{n}$ algebra}\footnote{In general, the quantum toroidal algebra has an extra central element $C$, see \cite{FJMM1}, \cite{FJMM2}. Our algebra $\E_n$ is a quotient of that  by the relation $C=1$. We consider only modules where
$C$ acts as $1$.} $\mathcal E_n = \mathcal E_{n}(q_1,q_2,q_3)$ is the unital associative algebra  with generators $E_{i,k}, F_{i,k}, H_{i,l}, K_i^{\pm 1}$, where $i\in \hat I, k\in \mathbb Z, l \in \mathbb Z\setminus\{0\},$ and relations as follows. Introduce the following generating series in a formal variable $z$:
$$E_i(z) = \sum_{k\in \mathbb Z} E_{i,k}\,z^{-k}, \quad F_i(z) = \sum_{k\in \mathbb Z} F_{i,k}\,z^{-k}, \quad K^{\pm}_i(z) = K_i^{\pm1}\, \exp \big(\pm(q-q^{-1})\sum_{l=1}^{\infty}H_{i,\pm l}z^{\mp l}\big), \quad i \in \hat I. $$
For $i,j \in \hat I$, define further $g_{i,j}(z,w)$ in formal variables $z,w$ by: 
\begin{align*}
    g_{i,j}(z,w) = \begin{cases} z-q_2w, \quad &\mbox{if} \quad i=j, \\
z-q_1w, \quad &\mbox{if} \quad i = j-1,\\
d(z-q_3w), \quad &\mbox{if} \quad i = j+1,\\
z-w, \quad &\mbox{otherwise}.
\end{cases}
\end{align*}
Let $\delta(z) = \displaystyle \sum_{k\in \mathbb Z} z^k$ be the formal delta distribution. 
The algebra $\E_n$ has the following defining relations: 

\begin{alignat}{2} 
&  \qquad \qquad K_iK_i^{-1}=K_i^{-1}K_i = 1,\label{equations 2.1} \\
 &[K_i^{\pm}(z),K_j^{\pm}(w)] = [K_i^+(z),K_j^{-}(w)] = 0,   \label{relations 2.2}
\\
&g_{i,j}(z,w)K_i^{\pm }(z)E_j(w) +g_{j,i}(w,z)E_j(w)K_i^{\pm}(z) = 0, \\
&g_{j,i}(w,z)K_i^{\pm }(z)F_j(w) +g_{i,j}(z,w)F_j(w)K_i^{\pm}(z) = 0, \\
   [E_i&(z),F_j(w)] =  \frac{\delta_{i,j}}{q-q^{-1}} \big(\delta(w/z)K^+_i(w)-\delta(z/w)K_i^-(z)\big), \label{EFK}  \\
  &g_{i,j}(z,w)E_i(z)E_j(w)+g_{j,i}(w,z)E_j(w)E_i(z)= 0, \quad \mbox{if}\hspace{0.5cm} i=j, j\pm 1, \label{relations 2.5} \\
   &g_{j,i}(w,z)F_i(z)F_j(w)+g_{i,j}(z,w)F_j(w)F_i(z)= 0, \quad \mbox{ if} \hspace{0.5cm}  i=j, j\pm 1,\label{relations 2.6} \\
   & \hspace{1.5cm} [E_i(z),E_j(w)] = [F_i(z), F_j(w)] = 0, \qquad \hspace{0.4cm} \mbox{if} \hspace{0.5cm} i\neq j, j\pm 1, \label{relations 2.7} \\
   E_i(z_1)E_i(&z_2)E_{i\pm 1}(w)-(q+q^{-1})E_i(z_1)E_{i\pm 1}(w)E_i(z_2)+E_{i\pm 1}(w)E_i(z_1)E_i(z_2)+(z_1\leftrightarrow z_2) = 0,\\
   F_i(z_1)F_i&(z_2)F_{i\pm 1}(w)-(q+q^{-1})F_i(z_1)F_{i\pm 1}(w)F_i(z_2)+F_{i\pm 1}(w)F_i(z_1)F_i(z_2)+(z_1\leftrightarrow z_2) = 0, \label{relations 2.9}
\end{alignat}
for all $i,j\in\hat I$.

The algebra $\mathcal E_n$ has an invertible central element $K$ given by \begin{equation*}K = \label{K central} K_0K_1\cdots K_{n-1}.\end{equation*}

Let $\mathcal E_n^+, \mathcal E_n^0, \mathcal E_n^-$ be the subalgebras of $\mathcal E_n$ generated by the coefficients of the generating series of $E_i(z), K_i^{\pm}(z), F_i(z), i\in \hat I,$ respectively. 
    We have an isomorphism of vector spaces $$\mathcal E_n \cong \mathcal E_n^-\otimes \mathcal E_n^0\otimes \mathcal E_n^+.$$

For each $i \in \hat I$, the subalgebra of $\mathcal E_n$ generated by the coefficients of $E_i(z), K_i^{\pm}(z), F_i(z)$ is isomorphic to the quantum affine algebra $U_q\,\hat{ \sL}_2$ at level $0$ in the Drinfeld new realization. We denote this subalgebra by $\hat{\mathcal U}_i$. The subalgebra of $\hat{\mathcal U}_i$ generated by $K_i, E_{i,0},F_{i,0}$ is isomorphic to the quantum algebra $U_q\,\sL_2$. We denote this subalgebra by $\mathcal U_i$. 

For each $i \in \hat I$, the subalgebra of $\mathcal E_n$ generated by $\hat{\mathcal U}_j, j \in \hat I\setminus\{i\},$ is isomorphic to the quantum affine algebra $U_q\,\hat{ \sL}_n$ at level $0$ in the Drinfeld new realization. The subalgebra of $\mathcal E_n$ generated by $\mathcal U_j, j \in \hat I\setminus\{i\},$ is isomorphic to the quantum algebra $U_q\,\sL_n$.

We introduce the formal  \textit{coproduct}\footnote{Our choice of coproduct coincides with \cite{BM}  and is opposite to \cite{FJMM1}, \cite{FJMM2}.} $\Delta$ on $\mathcal E_n$:
\begin{align}     \label{coproduct}
\begin{split}
&\Delta E_i(z) = 1 \otimes E_i(z) + E_i(z) \otimes K_i^-(z), \\
&\Delta F_i(z) = F_i(z)\otimes 1 + K^{+}_i(z)\otimes F_i(z),\\
&\Delta K^{\pm}_i(z) = K^{\pm}_i(z) \otimes K_i^{\pm}(z), 
\end{split}
\end{align}
for all  $i\in \hat I$.

The above formulas do not define a coproduct of Hopf algebras in the usual sense, as the right-hand sides involve infinite sums.  We use the coproduct $\Delta$ only in situations where these sums are convergent.

Let $\iota:  \mathcal E_n(q_1,q_2,q_3)\rightarrow \mathcal E_n(q_3,q_2,q_1)$ be the isomorphism of algebras given by
\begin{equation}  \label{transpose isomorphism}
    E_i(z) \mapsto E_{n-i}(z), \quad F_i(z) \mapsto F_{n-i}(z), \quad K^{\pm}_i(z) \mapsto K^{\pm}_{n-i}(z), \quad i\in \hat I.
\end{equation}

Let $\upsilon: \mathcal E_n(q_1,q_2,q_3) \rightarrow \mathcal E_n(q_1^{-1},q_2^{-1},q_3^{-1})$ be the isomorphism of algebras given by
$$E_i(z) \mapsto F_{n-i}(z), \quad F_i(z) \mapsto E_{n-i}(z), \quad K^{\pm}_i(z)\mapsto K^{\pm}_{n-i}(z), \quad i \in \hat I.$$ 

For $a\in \mathbb C^\times, $ let $\sigma_a:\mathcal E_n\rightarrow \mathcal E_n$ be the shift of spectral parameter automorphism given by
$$E_i(z)\mapsto E_i(z/a),\quad F_i(z)\mapsto F_i(z/a), \quad K_i^{\pm}(z)\mapsto K_i^{\pm}(z/a), \quad i\in \hat I. $$

Let $\tau: \mathcal E_n\rightarrow \mathcal E_n$ be the  shift of color automorphism given by
$$ E_i(z) \mapsto E_{i+1}(z), \quad F_i(z)\mapsto F_{i+1}(z), \quad K_i^{\pm}(z)\mapsto K^{\pm}_{i+1}(z), \quad i \in \hat I.$$

\subsection{Preliminaries  on \texorpdfstring{$\mathcal E_{n}$}{}-modules} Let $V$ be an $\mathcal E_n$-module. 

We say that $V$ has \textit{level $c\in \mathbb C^\times$}  if the central element $K\in \mathcal E_n$ acts on $V$ as the scalar $c$. 

Let $\varphi= (\varphi_i)_{i\in \hat I}$ be an $\hat I$-tuple of complex numbers.  The tuple $\varphi$ is called a \textit{weight} of $V$ if there exists a nonzero  vector $v\in V$ such that $K_iv = q^{\varphi_i}v$ for all $i\in \hat I.$ Such a vector $v$ is called a \textit{weight vector} of weight $\varphi.$

Let $ \bm \phi^{\pm} = (\phi^{\pm}_i(z))_{i\in \hat I}$ be an $\hat I$-tuple such that $\phi_i^{\pm}(z) \in \mathbb C[[z^{\mp 1}]]$. The tuple $\bm \phi^\pm$ is called an \textit{$\ell$-weight} of $V$ if there exists  a nonzero vector $v\in V$ such that
$$K^{\pm}_i(z)v = \phi^{\pm}_i(z)v, \quad i \in \hat I.$$
Such a vector  $v$ is called an \textit{$\ell$-weight vector} of $\ell$-weight $\bm \phi^{\pm}$.

Notice that if $\bm \phi^{\pm}$ is an $\ell$-weight, then  $\phi_i^+(\infty)\phi_i^-(0)=1$ for $i\in \hat I$, cf. \eqref{equations 2.1}.

A nonzero vector $v\in V$ is called  a \textit{generalized $\ell$-weight vector} of $\ell$-weight $\bm \phi^{\pm}$ if
$$(K^{\pm}_i(z)- \phi^{\pm}_i(z))^Nv=0, \ N\gg0, \quad i \in \hat I.$$

The associated vector subspace of generalized eigenvectors of an $\ell$-weight $\bm \phi^{\pm}$ is called the \textit{$\ell$-weight space} of $\ell$-weight $\bm \phi^\pm$.

An $\ell$-weight vector $v$ is called a \textit{highest $\ell$-weight vector} if $E_{i,k}v = 0$, $i\in\hat I, k\in\Z$. Similarly, an $\ell$-weight vector $v$ is called  a \textit{lowest $\ell$-weight vector} if $F_{i,k}v = 0,$ $i\in\hat I, k\in\Z$. An $\mathcal E_n$-module is called a \textit{highest} (resp. \textit{lowest}) \textit{$\ell$-weight module} if it is generated by a highest (resp. lowest) $\ell$-weight vector.

An $\mathcal E_n$-module $V$ is called \textit{$\ell$-weighted} if it admits a basis of generalized $\ell$-weight vectors with finite-dimensional $\ell$-weight spaces. An $\ell$-weighted module is called \textit{thin} if all  $\ell$-weight spaces are one-dimensional.

We note that any submodule or quotient module of an $\ell$-weighted (resp. thin) module is again $\ell$-weighted (resp. thin).

For $ i\in \hat I, a\in \mathbb C^\times,$ let $\bm i_a$ be the $\hat I$-tuple of rational functions\footnote{The  $i$-th rational function of $\bm i_a$ is written in the variable $1/z$ in accordance with the convention that $K_i^+(z)$ is a generating series in $1/z.$} given by
\begin{equation} \label{fundamental l weights}
    \bm i_a = \big (1,\dots,1,\underbrace{q\frac{1-q^{-1}a/z}{1-qa/z}}_{\mbox{ \tiny \textit{i}-th place}}, 1, \dots, 1\big ).
\end{equation}

Let $\hat{\mathcal P}$ be the abelian group generated by the $\hat I$-tuples  $\bm i_a, i\in \hat I, a\in \mathbb C^\times,$ where the product is given by coordinate-wise multiplication of rational functions. Then, the elements of $\hat{\mathcal P}$ 
are monomials in the variables $\bm i_a$ of the form: \begin{equation} \label{monomial}
    \bm m = \prod_{i\in \hat I,\, a \in\mathbb C^\times}\bm i_{a}^{k_{i,a}(\bm m)}, \quad k_{i,a} (\bm m)\in \mathbb Z,\end{equation} where only finitely many $k_{i,a}(\bm m)$ are nonzero. 

We call the collection of integers $(k_{i,a}(\bm m))_{i\in \hat I, a\in \mathbb C^\times}$ as above the \textit{exponents of $\bm m.$} Whenever $\bm m$ is clear from the context, we simply write $k_{i,a}$ for the exponents $k_{i,a}(\bm m)$.

Let $\bm m\in \hat{\mathcal P}$ be as in (\ref{monomial}). For $i\in \hat I,$ we say that $\bm m$ is \textit{$i$-dominant} (resp. \textit{$i$-antidominant}) if $k_{i,a}\in \mathbb Z_{\geq 0}$ (resp. $k_{i,a}\in \mathbb Z_{\leq 0}$) for all $a\in \mathbb C^\times.$ If $\bm m$ is $i$-dominant (resp. $i$-antidominant) for all $i\in \hat I$, we simply say that $\bm m$ is \textit{dominant} (resp. \textit{antidominant}).

\iffalse

Let $\mbox{wt}:\mathcal P\rightarrow P$ be the weight homomorphism defined on generators by
$$\bm i_a \mapsto  \mbox{wt}(\bm i_a) = \omega_i, \quad i \in \hat I, a\in \mathbb C^\times.$$

We call an $\mathcal E_n$-module $V$ \textit{rational} if it is $\ell$-weighted and all $\ell$-weights are of the form $ \bm m^\pm = (m^{\pm}_i(z))_{i\in \hat I}$ for some $\bm m\in \hat{\mathcal P}$.
\fi

\begin{dfn} \label{locally finite}

We call an $\mathcal E_n$-module $V$  \textit{locally finite}  if for every $i\in \hat I,$ the restriction of $V$ to $\hat{\mathcal U}_i$ decomposes as a sum of countably many finite-dimensional type one $\hat{\mathcal U}_i$-modules.
\end{dfn}

 Equivalently,  an $\mathcal E_n$-module V is locally finite if and only if $V$ has a countable basis and, for every $v\in V$ and $i\in\hat I$, there exists a finite-dimensional type one $\hat{\mathcal U}_i$-submodule of $V$ containing $v$.

Note that in Definition \ref{locally finite} we do not require $V$ to be a direct sum of finite-dimensional $\hat{\mathcal U}_i$-modules, although all modules we study here shall have this property.

For $\bm m\in \hat{\mathcal P}$, denote by $m_i^+(z)$ (resp. $m_i^-(z)$) the power series expansion of the $i$-th rational function of $\bm m$ at $z=\infty$ (resp. at $z=0)$. We set $\bm m^{\pm}  = (m^{\pm}_i(z))_{i\in \hat I}$. 

We call a rational function $f(z)$ balanced if $f(0)$, $f(\infty)$ are well defined and $f(0)f(\infty)=1$. All coordinates $m_i(z)$ of $\bs m\in\hat{\mathcal{P}}$ are balanced rational functions. Therefore, the constant terms of $m^{+}_i(z)$ and  $m^{-}_i(z)$ are inverse to each other.

\begin{lem} If $V$ is a locally finite $\mathcal E_n$-module, then the $\ell$-weights of $V$ are of the form $ \bm m^\pm = (m^{\pm}_i(z))_{i\in \hat I}$ for some $\bm m\in \hat{\mathcal P}$. In particular, for each $i\in \hat I$, the eigenvalues of $K_i$ are in $q^{\mathbb Z}.$
\end{lem}
\begin{proof}
The lemma follows from the corresponding $U_q\,\hat{\sL}_2$ statement, see for example \cite{CP}.
\end{proof}
All modules in this paper are going to be $\ell$-weighted. Except in Section \ref{sec: shifted}, they are also locally finite. In that case,  each $\ell$-weight space is finite-dimensional and corresponds to some $\bs m \in \hat{\mathcal P}$. We denote the corresponding $\ell$-weight space of $V$ simply by $V_{\bs m}.$

Let $\mathbb Z\langle \hat{\mathcal P}\rangle$ be  the abelian group whose elements are countable sums of the form $\sum_{k=0}^\infty c_k \bm m_k$, where $c_k\in\Z$ and $\bm m_k\in\hat{\mathcal P}$ are distinct monomials.

The \textit{$q$-character}  of a locally finite $\ell$-weighted $\mathcal E_n$-module $V$ is the formal sum: 
$$ \chi_q(V) = \sum_{\bm m\in \hat{\mathcal P}} (\dim V_{ \bm m}) \, \bm m\in \mathbb Z\langle \hat{\mathcal P}\rangle.$$

For $S\subset \hat{\mathcal P} $,  we write  $S\subset \chi_q(V)$ if every $\bs m\in S$ is a summand in $\chi_q(V)$.   

For $c\in \C^\times$, let $\sigma_c: \mathbb Z\langle\hat{\mathcal P}\rangle\rightarrow \mathbb Z\langle \hat{\mathcal  P} \rangle $ denote the group homomorphism mapping $\bm i_a \mapsto \bm i_{ca}$ for all $i\in\hat I$, $a\in \mathbb C^\times.$ Then, for any locally finite $\ell$-weighted module $V$, we have $$\chi_q(\sigma_c(V)) = \sigma_c(\chi_q(V)).$$

For $i\in \hat I,$ let $\rho_i: \hat{\mathcal P}\rightarrow  \hat{\mathcal  P} $ denote the group homomorphism obtained by setting $$\bm j_a = 1 \quad \mbox{for all} \quad  j\neq i, a\in \mathbb C^\times.$$

\begin{lem} \label{thin exponents}
If $V$ is a locally finite thin $\mathcal E_n$-module and $\bs m\in\hat{\mathcal P}$ is an $\ell$-weight of $V$, then the exponents of $\bm m$ satisfy $k_{i,a}\in\{-1,0,1\}$ for all $i\in\hat I$, $a\in \C^\times$. In addition, if $k_{i,a}=1$, then   $k_{i,aq^2}\neq -1$. 

In particular, all poles of the rational functions $m_i(z)$ are simple.
\end{lem}
\begin{proof}

Fix $i\in \hat I$ and regard $V$ as a module for the algebra $\hat{\mathcal U}_i\cong U_q\,\hat{\sL}_2$. Then,  $\rho_i(\bm m) = \prod_{a\in \mathbb C^\times} \bm i_a^{k_{i,a}},$ where $k_{i,a}\neq0$ only for finitely many $a\in \mathbb C^\times.$ 

An $\mathcal E_n$-module $V$ is called tame if the Cartan currents $K_i^\pm(z)$ are diagonalizable in $V$. Thin modules are clearly tame.
Restricting to $\hat{ \mathcal U}_i$, taking a submodule, or quotient of a locally finite tame module results in a locally finite and tame module. 

Therefore, $\rho_i(\bm m)$ is a monomial in the $q$-character of a tame irreducible finite-dimensional $\U_q\hat{\sL}_2$-module.

A tame irreducible finite-dimensional $\U_q\hat{\sL}_2$-module must also be thin. In particular, we have $|k_{i,a}|\leq 1$. Moreover, there are no double poles, so $k_{i,a}=1$  implies  $k_{i,aq^2}\neq -1$,
see, for example, Lemmas 3.5 and 3.6 in \cite{BrM}. 
\end{proof}

We call  a \textit{thin monomial} any monomial $\bm m \in \hat{\mathcal P}$ such that $k_{i,a} \in \{-1,0,1\}$  and such that if $k_{i,a}=1$, then   $k_{i,aq^2}\neq -1$  for every $i\in \hat I, a \in \mathbb C^\times$.

\begin{dfn}\label{definition}
A locally finite, $\ell$-weighted $\mathcal E_n$-module $V$ is called   \textit{bounded} if, for all $i\in \hat I$, 
     the set of eigenvalues of $K_i$ in $V$ is a finite subset of $q^{\mathbb Z}$. 
 \end{dfn}
We remark that it is likely that the property of being locally  finite for an $\ell$-weighted module
 may follow from the bounding condition on the eigenvalues of $K_i$ in
the Definition \ref{definition}. We do not pursue this question in the current paper.

\medskip 

We give straightforward reformulations of  Definition \ref{definition} in the following lemma.

\begin{lem}  \label{bounded lemma}
Let $V$ be a locally finite, $\ell$-weighted $\mathcal E_n$-module. The following are equivalent:
    
    \begin{enumerate}
        \item The module $V$ is bounded.
        \item  There exists $N\in\Z_{>0}$ such that if $\bm m \in \hat{\mathcal P}$ is an $\ell$-weight of $V$ with exponents $(k_{i,a})_{i\in \hat I, a\in \mathbb C^\times}$, then $|\sum_ak_{i,a}|\le N$ for all $i\in \hat I.$
        \item There exists $N\in\Z_{>0}$ such that if $L\subset V$ is an irreducible  finite-dimensional $\mathcal U_i$-module for some $i\in\hat I$, then $\dim L\leq N$. 

    \end{enumerate}

\end{lem}
\begin{proof}
For each $i\in \hat I$,  the corresponding eigenvalue of $K_i$ is given by $q^{\sum_a k_{i,a}}$, which shows that  the first and second parts are equivalent.

 If $L$ is an irreducible $U_q\sL_2$-module and $\dim L=r$, then the eigenvalues of $K$ in $L$ are of the form $q^k$, where $k=-r+1,-r+3,\dots,r-1$. It follows that the last part is equivalent to the first.
\end{proof}

\subsection{The FM \texorpdfstring{$q$-character algorithm}{}}  \label{strings section}
We describe a version of the FM algorithm adapted to the setting  of locally finite thin modules, see Theorem \ref{expansion is contained} below.

For $i\in \hat I, a\in \mathbb C^\times,$ define the $\ell$-root $\bm A_{i,a}\in \hat{\mathcal P}$ by setting
$$\bm A_{i,a} = \bm i_{aq^{-1}}\bm i_{aq}\bm{(i+1)}^{-1}_{aqq_3}\bm{(i-1)}^{-1}_{aqq_1}. $$

Let $\bm m\in \hat{\mathcal P}$ be a thin monomial with exponents $(k_{i,a})_{i\in \hat I, a\in \mathbb C^\times}$. For $i\in \hat I,$ let $\Sigma_i^{\pm}(\bm m) \subset \mathbb C^\times $ be the set of simple poles of the rational function
$$ \prod_{a : k_{i,a} = \pm 1} \bigg(q\frac{1-aq^{-1}/z}{1-aq/z}\bigg)^{\pm 1}$$ corresponding to the $i$-th component of the monomial $\prod_{a: k_{i,a}=\pm 1}\bs i^{\pm 1}_a$.

The set of poles $\Sigma^{\pm}_{i}(\bm m)$ can be described in terms of $i$-strings as follows. For $a\in \mathbb C^\times, r\in \mathbb Z_{\geq 0}$, we call  a monomial $\bm i_a\bm i_{aq_2}\cdots \bm i_{aq_2^r}$ a \textit{dominant $i$-string of $\bm m$} if $k_{i,aq_2^s} = 1$ for $s=0,\dots, r$ and $k_{i,aq_2^{-1}}, k_{i,aq_2^{r+1}} \neq  1.$ 
The $i$-th rational function of  $\bm i_a\bm i_{aq_2}\cdots \bm i_{aq_2^r}$ has the unique simple pole   $aqq_2^r$. Then, the set $\Sigma_i^+(\bm m)$ is the set of simple poles of all dominant $i$-strings of $\bs m$. 

Similarly, for $a\in \mathbb C^\times, r\in \mathbb Z_{\geq 0}$, we call a monomial $\bm i_a^{-1}\bm i_{aq_2}^{-1}\cdots \bm i_{aq_2^r}^{-1}$ an \textit{antidominant $i$-string of $\bm m$} if $k_{i,aq_2^s}= -1$ for $s=0,\dots, r$ and $k_{i,aq_2^{-1}}, k_{i,aq_2^{r+1}} \neq -1.$ The $i$-th rational function of $\bm i_a^{-1}\bm i_{aq_2}^{-1}\cdots \bm i_{aq_2^r}^{-1}$ has the unique simple pole $aq^{-1}$. Then, the set $\Sigma_i^-(\bm m)$ is the set of simple poles of all antidominant $i$-strings of $\bs m$.

From  $\Sigma_i^{\pm}(\bm m)$, we define the collection of  monomials $e_i(\bm m), f_i(\bm m) \subset \hat{\mathcal P}$  given by
\begin{align*}
&f_i(\bm m) = \{\bm m\cdot \bm A_{i,a}^{-1}: a\in \Sigma_i^+(\bm m)\},\\ &e_i(\bm m) = \{\bm m\cdot \bm A_{i,a}: a\in \Sigma_i^-(\bm m)\}. 
\end{align*}

We call the collection of monomials $f_i(\bm m)$ and $e_i(\bm m)$ the \textit{positive} and \textit{negative $i$-expansions of $\bm m$}, respectively.

Let $\bs m, \bs m'$ be thin monomials. It is easy to see that $\bs m'\in e_i (\bs m)$ if and only if  $\bs m\in f_i (\bs m')$.

\begin{example} \label{thin example}
Let $\hat I = \{0,1,2,3\}$, corresponding to $n=4$. Consider the thin monomial: \begin{equation*}  \bm m = \bm 0_{q^{-1}}^{-1}\bm 0_{qq_2^2}^{-1}\bm 1_{q^{-1}q_3}\bm 1_{qq_3} \bm 1_{qq_1^{-1}q_2} \bm 3_{qq_3^{-1}}^{-1}. \label{monomial vertical wiggle}
\end{equation*}

The dominant $1$-strings of $\bm m$ are $\bm 1_{q^{-1}q_3}\bm 1_{qq_3}$ and $\bm 1_{qq_1^{-1}q_2}$. The antidominant $0$-strings of $
\bm m$ are $\bm 0_{q^{-1}}^{-1}$ and  $\bm 0_{qq_2^2}^{-1}$. The antidominant $3$-string of $
\bm m $ is $\bm 3_{qq_3^{-1}}^{-1}$. The corresponding nonempty sets of poles are: 
\begin{align*}
 \Sigma_1^+(\bm m) = \{q_1^{-1}, q_1^{-1}q_2^2\}, \quad 
\Sigma^-_0(\bm m) = \{q_2^{-1}, q_2^2\}, \quad \Sigma_3^-(\bm m)  = \{q_3^{-1}\}.
\end{align*}
The corresponding nonempty expansions of $\bs m$ are:
\begin{align*}
f_1(\bs m)= \{ \bm m \bm A_{1,q_1^{-1}}^{-1},\bm m \bm A_{1,q_1^{-1}q_2^2}^{-1}  \}, \quad  
e_0(\bs m)=  \{ \bm m \bm A_{0, q_2^{-1}}, \bm m \bm A_{0,q_2^2}   \}, \quad  e_3(\bm m) = \{\bm m\bm A_{3,q_3^{-1}}\},
\end{align*}
where the monomials are given as follows: \begin{align*}
& \bm m \bm A_{1,q_1^{-1}}^{-1} = \bm 0_{q^{-1}}^{-1}\bm 0_{q}\bm 0_{qq_2^2}^{-1}\bm 1_{q^{-1}q_3}\bm 1_{qq_1^{-1}}^{-1} \bm 1_{qq_1^{-1}q_2} \bm 2_{qq_1^{-1}q_3} \bm 3_{qq_3^{-1}}^{-1}, \\
&\bm m \bm A_{1,q_1^{-1}q_2^2}^{-1} = \bm 0_{q^{-1}}^{-1}\bm 1_{q^{-1}q_3}\bm 1_{qq_3} \bm 1_{qq_1^{-1}q_2^2}^{-1} \bm 2_{qq_1^{-2}q_2}   \bm 3_{qq_3^{-1}}^{-1}, \\
& \bm m \bm A_{0,q_2^{-1}} =     \bm 0_{q^{-1}q_2^{-1}}\bm 0_{qq_2^2}^{-1}\bm 1_{qq_3} \bm 1_{qq_1^{-1}q_2} \bm 3_{qq_3^{-1}}^{-1} \bm 3_{q^{-1}q_1}^{-1},  \\
& \bm m \bm A_{0,q_2^{2}} =  \bm 0_{q^{-1}}^{-1}\bm 0_{qq_2}\bm 1_{q^{-1}q_3}\bm 1_{qq_3}  \bm 3_{qq_3^{-1}}^{-1} \bm 3_{qq_2q_3^{-1}}^{-1} , \\
&\bm m \bm A_{3,q_3^{-1}} = \bm 0_{q^{-1}}^{-1}\bm 0_{q}^{-1}\bm 0_{qq_2^2}^{-1} \bm 1_{q^{-1}q_3}\bm 1_{qq_3} \bm 1_{qq_1^{-1}q_2} \bm 2_{qq_1q_3^{-1}}^{-1} \bm 3_{qq_1}.
\end{align*}
\end{example}

\begin{dfn}    \label{good monomial}
We call a thin monomial $\bs m$ a seed monomial if there exists a set $S\subset \hat{\mathcal P}$ such that
\begin{itemize}
\item $\bs m\in S$,
\item $\bm t$ is thin for all $\bm t\in S$, 
\item  for all $\bs t\in S,  i\in \hat I$, we have $e_i(\bs t)\subset S$, $f_i(\bs t)\subset S$.
\end{itemize}
\end{dfn}

Given a seed monomial $\bs m,$ we denote by $S(\bm m)$ the smallest set $S$ with respect to inclusion satisfying the conditions in Definition \ref{good monomial}. In other words, $\mbox{S}(\bm m)$ is the smallest set of thin monomials containing $\bm m$ which is closed under positive and negative expansions in all directions.
Note that all monomials $\bs t\in S(\bm m)$ are also seed monomials and  $S(\bm t ) = S(\bm m)$.

We call a seed monomial $\bs m$ a {\it bounded seed monomial} if 
$S(\bs m)$ is bounded in the sense of Lemma \ref{bounded lemma} (ii). In other words, $\bs m$ is a bounded seed monomial if there exists $N\in\Z_{\geq 0}$ such that,  for all $\bm t\in S(\bm m)$   and for all $i\in \hat I$,  we have $|\sum_a k_{i,a}(\bm t)|\leq N$.

For a seed monomial $\bm m$, we call the formal sum
$$\mbox{FM}(\bm m) = \sum_{\bm t\in S(\bm m)} \bm t\in\Z\langle \hat{\mathcal P}\rangle$$
  the FM \textit{expansion} of $\bm m$.

\begin{thm} \label{expansion is contained} Let $V$ be a locally finite, irreducible and thin $\mathcal E_n$-module. Let $\bm m\in  \chi_q(V).$ Then, $\bm m$ is a seed monomial and \normalfont{{\textrm{FM}}}$(\bm m)\subseteq \chi_q(V).$
\end{thm}

\begin{proof}
Let $\bs m\in \chi_q(V)$ be a thin monomial and let $i\in\hat I$. We claim that $e_i(\bs m)\subset \chi_q(V)$ and $f_i(\bs m)\subset \chi_q(V)$. 
We argue for the case of $e_i(\bm m)$, the case of $f_i(\bm m)$ is completely parallel.

For $\bs t\in \chi_q(V)$, let $v_{\bs t}\in V$ be a nonzero eigenvector for the Cartan currents $K_i^\pm(z)$ with eigenfunctions $\bs t_i^\pm(z)$. Since $V$ is locally finite and thin, the set $\{v_{\bs t}\}_{\bs t\in\chi_q(V)}$ is a basis of $V$.

Since $V$ is locally finite, we have $\hat{\mathcal U}_i$ finite-dimensional submodules $U\subset  W$ of  $V$ such that $W/U$ is irreducible, $v_{\bm m}\in W$ and $v_{\bm m}\not \in U$. For any $a_0\in \Sigma_i^-(\bs m)$, we know that the restriction  $\rho_i(\bm m  \bm  A_{i,a_0})$ belongs to the  $q$-character of $W/U$. Therefore, $E_i(z) v_{\bm m}$ is non-zero.

By Proposition 3.3 in \cite{MY}, $E_i(z) v_{\bs m}$ belongs to the span of vectors of the form $v_{\bs m \bm A_{i,a}}$, $a\in \C^\times$. Since the restrictions  $\rho_i(\bs A_{i,a})$ are distinct for distinct $a\in \C^\times$, it follows that the matrix coefficient of $E_i(z)$ between $v_{\bs m}$ and $v_{\bs m  \bm  A_{i,a_0}}$ must be non-zero. Thus, $\bs m  \bm  A_{i,a_0}\in \chi_q(V)$.

The claim immediately implies the theorem.
\end{proof}
We will prove that if $\bs m$ is a seed monomial, then there exists a locally finite, irreducible,  thin $\mathcal E_n$-module $V$ such that \textrm{FM}$(\bm m)=\chi_q(V)$, see Theorem \ref{H and V modules}.

\medskip

We call a monomial $\bs m\in \hat{\mathcal P}$  of {\it integer type} if, for any nonzero exponent $k_{i,a}$ of $\bs m$, the spectral parameter $a\in \mathbb C^\times$ is of the form $$a=q q_1^{-x}q_3^{-y},$$ where $x,y\in \mathbb Z$ and $i\equiv x-y \in \hat I.$

Equivalently,  a thin monomial $\bm m$ is of integer type  if, for   $i\in \hat I$,  every pole $p\in \Sigma_i^{\pm}(\bm m)$ is of the form $p = q_1^{-x}q_3^{-y}$, where $x,y \in \mathbb Z$ and $x-y \equiv i.$

Similarly to \cite{CM}, \cite{EM}, \cite{FM}, in Lemma \ref{split lemma} we claim that we can restrict our attention to seed monomials of integer type.

We say that the spectral parameters $a,b\in\C^\times$ are \textit{$n$-independent} over $q_1,q_2$ if $a/b\not\in q_1^{n\mathbb Z}q_2^{\mathbb Z}=q_3^{n\mathbb Z}q_2^{\mathbb Z}$.

\begin{lem}\label{split lemma}
Let $\bs m\in\hat{\mathcal P}$ be a non-trivial monomial. There exist $s\in \mathbb Z_{>0},$ $a_j\in \mathbb C^\times$, and unique non-trivial monomials $\bs m_j\in \hat{\mathcal P}$, $j=1,\dots,s,$ such that
\begin{itemize}
\item  $\bs m=\prod_{j=1}^s\bs m_j,$
\item Each $\sigma_{a_j^{-1}}(\bs m_j)$ is a  monomial of integer type,
\item For all $i\neq j$, the spectral parameters $a_i, a_j$ are $n$-independent over $q_1,q_2$.
\end{itemize}
\end{lem} 
\begin{proof}
For every variable $\bm i_b$ appearing in $\bm m$, choose $a\in \mathbb C^\times$ in the coset $\frac{bq_1^i}{q} H, $ where $H = q_1^{n\mathbb Z}q_2^{\mathbb Z}.$ Then, $b/a = qq_1^{-x}q_3^{-y}$ with $x,y \in \mathbb Z, x-y\equiv i.$ This allows us to break $\bm m$ into disjoint factors $\bm m_a$, where each $\bm m_a$ is a product of variables $\bm i_b,$ where $bq_1^{i}/q$ belongs to the same coset $a H$.
\end{proof} 

\begin{prop}    \label{independent monomials}
Let $\bs m=\prod_{j=1}^s\bs m_j$ be as in Lemma \ref{split lemma}. Assume that $\bm m_j$ are seed monomials. Then, for every $j$ and $\bs t\in S(\bs m_j)$,  $\sigma_{a_j^{-1}}(\bs t)$ is a seed monomial of integer type. Moreover, $$\normalfont{\mbox{FM}}(\bm m)=\prod_{j=1}^s \normalfont{\mbox{FM}}(\bs m_j).$$
\end{prop}
\begin{proof}
Since every pole $p\in \Sigma_i^{\pm}(\sigma_{a_j^{-1}}(\bm m_j))$ is of the form $p=q_1^{-x}q_3^{-y}$, where  $x,y\in \mathbb Z$ and $x-y\equiv i$, then  $\bm A^{\pm 1}_{i,p}$ is  also  a monomial of integer type. The first part of the proposition follows.

Note that if $\sigma_{a}(\bs t_1)$ and  $\sigma_{b}(\bs t_2)$ are of integer type and $a,b$ are $n$-independent over $q_1,q_2$, then the monomials 
$\bs t_1$ and $\bs t_2$ do not share variables, that is, $k_{i,c}(\bs t_1)\neq 0$ implies  $k_{i,c}(\bs t_2)= 0$. The last part of the proposition follows.
\end{proof}

Under the assumptions of Proposition \ref{independent monomials},  in the case of bounded modules, we will prove that  the tensor product of locally finite, irreducible, thin $\mc E_n$-modules corresponding to $\bs m_i$ is locally finite, irreducible, and thin.

\section{Seed monomials}    \label{seed monomials section}
We classify the seed monomials of integer type, see Theorem \ref{no intersections}  below. We also classify the bounded seed monomials of integer type, see Theorem \ref{bounded seed}.

\subsection{Pictures of monomials} \label{monomial pic}
We start by developing a visual interpretation of monomials and their expansions.

Consider the integer plane $\mathbb Z\times \mathbb Z.$ 
At each point $(x,y)$, we place a unit box $\square$ centered at $(x,y)$ and write  $\square=(x,y).$ For $\square = (x,y)$ and $r,s\in \mathbb Z,$ we write $\square \pm  (r,s)$ for the box at $(x\pm r, y\pm s).$

To each box $\square = (x,y),$ we assign a {\it color}  $i\in \hat I$ given by $i\equiv x-y$. Notice that boxes aligned diagonally along the line $y=x$ share the same color.

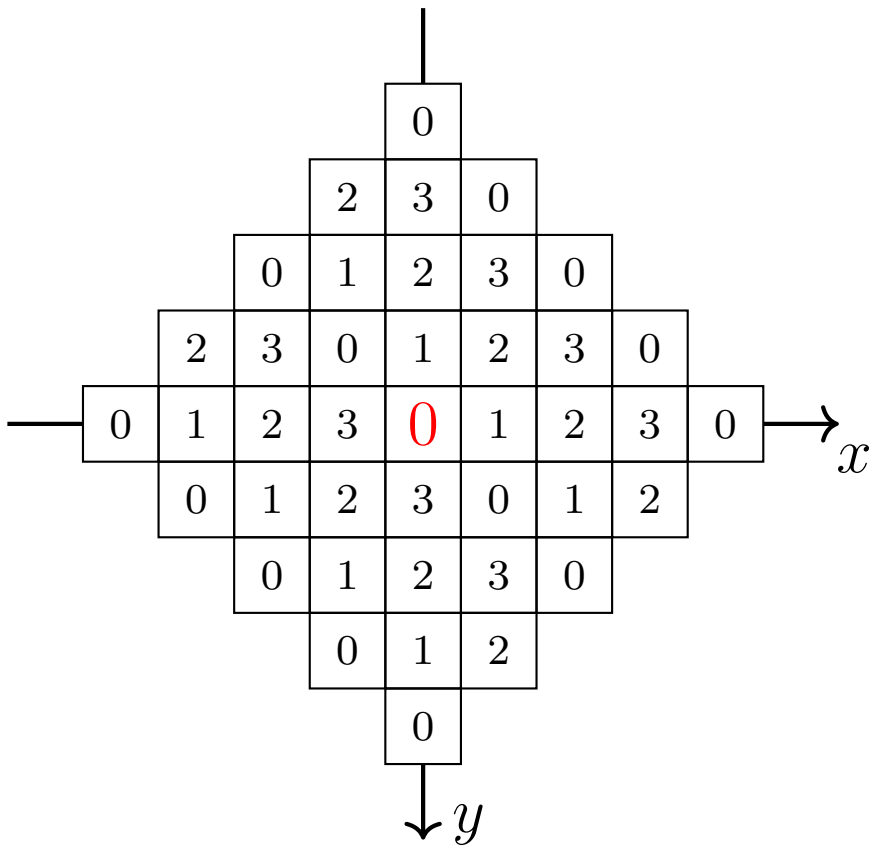
\begin{figure}[ht!]  
\begin{tikzpicture}

\rzerobox at (0,0);
\onebox at (0.5,0);
\twobox at (1,0);
\threebox at (1.5,0);
\zerobox at (2,0);

\threebox at (0,-0.5);
\zerobox at (0.5,-0.5);
\onebox at (1,-0.5);
\twobox at (1.5,-0.5);

\twobox at (0,-1);
\threebox at (0.5,-1);
\zerobox at (1,-1);

\onebox at (0,-1.5);
\twobox at (0.5,-1.5);

\zerobox at (0,-2);

\threebox at (-0.5,0);
\twobox at (-1,0);
\onebox at (-1.5,0);
\zerobox at (-2,0);

\twobox at (-0.5,-0.5);
\onebox at (-1,-0.5);
\zerobox at (-1.5,-0.5);

\onebox at (-0.5,-1);
\zerobox at (-1,-1);

\zerobox at (-0.5,-1.5);

\onebox at (0,0.5);
\twobox at (0.5,0.5);
\threebox at (1,0.5);
\zerobox at (1.5,0.5);

\twobox at (0,1);
\threebox at (0.5,1);
\zerobox at (1,1);

\threebox at (0,1.5);
\zerobox at (0.5,1.5);

\zerobox at (0,2);

\zerobox at (-0.5,0.5);
\threebox at (-1,0.5);
\twobox at (-1.5,0.5);

\onebox at (-0.5,1);
\zerobox at (-1,1);

\twobox at (-0.5,1.5);

\node at (3.1, 0.02) {$x$};
\node at (0.55, -2.4) {$y$};
\draw[->, thick] (0.25,-2) -- (0.25,-2.5);
\draw[ thick] (0.25,2.5) -- (0.25,3);
\draw[->, thick] (2.5,0.25) -- (3,0.25);
\draw[thick] (-2,0.25) -- (-2.5,0.25);

\end{tikzpicture}  
\caption{A finite portion of $\mathbb Z \times \mathbb Z$ with colored boxes for $n=4.$}\label{Figure grid}
\end{figure}

In Figure \ref{Figure grid}, we show  a finite portion of the plane $\mathbb Z\times \mathbb Z$ with boxes colored by the set $\hat I = \{0,1,2,3\}$ corresponding to $n=4.$ The box at $(0,0)$  has color $0$ shown in red. In all our pictures, we shall always orient the coordinate axes as in Figure \ref{Figure grid}.

To each box $\square = (x,y),$ we assign a \textit{content}   $c(\square)\in q_1^{\mathbb Z}q_3^{\mathbb Z}$ given by 
\begin{equation}c(\square) = q_1^{-x}q_3^{-y}. \label{content 2d}
\end{equation}
In particular, we have $c(\square+(r,r)) = c(\square)q_2^r$, for $r\in \mathbb Z.$ Since the parameters $q_1,q_2$ are generic, boxes are distinguished by their content, i.e.,  $c(\square_1)=c(\square_2)$  if and only if $\square_1= \square_2$.

Let $\mathbb Z +\frac{1}{2}$ be the set of half-integers. The set $(\mathbb Z+\frac{1}{2})\times (\mathbb Z+\frac{1}{2})$ is naturally identified with the set of corners of all unit boxes in the integer plane $\mathbb Z \times \mathbb Z.$  We also define the content of a corner $(x+ \frac{1}{2},y+\frac{1}{2}) \in (\mathbb Z+\frac{1}{2})^2$ by setting $$c\big(x+\textstyle\frac{1}{2},y+\textstyle \frac{1}{2}\big) = qq_1^{-x}q_3^{-y}.$$

We identify monomials of integer type with finitely supported maps $(\mathbb Z+\frac{1}{2})\times (\mathbb Z+\frac{1}{2}) \rightarrow \mathbb Z$ as follows. A monomial $\bm m \in \hat{\mathcal P}$  of integer type with exponents $(k_{i,a})_{i\in \hat I, a\in \mathbb C^\times}$ corresponds to a finitely supported map $\kappa_{\bm m}: (\mathbb Z+\frac{1}{2})\times (\mathbb Z+\frac{1}{2}) \rightarrow \mathbb Z$ given by
$$\textstyle \kappa_{\bm m}(x+\frac{1}{2},y+\frac{1}{2}) =  k_{i,qq_1^{-x}q_3^{-y}},$$ where  $ i\equiv  x-y \in \hat I$.

%Conversely, a finitely supported map  $\kappa:  (\mathbb Z+\frac{1}{2})\times (\mathbb Z+\frac{1}{2}) \rightarrow \mathbb Z$ is identified with a monomial $\bm m_{\kappa}\in \hat{\mathcal P}$ of integer type  given by
%$$\bm m_{\kappa} = \prod_{(x,y)\in\Z\times\Z}\bm i_{qq_1^{-x}q_3^{-y}}^{\kappa({x+\frac{1}{2},y+\frac{1}{2})}}, $$ where   $i\equiv x-y \in \hat I$.

Every thin monomial $\bm m$ of integer type  determines a  finitely supported map $\kappa_{\bs m}$ with image in the set $\{-1,0,1\}$.
  We picture such a monomial by placing a bullet ``\,\textbullet\," and a cross ``$\bm \times$" at every corner   mapped by  $\kappa_{\bm m}$ to $1$ and $-1$, respectively.

\begin{figure}[H]  
\begin{tikzpicture}

\onebox at (0,0.5);
\twobox at (0.5,0.5);
\threebox at (1,0.5);
\zerobox at (1.5,0.5);
\rzerobox at (0,0);
\onebox at (0.5,0);
\twobox at (1,0);
\threebox at (1.5,0);
\threebox at (0,-0.5);
\zerobox at (0.5,-0.5);
\onebox at (1,-0.5);
\twobox at (1.5,-0.5);
\twobox at (0,-1);
\threebox at (0.5,-1);
\zerobox at (1,-1);
\onebox at (1.5,-1);

\node at (0,1) {\textbullet};
\node at (0,0.5) {$\bm \times$};
\node at (0.5,0.5) {\textbullet};
\node at (0.5,-0.5) {$\bm \times$};
\node at (1.5,-0.5) {\textbullet};
\node at (1.5,-1) {$\bm \times$};

\end{tikzpicture}
\caption{A visualization of a thin monomial of integer type.} \label{thimonomialviz}
\end{figure}
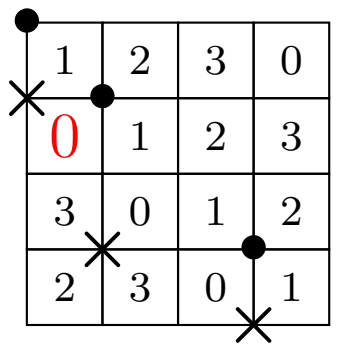

The thin monomial of integer type 
$\bm m = \bm 0_{q^{-1}}^{-1}\bm 0_{qq_2^2}^{-1}\bm 1_{q^{-1}q_3}\bm 1_{qq_3} \bm 1_{qq_1^{-1}q_2} \bm 3_{qq_3^{-1}}^{-1}$ in Example \ref{thin example} is pictured in Figure \ref{thimonomialviz}.

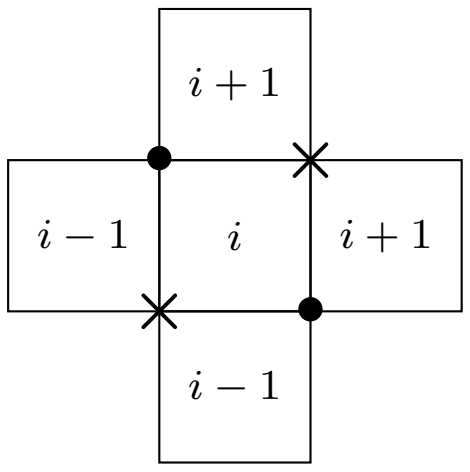
\begin{figure}[H]  
\begin{tikzpicture}

\ibox at (0,0);
\inextbox at (1,0);
\inextbox at (0,1);
\iprevbox at (-1,0);
\iprevbox at (0,-1);

\node at (0,1) {\textbullet};
\node at (1,0) {\textbullet};
\node at (1,1) {$\bm \times$};
\node at (0,0) {$\bm \times$};

\end{tikzpicture}
\caption{A visualization of the $\ell$-root $\bm A_{i,c(\square)}$.}
\label{lrootpic}
\end{figure}

Let $\square = (x,y)$ be a box of color $i\in \hat I.$ Then, the $\ell$-root
$$\bm A_{i,c(\square)} = \bm i_{q^{-1}c(\square)}\bm i_{qc(\square)}\bm{(i+1)}^{-1}_{qq_3c(\square)}\bm{(i-1)}^{-1}_{qq_1c(\square)}$$
is pictured by placing bullets and crosses at the four corners of $\square$ as in  Figure \ref{lrootpic}. The picture of $\bm A_{i,c(\square)}^{-1}$ is obtained from the picture of  $\bm A_{i,c(\square)}$ by exchanging crosses and bullets.

We often say that the $\ell$-root $\bm A_{i,c(\square)}$ is centered at the box $\square = (x,y).$

\begin{figure}[ht!]  
\begin{tikzpicture}

\squarebox at (0,0);
\squarebox at (0.5,-0.5);
\squarebox at (1, -1);
\diagonalbox at (1.5,-1.5);

\node at (0.5,0) {\textbullet};
\node at (1,-0.5) {\textbullet};
\node at (1.5,-1) {\textbullet};

\draw[->, thick] (2.5,-0.5) -- (4,-0.5);
\node at (3.2,-0.1) {$\bm A^{-1}_{i,c(\boxxslash)}$};

\squarebox at (4.5,0);
\squarebox at (5,-0.5);
\squarebox at (5.5, -1);
\diagonalbox at (6,-1.5);

\node at (5,0) {\textbullet};
\node at (5.5,-0.5) {\textbullet};
\node at (6,-1.5) {\textbullet};
\node at (6.5,-1) {\textbullet};
\node at (6.5,-1.5) {$\bm \times$};

\node at (0.25,0.28) {$i$};
\node at (4.75,0.28) {$i$};

\end{tikzpicture}
\caption{Expansion at a bullet. }

\label{dominantexpansion}
\end{figure}
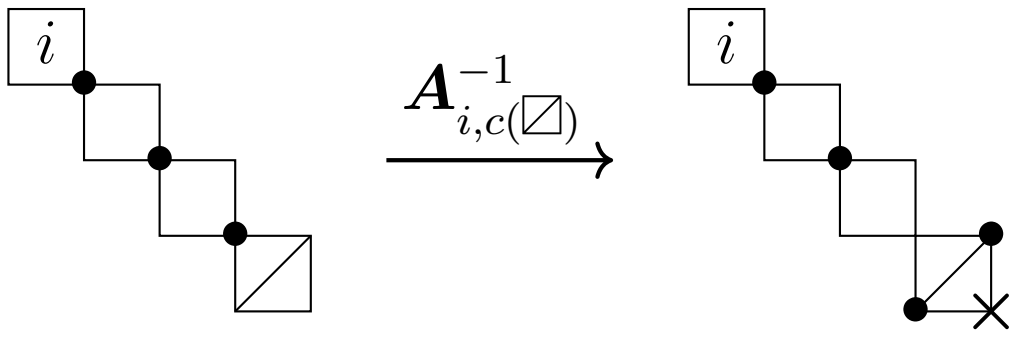

A dominant $i$-string  is illustrated in Figure \ref{dominantexpansion}. In the figure, the content $c(\boxxslash)$ of the box $\boxxslash = (x,y)$ gives the only pole coming from the dominant $i$-string. In this case, we say that the bullet at $(x-\frac{1}{2},y-\frac{1}{2})$ is \textit{expandable}. If such an $i$-string is part of a thin monomial $\bs m$, then we call the monomial $\bm m A^{-1}_{i,c(\boxxslash)}\in f_i(\bm m)$ \textit{the expansion of $\bm m$ at the bullet}   $(x-\frac{1}{2},y-\frac{1}{2}).$

\begin{figure}[H]  
\begin{tikzpicture}

\diagonalbox at (-0.5,0.5);
\squarebox at (0,0);
\squarebox at (0.5,-0.5);
\squarebox at (1, -1);

\node at (0,0.5) {$\bm \times$};
\node at (0.5,0) {$\bm \times$};
\node at (1,-0.5) {$\bm \times$};
\node at (1.5,-1) {$\bm \times$};

\draw[->, thick] (2,0) -- (3.5,0);
\node at (2.7,0.3) {$\bm A_{i,c(\boxxslash)}$};

\diagonalbox at (4,0.5);
\squarebox at (4.5,0);
\squarebox at (5,-0.5);
\squarebox at (5.5, -1);

\node at (4, 0.5) {$\bm \times$};
\node at (4.5,1) {$\bm \times$};
\node at (4,1) {\textbullet};
\node at (5,0) {$\bm \times$};
\node at (5.5,-0.5) {$\bm \times$};
\node at (6,-1) {$\bm \times$};

\node at (5.75,-0.75) {$i$};
\node at (1.25,-0.75) {$i$};

\end{tikzpicture}
\caption{Expansion at a cross. }
\label{antidominantexp}
\end{figure}
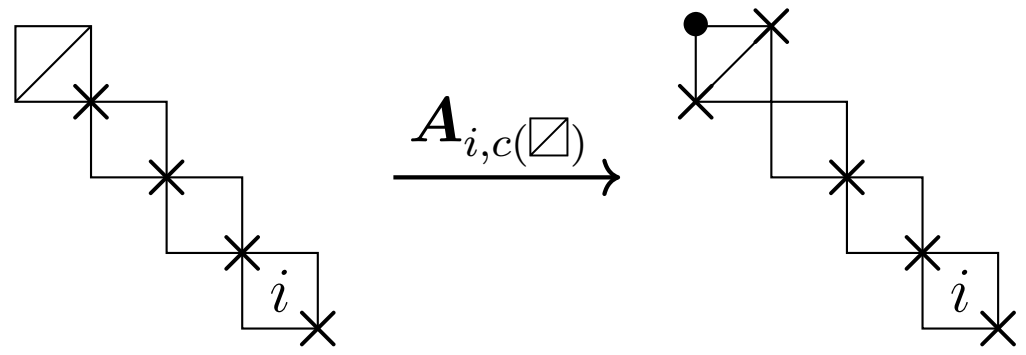

Analogously,  an antidominant $i$-string is illustrated in Figure \ref{antidominantexp}. The box $\boxxslash = (x,y)$ is such that $c(\boxxslash)$ gives the only pole coming from the antidominant $i$-string. We say that the cross at $(x+\frac{1}{2}, y+\frac{1}{2})$ is \textit{expandable}. If such an antidominant $i$-string is part of a thin monomial $\bs m$, then we call the monomial $\bm m A_{i,c(\boxxslash)}\in e_i(\bm m)$ \textit{the expansion of $\bm m$ at the cross  $(x+\frac{1}{2},y+\frac{1}{2}).$}

\subsection{The paths of a thin monomial} In this section, we view a thin monomial as a union of segments or rays.

Let $\bs m$ be a non-trivial thin monomial of integer type. From each bullet in the picture of $\bm m,$ we draw a horizontal segment extending to the right and a vertical segment extending downward until we reach  the first cross in that direction. If no such cross exists, we  draw a ray extending to the right or downward, respectively. Formally, given a bullet at the corner $(u,v),$ let $r_1, r_2\in \mathbb Z_{>0}$  be the smallest integers such that $(u+r_1,v)$ and $(u,v+r_2)$ are crosses. If such $r_1$ or $r_2$ does not exist, set it to be $+\infty$. Then, we draw straight segments connecting $(u,v)$ to $(u+r_1,v)$ and $(u,v+r_2)$.

Similarly, from each cross in the picture of $\bs m$, we draw a horizontal segment extending to the left and a vertical segment extending upward until we reach the first bullet in that direction. If no such bullet exists, we draw a  ray extending to the left or upward, respectively.  Formally, given a cross at $(u,v),$ let $r_1,r_2\in \mathbb Z_{>0}$ be the smallest integers such that $(u-r_1,v)$ and $(u,v-r_2)$ are bullets, otherwise set $r_1=+\infty$ or $r_2=+\infty.$ We draw straight segments connecting $(u,v)$ to $(u-r_1,v)$ and $(u,v-r_2)$.

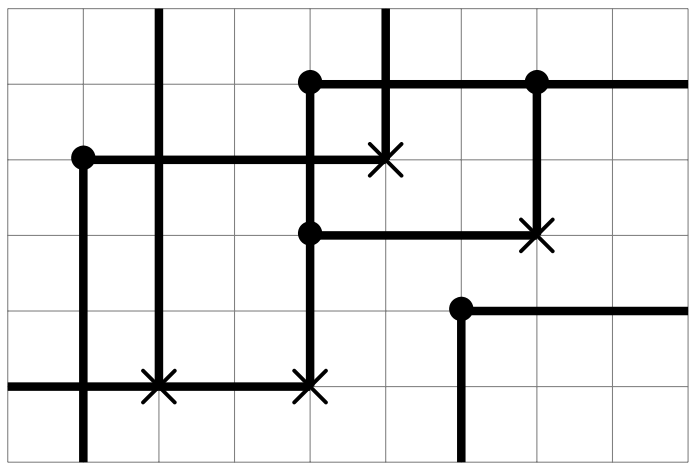
\begin{figure}[H]
\begin{tikzpicture}
 
 \draw[step=0.5cm, gray, very thin] (-2,-0.5) grid (2.5,2.5);

 \node at (-1.5,1.5) {\textbullet};
 \node at (-1,0) {$\bm \times$};
 \node at (0, 0) {$\bm \times$};
 \node at (0,1) {\textbullet};
 \node at (0,2) {\textbullet};
 \node at (1.5,1) {$\bm \times$};
 \node at (1.5,2) {\textbullet};
 \node at (0.5,1.5) {$\bm \times$};
 \node at (1,0.5) {\textbullet};
 \draw[ultra thick] (-1,0)--(-2,0);
 \draw[ultra thick] (-1,0)--(-1,2.5);
 \draw[ultra thick] (-1.5,1.5)--(0.5,1.5);
 \draw[ultra thick] (-1.5,1.5)--(-1.5,-0.5);
 \draw[ultra thick] (-1,0)--(0,0);
 \draw[ultra thick] (0,0)--(0,2);
 \draw[ultra thick] (0,1)--(1.5,1);
 \draw[ultra thick] (1.5,1)--(1.5,2);
 \draw[ultra thick] (0,2)--(2.5,2);
 \draw[ultra thick] (0.5,1.5)--(0.5,2.5);
 \draw[ultra thick] (1,0.5)--(2.5,0.5);
 \draw[ultra thick] (1,0.5)--(1,-0.5);

\end{tikzpicture}
\caption{The paths of a thin monomial.}
\label{paths thin}
\end{figure}

We call the resulting picture the \textit{paths of $\bs m$}. Figure \ref{paths thin} gives an example of the paths of a thin monomial. We notice that the resulting paths do not depend on $n$. Therefore, in the path pictures we often omit the color of boxes.

In the following lemma, we give necessary conditions on the paths for a monomial to be a seed monomial.

\begin{lem}   \label{subpaths lemma}
 Let $\bm m$ be a seed monomial of integer type. Then, the paths of $\bs m$ do not contain any of the following eight types of subpictures:
$$
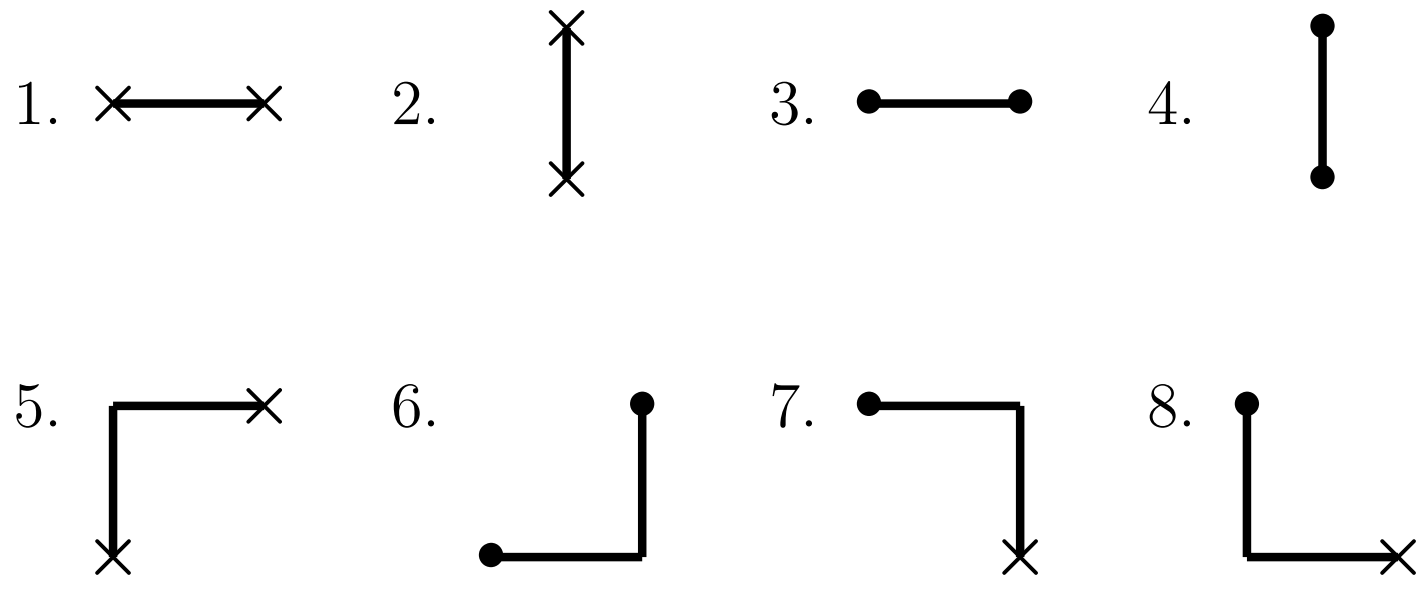
\begin{tikzpicture}
 \node at (-0.5,0) {$1.$};
 \node at (0,0) {$\bm \times$};
 \node at (1,0) {$\bm \times $};
 \draw[ultra thick] (0,0)--(1,0);
 \node at (2,0) {$2.$};
 \node at (3,0.5) {$\bm \times$};
 \node at (3,-0.5) {$\bm \times$};
 \draw[ultra thick] (3,0.5)--(3,-0.5);
 \node at (4.5,0) {$3.$};
 \node at (5, 0) {\normalfont \textbullet};
 \node at (6,0) {\normalfont \textbullet};
 \draw[ultra thick] (5,0)--(6,0);
  \node at (7,0) {$4.$};
  \node at (8, 0.5) {\normalfont \textbullet};
  \node at (8,-0.5) {\normalfont \textbullet};
  \draw[ultra thick] (8,0.5)--(8,-0.5);
  \node at (-0.5, -2) {$5.$};
  \node at (1, -2) {$\bm \times$};
  \node at (0,-3) {$\bm \times$};
  \draw[ultra thick] (0,-2)--(0,-3);
  \draw[ultra thick] (0,-2)--(1,-2);
  \node at (2,-2) {$6.$}; 
  \node at (2.5,-3) {\normalfont \textbullet};
  \node at (3.5,-2) {\normalfont \textbullet};
  \draw[ultra thick] (2.5, -3)--(3.5,-3);
  \draw[ultra thick] (3.5,-3)--(3.5,-2);
  \node at (4.5,-2) {$7.$};
   \node at (5, -2) {\normalfont \textbullet};
   \node at (6,-3) {$\bm \times$};
   \draw[ultra thick] (5,-2)--(6,-2);
   \draw[ultra thick] (6,-2)--(6,-3);
   \node at (7, -2) {$8.$};
   \node at (7.5,-2) {\normalfont \textbullet};
   \node at (8.5,-3) {$\bm\times$};
   \draw[ultra thick] (7.5,-2)--(7.5,-3);
   \draw[ultra thick] (7.5,-3)--(8.5,-3);
\end{tikzpicture}$$ 
In these pictures, except for the endpoints, we assume that there are no bullets and no crosses on the drawn segments.
\end{lem}

\begin{proof}

If the paths of $\bm m$ contain any of the pictures shown in the statement of the lemma, we perform successive expansions of $\bm m$ until we   arrive at a monomial with at least one exponent equal to $+2$ or $-2.$ Since such a monomial is not thin, $\bm m $ is not a seed monomial.

In picture 1, let $(u,v)$ and $(u',v), u<u',$ denote the coordinates of the crosses. After possibly performing successive expansions of antidominant strings of $\bm m$, we may assume that there are no crosses at the corners 
$$(u'-1,v-1), (u'-2,v-1),\dots,(u,v-1).$$

We first expand  the rightmost cross at $(u',v)$. That creates crosses at $(u'-1,v)$ and $(u', v-1)$. If the resulting monomial is not thin, then $\bs m$ is not a seed monomial. Otherwise, we expand the cross at $(u'-1,v)$. That  creates a cross at $(u'-2,v)$. Continuing in this way,  we either obtain a non-thin monomial or a monomial with crosses at $(u+1,v)$ and $(u,v)$. Finally, expanding the cross at $(u+1,v)$ yields a  monomial such that the exponent of the corner at $(u,v)$ is $-2.$ 

The cases of pictures 2, 3, and 4 are treated similarly. We successively expand crosses from bottom to top in picture 2, bullets from left to right in picture 3 and from top to bottom in picture 4. 

In the remaining pictures, we also expand along the edges. For example,
let $(u,v')$ and $(u',v), u<u', v<v'$ be the crosses in picture 5. We may assume that there are no crosses at  $(u'-1,v-1), (u'-2,v-1),\dots, (u,v-1)$ and at $(u-1,v'-1), (u-1,v'-2),\dots, (u-1,v)$.

We successively expand the crosses at $(u',v), (u'-1,v), \dots, (u+1,v)$, then the crosses located at $(u,v'), (u,v'-1),\dots, (u,v+1)$. The resulting  monomial has exponent -2 for the corner at $(u,v).$
\end{proof}

 We call the paths of $\bs m$ {\it non-intersecting } if  the following conditions hold: any two different horizontal segments or rays have no common points; any two different vertical segments or rays have no common points;  a vertical segment (or ray) can only intersect a horizontal segment (or ray)  at a cross or a bullet.

\begin{figure}[H]
\begin{tikzpicture}
 
 \draw[step=0.5cm, gray, very thin] (-2,-0.5) grid (2.5,2.5);

\node at (0.25,0.75) {\textcolor{red}{0}};
\node at (-1,0) {\textbullet};
 \node at (0.5,0) {$\bm \times$};
 \node at (0.5,0.5) {\textbullet};
 \node at (1.5,0.5) {$\bm \times$};
\node at (1.5,1) {\textbullet};
\node at (-1.5,1) {\textbullet};
\node at (0.5, 1) {$\bm \times$};
\node at (0.5,1.5) {\textbullet};
\node at (-0.5, 1.5) {$\bm \times$};
\node at (-0.5,2) {\textbullet};
\node at (1.5,0) {\textbullet};
\node at (-1,2) {$\bm \times$};
\node at (1.5,2) {$\bm \times$};
\node at (2,1.5) {$\bm \times$};

\draw[ultra thick] (1.5,0)--(1.5,-0.5);
\draw[ultra thick] (1.5,0)--(2.5,0);
\draw[ultra thick] (-1,0)--(-1,-0.5);
\draw[ultra thick] (-1,0)--(0.5,0);
\draw[ultra thick] (0.5,0)--(0.5,0.5);
\draw[ultra thick] (0.5,0.5)--(1.5,0.5);
\draw[ultra thick] (1.5,0.5)--(1.5,1);
\draw[ultra thick] (1.5,1)--(2.5,1);
\draw[ultra thick] (-1.5,1)--(-1.5,-0.5);
\draw[ultra thick] (-1.5,1)--(0.5,1);
\draw[ultra thick] (0.5,1)--(0.5,1.5);
\draw[ultra thick] (0.5,1.5)--(2,1.5);
\draw[ultra thick] (-0.5,1.5)--(-2,1.5);
\draw[ultra thick] (-0.5,1.5)--(-0.5,2);
\draw[ultra thick] (-0.5,2)--(1.5,2);
\draw[ultra thick] (-1,2)--(-2,2);
\draw[ultra thick] (-1,2)--(-1,2.5);
\draw[ultra thick] (1.5,2)--(1.5,2.5);
\draw[ultra thick] (2,1.5)--(2,2.5);

\end{tikzpicture}
\caption{A thin monomial with non-intersecting paths.}
\label{non intersecting}
\end{figure}
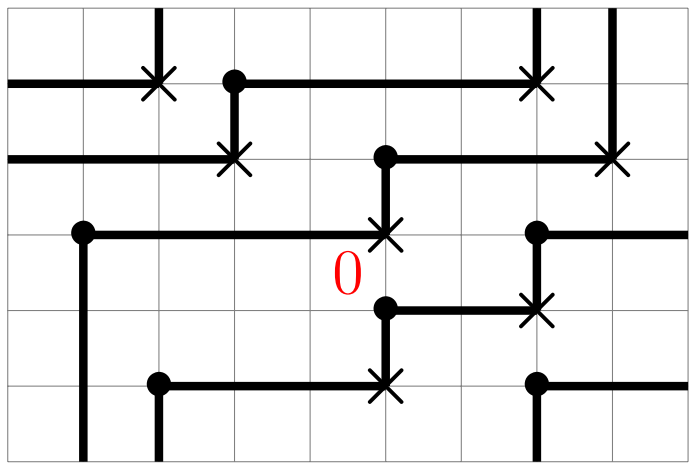

An example of a monomial $\bs m$ with non-intersecting paths is given in Figure \ref{non intersecting}. The paths of the monomial in Figure \ref{thimonomialviz} are also non-intersecting.

The following lemma is straightforward.
\begin{lem} \label{non intersecting is equivalent}
A thin monomial $\bs m$ of integer type does not contain any of the eight subpictures of Lemma \ref{subpaths lemma} if and only if the paths of $\bs m$ are non-intersecting.
\end{lem}

\begin{proof}
If two horizontal segments or rays of $\bm m$ have common points, then we have a subpicture of type 1 or 3. Similarly, if two vertical segments or rays of $\bm m$ have common points, then we  have a subpicture of type 2 or 4.

If a vertical segment or ray intersects a horizontal segment or ray at a point which is neither a bullet nor a cross, then we have a subpicture of type 5, 6, 7 or 8.
\end{proof}

\medskip

\subsection{Classification of seed monomials and their FM expansions}
We start by showing that thin monomials of integer type  with non-intersecting paths are exactly the seed monomials.

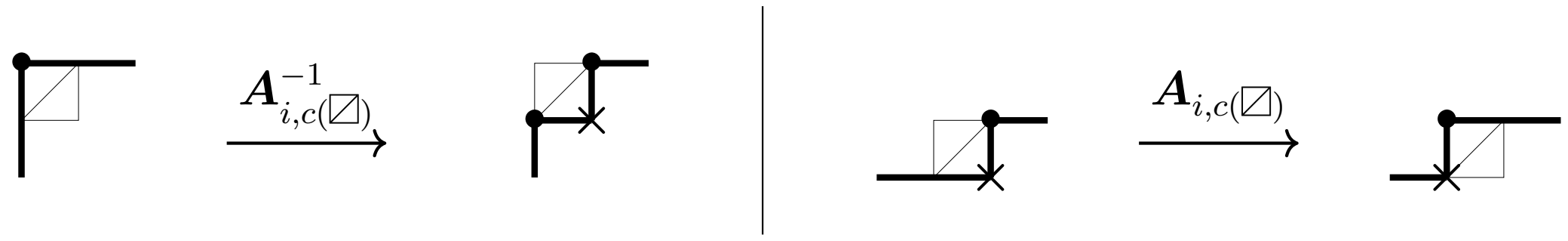
\begin{figure}[H]
\begin{tikzpicture}
 
 \thindiagonalbox at (0,0);
 \node at (0,0.5) {\textbullet};
 \draw[ultra thick] (0,0.5)--(0,-0.5);
\draw[ultra thick] (0,0.5)--(1,0.5);

 \node at (2.5,0.2) {$\bm A^{-1}_{i,c(\boxxslash)}$};
\draw[->, thick] (1.8,-0.2) -- (3.2,-0.2);

 \thindiagonalbox at (4.5,0);
 \node at (4.5,0) {\textbullet};
 \node at (5,0) {$\bm \times$};
 \node at (5,0.5) {\textbullet};
\draw[ultra thick] (4.5,0)--(4.5,-0.5);
\draw[ultra thick] (4.5,0)--(5,0);
\draw[ultra thick] (5,0)--(5,0.5);
\draw[ultra thick] (5,0.5)--(5.5,0.5);

\draw (6.5,-1)--(6.5,1);

\thindiagonalbox at (8,-0.5);
 \node at (8.5,-0.5) {$\bm \times$};
 \node at (8.5,0) {\textbullet};
 \draw[ultra thick] (7.5,-0.5)--(8.5,-0.5);
 \draw[ultra thick ] (8.5,-0.5)--(8.5,0);
\draw[ultra thick] (8.5,0)--(9,0);
 \node at (10.5,0.2) {$\bm A_{i,c(\boxxslash)}$};
\draw[->, thick] (9.8,-0.2) -- (11.2,-0.2);
\thindiagonalbox at (12.5,-0.5);
\node at (12.5,-0.5) {$\bm \times$};
\node at (12.5,0) {\textbullet};

\draw[ultra thick] (12,-0.5)--(12.5,-0.5);
\draw[ultra thick] (12.5,-0.5)--(12.5,0);
\draw[ultra thick] (12.5,0)--(13.5,0);

\end{tikzpicture}
\caption{Examples of local modifications of a path by expansion.}
\label{path expansion}
\end{figure}

The expansions of thin monomials with non-intersecting paths are illustrated in Figure \ref{path expansion}. The expansion at a bullet or a cross consists of a local modification of the path around the box 
$\boxxslash$, obtained by removing the two sides of the box that are currently part of the path and replacing them with the other two sides. Then, bullets and crosses  change accordingly.

\begin{thm}  \label{no intersections}
Let $\bm m$ be a thin monomial of integer type. Then, $\bm m$ is a seed monomial if and only if the paths of $\bs m$ are non-intersecting.
\end{thm}

\begin{proof}

The only if direction follows immediately from Lemmas \ref{subpaths lemma} and \ref{non intersecting is equivalent}.

Conversely, let the paths of $\bs m$ be non-intersecting. We claim that any expansion of $\bs m$ results in a thin monomial with non-intersecting paths. Indeed, let $(u,v)$ be a bullet which can be expanded. Then, $(u+1,v+1)$ has no bullet. In addition, by the non-intersecting property, the corner at $(u+1,v+1)$ does not belong to paths of $\bs m$. Therefore, expanding a bullet $(u,v)$ does not produce any intersections, see Figure \ref{path expansion}. The case of a cross expansion is similar.
\end{proof}

\begin{figure}[H]
\begin{tikzpicture}
 
 \draw[step=0.5cm, gray, very thin] (-2,-0.5) grid (2.5,2.5);
 \node at (-1,0) {$\bm \times$};
 \node at (-1,1) {\textbullet};
 \node at (0.5,1) {$\bm \times$};
 \node at (0.5,1.5) {\textbullet};
 \node at (1.5,1.5) {$\bm \times$};    
\draw[ultra thick] (-1,0)--(-2,0);
\draw[ultra thick] (-1,0)--(-1,1);
\draw[ultra thick] (-1,1)--(0.5,1);
\draw[ultra thick] (0.5,1)--(0.5,1.5);
\draw[ultra thick] (0.5,1.5)--(1.5,1.5);
\draw[ultra thick] (1.5,1.5)--(1.5,2.5);

\node at (0,-1) {$hv$ chain};

\draw[step=0.5cm, gray, very thin] (5.99,-0.5) grid (10.5,2.5);
\node at (6.5,0.5) {\textbullet};
\node at (7.5,0.5) {$\bm \times$};
\node at (7.5,1) {\textbullet};
\node at (8.5,1) {$\bm \times$};
\node at (8.5,2) {\textbullet};
\draw[ultra thick] (6.5,0.5)--(6.5,-0.5);
\draw[ultra thick] (6.5,0.5)--(7.5,0.5);
\draw[ultra thick] (7.5,0.5)--(7.5,1);
\draw[ultra thick] (7.5,1)--(8.5,1);
\draw[ultra thick] (8.5,1)--(8.5,2);
\draw[ultra thick] (8.5,2)--(10.5,2);
\node at (8,-1) {$vh$ chain};

 \draw[step=0.5cm, gray, very thin] (-2,-4.5) grid (2.5,-1.5);
 \node at (-1,-4) {$\bm \times$};
 \node at (-1,-3) {\textbullet};
 \node at (0.5,-3) {$\bm \times$};
 \node at (0.5,-2.5) {\textbullet};
 \node at (1,-2.5) {$\bm \times$};    
 \node at (1,-2) {\textbullet};
\draw[ultra thick] (-2,-4)--(-1,-4);
\draw[ultra thick] (-1,-4)--(-1,-3);
\draw[ultra thick] (-1,-3)--(0.5,-3);
\draw[ultra thick] (0.5,-3)--(0.5,-2.5);
\draw[ultra thick] (0.5,-2.5)--(1,-2.5);
\draw[ultra thick] (1,-2.5)--(1,-2);
\draw[ultra thick] (1,-2)--(2.5,-2);

\node at (0,-5) {$hh$ chain};

\draw[step=0.5cm, gray, very thin] (5.99,-4.5) grid (10.5,-1.5);
\node at (6.5,-3.5) {\textbullet};
\node at (7.5,-3.5) {$\bm \times$};
\node at (7.5,-3) {\textbullet};
\node at (8.5,-3) {$\bm \times$};
\node at (8.5,-2.5) {\textbullet};
\node at (9.5,-2.5) {$\bm \times$};
\draw[ultra thick] (6.5,-3.5)--(6.5,-4.5);
\draw[ultra thick] (6.5,-3.5)--(7.5,-3.5);
\draw[ultra thick] (7.5,-3.5)--(7.5,-3);
\draw[ultra thick] (7.5,-3)--(8.5,-3);
\draw[ultra thick] (8.5,-3)--(8.5,-2.5);
\draw[ultra thick] (8.5,-2.5)--(9.5,-2.5);
\draw[ultra thick] (9.5,-2.5)--(9.5,-1.5);
\node at (8,-5) {$vv$ chain};

\end{tikzpicture}
\caption{The four types of chains.}
\label{chains}
\end{figure}
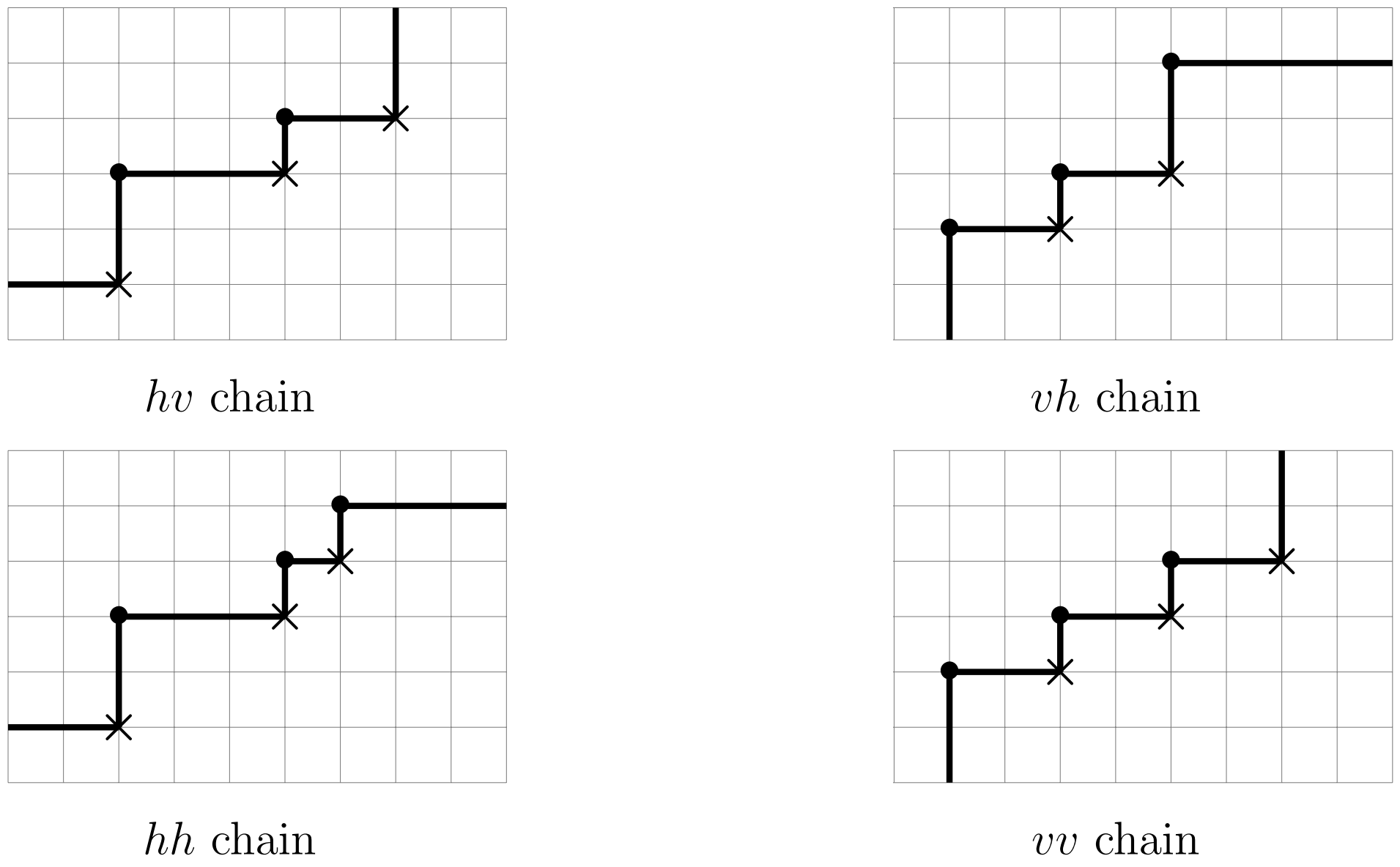

Let $\bm m$ be a seed monomial. We call a \textit{chain} of $\bm m$ a connected component of the paths of $\bm m $, so that $\bm m$ is a union of non-intersecting chains. A chain contains exactly two rays and divides the plane $\mathbb Z\times \mathbb Z$ into two disjoint regions. We have four types of chains depending on the rays: $hh$, $hv$, $vh$ and $vv$, as shown in Figure \ref{chains}. Here, $h$ stands for a horizontal ray and $v$ for a vertical ray.

We group  seed monomials according to how the rays of their chains are arranged. For this purpose, we do not distinguish between rays whose intersection is also a ray. We write, for example, the ray $x=u$, $y\gg 0$, for all rays of the form $\{(u,t): t\geq v\}$ with starting point $(u,v).$

For an $hh$ chain with rays $y=v$, $x\ll 0$, and  $y=u$, $x\gg 0$,  and for a $vv$ chain with rays $x=u$, $y\gg 0$, and  $x=v$, $y\ll 0$, we call $v-u\in\Z_{\geq 1}$ the {\it width} of the chain.

Let $l,m,p\in\Z_{\geq 0}$. Let
$\la,\mu,\nu, \tilde \la, \tilde \mu, \tilde \nu$ be finite strictly decreasing sequences of  half-integers such that 
 $|\la|=|\tilde \la|=l$, $|\mu|=|\tilde \mu|=m$, $|\nu|=|\tilde \nu|=p$.
 We write $\la=(\la_1,\dots,\la_l)$, where $\la_1>\la_2>\dots >\la_l$, $\tilde \la=(\tilde \la_1,\tilde \la_2,\dots ,\tilde \la_l)$, where $\tilde \la_1>\tilde \la_2>\dots >\tilde \la_l$, etc.

Notice that $hh$ and $vv$ chains cannot coexist in a seed monomial. Thus, we start with monomials without $hh$ chains. 

Let $\mathcal V_{\la,\mu,\nu}^{\tilde \la,\tilde \mu,\tilde \nu}$ be the set of all seed monomials having exactly the following $l+m+p$ chains:
\begin{itemize} 
\item $hv$ chains: there are $p$ chains of this type. Each chain has a horizontal ray  $y=\nu_i$, $x\ll 0,$ and a vertical ray  $x=\tilde \nu_i$, $y\ll 0$, for $i=1,\dots,p$;
\item $vv$ chains: there are $m$ chains of this type. Each chain has vertical rays  $x=\mu_i$, $y\gg 0$, and $x=\tilde \mu_i$, $y\ll 0$, $\tilde \mu_i>\mu_i$, for $i=1,\dots,m$; 
\item $vh$ chains: there are $l$ chains of this type. Each chain has a vertical ray  $x= \la_i$, $y\gg 0$, and a horizontal ray  $y=\tilde \la_i$, $x\gg 0$,  for $i=1,\dots,l$. 
\end{itemize}

We proceed to monomials without  $vv$ chains. 

Let $\mathcal H_{\la,\mu,\nu}^{\tilde \la,\tilde \mu,\tilde \nu}$ be the set of all seed monomials having exactly the following $l+m+p$ chains:
\begin{itemize}
\item $hv$ chains: there are $p$ chains of this type. Each chain has a horizontal ray  $y=\nu_i$, $x\ll 0,$ and a vertical ray  $x=\tilde \nu_i$, $y\ll 0$, for $i=1,\dots,p$;
\item $hh$ chains: there are $m$ chains of this type. Each chain has horizontal rays  $y=\mu_i$, $x\ll 0$, and $y=\tilde \mu_i$, $x\gg 0$, $\mu_i >\tilde \mu_i$, for $i=1,\dots,m$;
\item $vh$ chains: there are $l$ chains of this type. Each chain has a vertical ray  $x= \la_i$, $y\gg 0$, and a horizontal ray  $y=\tilde \la_i$, $x\gg 0$,  for $i=1,\dots,l$.
\end{itemize} 

If $m=0$, we have $\mathcal V_{\la,\varnothing,\nu}^{\tilde \la,\varnothing,\tilde \nu}=\mathcal H_{\la,\varnothing,\nu}^{\tilde \la,\varnothing,\tilde \nu}$. Here, $\varnothing$ stands for the empty sequence of half-integers. This case corresponds to monomials with $hv$  or $vh$ chains only. Henceforth, we reserve the notation $\mathcal H_{\la,\mu,\nu}^{\tilde \la,\tilde \mu,\tilde \nu}$ for the case $m>0$.

\begin{example}
The monomial in Figure \ref{thimonomialviz} 
belongs to the set  $$\mathcal H_{\varnothing,(5/2,3/2,-1/2) , \varnothing}^{\varnothing, (3/2,-1/2,-3/2), \varnothing}.$$
The monomial in Figure \ref{non intersecting} belongs to the set $$\mathcal V^{(3/2,-1/2),(7/2),(5/2,-5/2)}_{(5/2,-5/2),(-7/2),(-3/2,-5/2)}.$$
\end{example}

We state simple properties of $\mathcal V_{\la,\mu,\nu}^{\tilde \la,\tilde \mu,\tilde \nu}$  and $\mathcal H_{\la,\mu,\nu}^{\tilde \la,\tilde \mu,\tilde \nu}$ 
in the following lemma.

\begin{lem} \label{h and v sets} The sets of monomials
$\mathcal V_{\la,\mu,\nu}^{\tilde \la,\tilde \mu,\tilde \nu}$  and $\mathcal H_{\la,\mu,\nu}^{\tilde \la,\tilde \mu,\tilde \nu}$ satisfy the following properties.
\begin{enumerate}
\item The set $\mathcal V_{\la,\mu,\nu}^{\tilde \la,\tilde \mu,\tilde \nu}$ is non-empty if and only if
\begin{align*}   %\label{v assumptions}
 \tilde \mu_i>  \mu_i, \quad i=1,\dots,m, \qquad \la_l>\mu_1\ (l,m\neq 0), \qquad 
\tilde \mu_m>\tilde \nu_1\ (m,p\neq 0), 
\end{align*}

\item The set $\mathcal H_{\la,\mu,\nu}^{\tilde \la,\tilde \mu,\tilde \nu}$ is non-empty if and only if 
\begin{align*} %\label{h assumptions}
\mu_i > \tilde \mu_i, \quad i=1,\dots,m, \qquad \tilde \la_l>\tilde \mu_1\  (l,m\neq 0), \qquad \mu_m>\nu_1 \ (m,p\neq 0).   
\end{align*}

\item All sets of monomials  $\mathcal V_{\la,\mu,\nu}^{\tilde \la,\tilde \mu,\tilde \nu}$ and  $\mathcal H_{\la,\mu,\nu}^{\tilde \la,\tilde \mu,\tilde \nu}$ ($m\neq 0$) are disjoint.

\item Any seed monomial of integer type belongs to either  $\mathcal V_{\la,\mu,\nu}^{\tilde \la,\tilde \mu,\tilde \nu}$   or  $\mathcal H_{\la,\mu,\nu}^{\tilde \la,\tilde \mu,\tilde \nu}$.
\qed
\end{enumerate}
\end{lem}

The next theorem explicitly describes the expansions of any seed monomial of integer type.

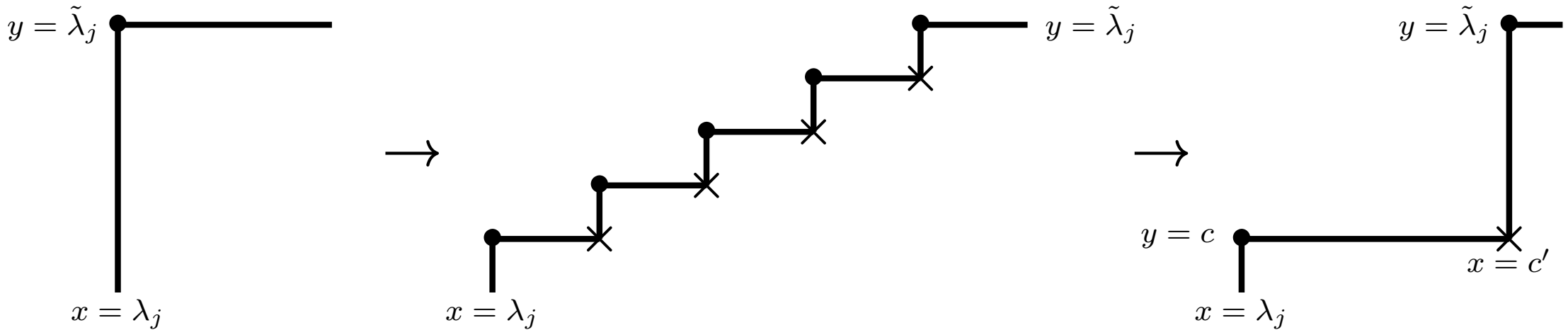
\begin{figure}[H]
\begin{tikzpicture}

\node at (0,0) {\textbullet};
\draw[ultra thick] (0,0)--(0,-2.5);
\draw[ultra thick] (0,0)--(2,0);

\node at (0,-2.7) {\text{\tiny $x=\lambda_j$}};
\node at (-0.6,0) {\text{\tiny $y=\tilde  \lambda_j$}};

\draw[->, thick] (2.5,-1.2)--(3,-1.2);

\node at (3.5,-2) {\textbullet};
\node at (4.5,-2) {$\bm \times$};
\node at (4.5,-1.5) {\textbullet};
\node at (5.5, -1.5) {$\bm \times$};
\node at (5.5,-1) {\textbullet};
\node at (6.5, -1) {$\bm \times$};
\node at (6.5, -0.5) {\textbullet};
\node at (7.5,-0.5) {$\bm \times$};
\node at (7.5,0) {\textbullet};

\draw[ultra thick] (3.5,-2)--(3.5,-2.5);
\draw[ultra thick] (3.5,-2)--(4.5,-2);
\draw[ultra thick] (4.5,-2)--(4.5,-1.5);
\draw[ultra thick] (4.5,-1.5)--(5.5,-1.5);
\draw[ultra thick] (5.5,-1.5)--(5.5,-1);
\draw[ultra thick] (5.5,-1)--(6.5,-1);
\draw[ultra thick] (6.5,-1)--(6.5,-0.5);
\draw[ultra thick] (6.5,-0.5)--(7.5,-0.5);
\draw[ultra thick] (7.5,-0.5)--(7.5,-0);
\draw[ultra thick] (7.5,-0)--(8.5,-0);

\node at (3.5,-2.7) {\text{\tiny $x=\lambda_j$}};
\node at (9.1,0) {\text{\tiny $y=\tilde \lambda_j$}};

\draw[->, thick] (9.5,-1.2)--(10,-1.2);

\node at (10.5,-2) {\textbullet}; 
\node at (13,-2) {$\bm \times$};
\node at (13,0) {\textbullet};
\draw[ultra thick] (10.5,-2)--(10.5,-2.5);
\draw[ultra thick] (10.5,-2)--(13,-2);
\draw[ultra thick] (13,-2)--(13,0);
\draw[ultra thick] (13,0)--(13.5,0);

\node at (10.5,-2.7) {\text{\tiny $x=\lambda_j$}};
\node at (12.4,0) {\text{\tiny $y=\tilde  \lambda_j$}};
\node at (9.9,-2) {\text{\tiny $y=c$}};
\node at (13,-2.2) {\text{\tiny $x=c'$}};

\end{tikzpicture}
\caption{Successive expansions of a $vh$ chain.}
\label{vh exp}
\end{figure}

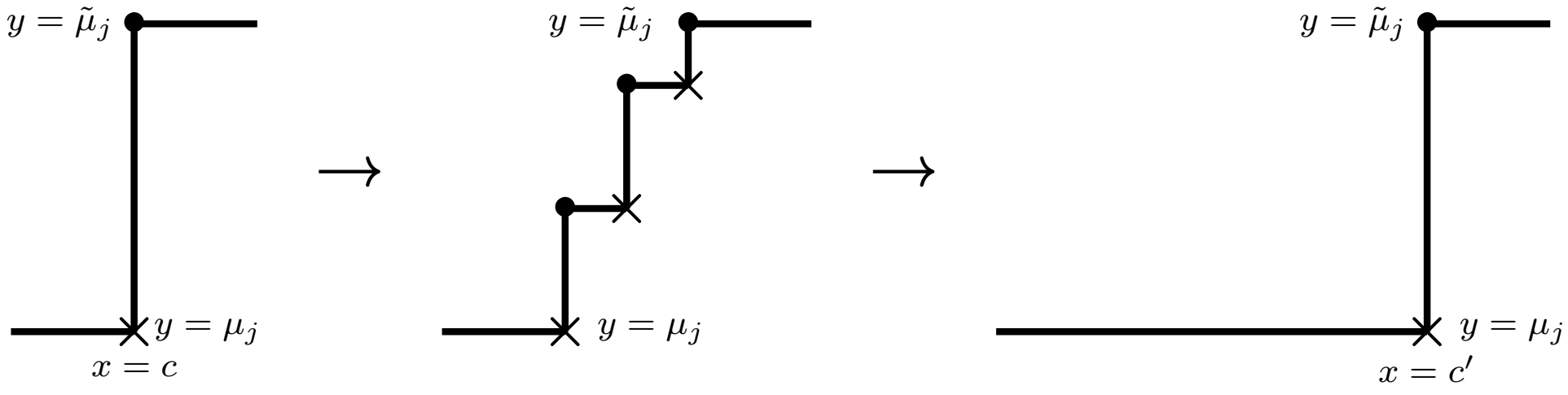
\begin{figure}[H]
\begin{tikzpicture}

 \node at (0,0) {$\bm \times$};
\node at (0,2.5) {\textbullet};
\draw[ultra thick] (0,0)--(-1,0);
\draw[ultra thick] (0,0)--(0,2.5);
\draw[ultra thick] (0,2.5)--(1,2.5);

\node at (0.6,0) {\text{\tiny $y=\mu_j$}};
\node at (-0.6,2.5) {\text{\tiny $y=\tilde \mu_j$}};
\node at (0,-0.3) {\text{\tiny $x=c$}};

\draw[->, thick] (1.5,1.3)--(2,1.3);

\node at (3.5,0) {$\bm \times$};
\node at (3.5,1) {\textbullet};
\node at (4,1) {$\bm \times $};
\node at (4,2) {\textbullet};
\node at (4.5,2) {$\bm \times$};
\node at (4.5,2.5) {\textbullet};

\draw[ultra thick] (3.5,0)--(2.5,0);
\draw[ultra thick] (3.5,0)--(3.5,1);
\draw[ultra thick] (3.5,1)--(4,1);
\draw[ultra thick] (4,1)--(4,2);
\draw[ultra thick] (4,2)--(4.5,2);
\draw[ultra thick] (4.5,2)--(4.5,2.5);
\draw[ultra thick] (4.5,2.5)--(5.5,2.5);

\node at (4.2,0) {\text{\tiny $y=\mu_j$}};
\node at (3.8,2.5) {\text{\tiny $y=\tilde  \mu_j$}};

\draw[->, thick] (6,1.3)--(6.5,1.3);

\node at (10.5,0) {$\bm \times$};
\node at (10.5,2.5) {\textbullet};

\draw[ultra thick] (7,0)--(10.5,0);
\draw[ultra thick] (10.5,0)--(10.5,2.5);
\draw[ultra thick] (10.5,2.5)--(11.5,2.5);

\node at (11.2,0) {\text{\tiny $y=\mu_j$}};
\node at (10.5,-0.3) {\text{\tiny $x=c'$}};
\node at (9.9,2.5) {\text{\tiny $y=\tilde \mu_j$}};

\end{tikzpicture}
\caption{Successive expansions of an $hh$ chain.}
\label {hh exp}
\end{figure}

\begin{thm}  \label{expansions of seed } Let $\bs m$ be a seed monomial of integer type. 
\begin{enumerate}
\item If $\bs m \in  \mathcal H_{\la,\mu,\nu}^{\tilde \la,\tilde \mu,\tilde \nu}$ for some $\lambda, \tilde \lambda, \mu, \tilde \mu, \nu, \tilde \nu,$  then $S(\bs m)= \mathcal H_{\la,\mu,\nu}^{\tilde \la,\tilde \mu,\tilde \nu}$.
\vspace{0.25cm}
\item If $\bs m \in  \mathcal V_{\la,\mu,\nu}^{\tilde \la,\tilde \mu,\tilde \nu}$ for some $\lambda, \tilde \lambda, \mu, \tilde \mu, \nu, \tilde \nu,$ then $S(\bs m)= \mathcal V_{\la,\mu,\nu}^{\tilde \la,\tilde \mu,\tilde \nu}$.
\end{enumerate}
\end{thm}
\begin{proof}

We just prove item $(i).$ The proof of $(ii)$ is completely analogous.

Let $\bs m \in  \mathcal H_{\la,\mu,\nu}^{\tilde \la,\tilde \mu,\tilde \nu}$. For every $i\in \hat I,$ we have $e_i(\bm m)\cup f_i(\bm m)\subset  \mathcal H_{\la,\mu,\nu}^{\tilde \la,\tilde \mu,\tilde \nu}$, since any ray of $\bm m$  does not change  under expansions, see Figure  \ref{path expansion}. Thus, $S(\bm m)\subseteq  \mathcal H_{\la,\mu,\nu}^{\tilde \la,\tilde \mu,\tilde \nu}.$

We prove the opposite inclusion by induction on the number of chains in $\bs m$.  The case of a single chain is straightforward. We show some expansions in Figures \ref{vh exp} and \ref{hh exp}. 

Let us show the induction step. Let $ \bs m'\in  \mathcal H_{\la,\mu,\nu}^{\tilde \la,\tilde \mu,\tilde \nu}$.
If $\la\neq \varnothing$, we modify the ``lowest" $vh$ chain (containing the ray $x=\la_1$, $y\gg 0$) in $\bs m$ and $\bs m'$ until it has the form shown in the right picture of Figure \ref{vh exp}, where the coordinates $c,c'$ are sufficiently large.
After that, the remaining chains can be  successively expanded to each other by the  induction hypothesis.

If $\la=\varnothing$, we modify the ``lowest" $hh$ chain (containing the ray $y=\mu_1$, $x\ll 0$) until it has the form shown on the right in Figure \ref{hh exp}, where $c'$ is sufficiently large, and use the induction hypothesis.

The other cases are similar.
\end{proof}

\subsection{Bounded seed monomials}

\begin{thm} 
Let $\bm m$ be a seed monomial of integer type. Then, $\bm m$ is bounded if and only if  $\bm m \in \mathcal V^{\varnothing, \tilde \mu, \varnothing}_{\varnothing, \mu, \varnothing} \sqcup \mathcal H^{\varnothing, \tilde \mu, \varnothing}_{\varnothing, \mu, \varnothing}.$ In other words, $\bs m$ is bounded if and only if the paths of $\bs m$ contain neither $hv$ nor $vh$ chains.
\label{bounded seed}
\end{thm}

\begin{proof}
Suppose that $\bm m$ has a $vh$ chain. Then, one can expand  the  ``lowest" $vh$ chain (containing the ray $x=\la_1$, $y\gg 0$) of $\bm m$ into a $vh$ chain  which has an arbitrarily large number of bullets of color $i$ and no crosses of color $i$, see the second picture in Figure \ref{vh exp}.  Therefore, $\bs m$ is not a bounded seed monomial.

The case of $hv$ chains is similar.

Suppose that $\bs m$ has only chains of type $hh$. Then, for the $hh$ chain with rays $y=\mu_j, x\ll 0$ and $y=\tilde \mu_j, x\gg 0,$
the number of bullets and the number of crosses are at most $\mu_j - \tilde \mu_j$, see the second picture in Figure \ref{hh exp}. Therefore, for $i\in \hat I,$ the sum of absolute values of the exponents is bounded: $\sum_a|k_{i,a}|\leq \sum_{j=1}^m (\mu_j  - \tilde\mu_j).$  It is easy to see that this bound is sharp.

The case of $vv$ chains is similar.
\end{proof}

\section{The construction of chain modules}  \label{construction chain}
In this section,  we recall the construction of various known $\mathcal E_n$-modules, namely the vector and  covector representations, the Fock and dual Fock representations, see \cite{FJMM1}, \cite{FJMM2}. We shall use these representations to construct the  $\mathcal E_n$-modules $H_{\la,\mu,\nu}^{\tilde \la,\tilde \mu,\tilde \nu}$ and $ V_{\la,\mu,\nu}^{\tilde \la,\tilde \mu,\tilde \nu}$, whose basis vectors are parametrized by the monomials in $\mathcal H_{\la,\mu,\nu}^{\tilde \la,\tilde \mu,\tilde \nu}$ and $\mathcal V_{\la,\mu,\nu}^{\tilde \la,\tilde \mu,\tilde \nu}$, where $\lambda, \mu, \nu, \tilde \lambda , \tilde \mu, \tilde \nu$ are as in Lemma \ref{h and v sets}.

\subsection{The cases of a single \texorpdfstring{$hh$}{} chain or \texorpdfstring{$vv$}{} chain}   \label{hh vv chains}
 We start by recalling the vector and covector representations $W(u)$ and $\check{W}(u),$ which in our language correspond to a single $hh$ and $vv$ chain of width $1$, respectively.

We denote by $Y(z)$ the rational function (cf.  (\ref{fundamental l weights})):
$$Y(z) =  q\frac{1-q^{-1}z}{1-qz}. $$

  For $u\in\C^\times$, the vector representation $W(u)$ is an $\mathcal E_n$-module with basis $\{[u]_j\}_{j\in \mathbb Z}$ and  action given by the formulas
\begin{align*}
&F_i(z) \cdot [u]_j  = \begin{cases}  \delta(q_1^{-j-1}u/z)[u]_{j+1}, & i\equiv j+1,
\\
0, & \mbox{otherwise},
\end{cases}
\\\\
&E_i(z) \cdot [u]_j = \begin{cases} \delta(q_1^{-j}u/z)[u]_{j-1}, & i\equiv j,
\\
0, & \mbox{otherwise},
\end{cases}
\\\\
&K^{\pm}_i(z)\cdot [u]_j = \begin{cases} Y(q^{-1}q_1^{-j-1}u/z)[u]_j, &i\equiv j+1,
\\
Y(qq_1^{-j}u/z)^{-1}[u]_j, &i\equiv j,
\\
[u]_j, & \mbox{otherwise},
\end{cases}
\end{align*}
where $i\in \hat I, j\in \mathbb Z$.

 For $u\in\C^\times$, the covector representation $\check W(u)$ is an $\mathcal E_n$-module with basis $\{[u]^\vee_j\}_{j\in \mathbb Z}$  and action given by 
\begin{align*}
&F_i(z) \cdot {[ u]_j^\vee}  = \begin{cases}  \delta(q_3^{-j-1}u/z){ [ u]_{j+1}^\vee}, & i\equiv -(j+1),
\\
0, &  \mbox{otherwise},
\end{cases}
\\\\
&E_i(z) \cdot [ u]_j^\vee = \begin{cases} \delta(q_3^{-j}u/z){[  u]}_{j-1}^\vee, & i \equiv -j,
\\
0, & \mbox{otherwise},
\end{cases}
\\\\
&K^{\pm}_i(z)\cdot {[u]}_j^\vee = \begin{cases} Y(q^{-1}q_3^{-j-1}u/z){[ u]}_j^\vee, &  i \equiv -(j+1),
\\
Y(qq_3^{-j}u/z)^{-1}{[u]}_j^\vee, &  i \equiv-j,
\\
{[u]}_j^\vee, & \mbox{otherwise},
\end{cases}
\end{align*}
where $i\in \hat I, j\in \mathbb Z$.

We remark that the covector representation $\check W(u)$ is the $\mathcal E_n(q_1,q_2,q_3)$-module obtained from the vector representation $W(u)$ of $\mathcal E_n(q_3,q_2,q_1)$ by twisting with the isomorphism $\iota: \mathcal E_n(q_1,q_2,q_3)\rightarrow \mathcal E_n(q_3,q_2,q_1)$, see $(\ref{transpose isomorphism}).$

Notice that we have an isomorphism of modules $W(u)$ and $W(uq_1^n)$ sending $[u]_j$ to  $[uq_1^n]_{j+n}$. Similarly, 
we have an isomorphism of modules $\check W(u)$ and $\check W(uq_3^n)$ sending $[u]_j^\vee$ to  $[uq_3^n]^\vee_{j+n}$.

The $q$-characters of $W(u)$ and $\check W (u)$ consist of monomials of integer type if and only if \begin{equation} \label{uinteger} u=q_1^{-a}q_3^{-b}\in q_1^{\mathbb Z}q_3^{\mathbb Z},  \quad \text{where} \ \ a-b\in n\Z. \end{equation}  

\begin{figure}[H]
\begin{tikzpicture}

\rzerobox at (0,0);
\threebox at (-0.5,0);
\twobox at (-1,0);
\onebox at (-1.5,0);
\zerobox at (-2,0);
\threebox at (-2.5,0);

\node at (-1.5,0.5) {\textbullet};
\node at (-1.5,0) {$\bm \times$};

\draw[ultra thick] (-2.5,0)--(-1.5,0);
\draw[ultra thick] (-1.5,0)--(-1.5,0.5);
\draw[ultra thick] (-1.5,0.5)--(0.5,0.5);

\node at (-1,-1) {$[u]_{-4} \in W_{1/2}^{-1/2}$};

\draw (2.5,-1.5)--(2.5,1.3);

\zerobox at (5,-1.5);
\onebox at (5,-1);
\twobox at (5,-0.5);
\threebox at (5,0);
\rzerobox at (5,0.5);

\node at (5, -0.5) {\textbullet};
\node at (5.5,-0.5) {$\bm \times$};
\draw[ultra thick] (5,-1.5)--(5,-0.5);
\draw[ultra thick] (5,-0.5)--(5.5,-0.5);
\draw[ultra thick] (5.5,-0.5)--(5.5,1);
 w
\node at (7.5, 0) {$[u ]_2^\vee \in \check W_{-1/2}^{1/2}$};

\end{tikzpicture}
\caption{The monomials corresponding to the basis vectors $[u]_{-4}$ and $[u]_2^\vee$ for $n=4$ and $ u =1.$}
\label{width 1}
\end{figure}
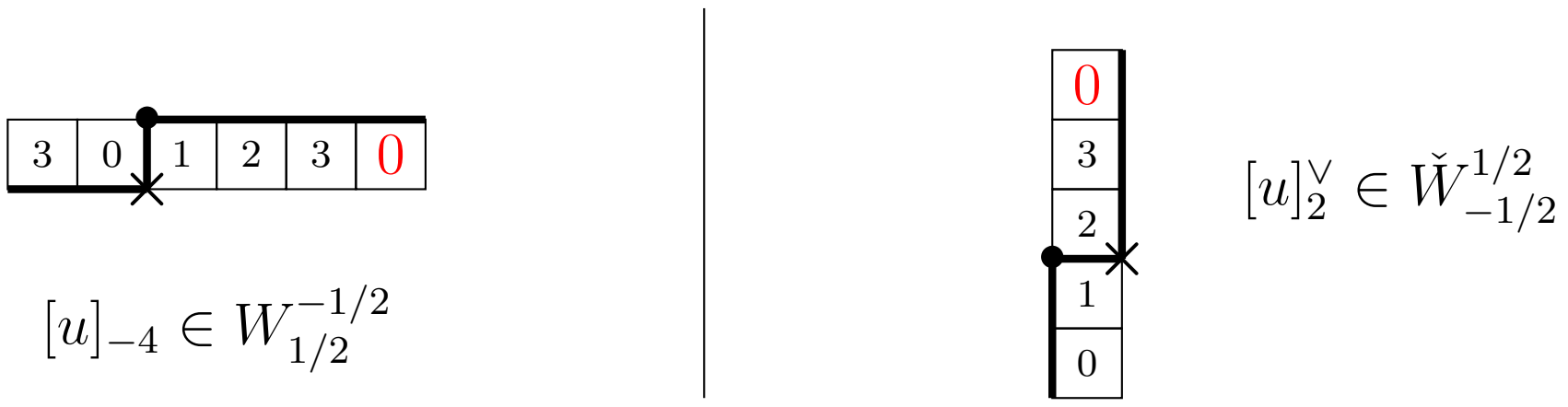

Let $u$ be as in (\ref{uinteger}). We  picture the monomials corresponding to the vectors $[u]_j$ and $[u]_j^\vee$, $j\in\Z$, using our conventions, see Figure \ref{width 1}.

The monomial corresponding to  $[u]_j$ is an $hh$ chain of width $1,$ where we have a bullet of color $j+1 \ (\operatorname{mod}\ n)$ at $(a,b) + (j+\frac{1}{2},-\frac{1}{2})$ and a cross of color $j \ (\operatorname{mod}\ n)$ at $(a,b)+ (j+\frac{1}{2}, \frac{1}{2})$. In particular, this $hh$ chain has the rays $y=b+\frac{1}{2}$, $x\ll 0$ and $y=b-\frac{1}{2}$, $x\gg 0$. The $q$-character of $W(u)$ is the sum of all monomials in $\mathcal H^{\varnothing,(b-1/2),\varnothing}_{\varnothing,(b+1/2),\varnothing}$. We often  denote the module $W(u)$ by  $W_{b+1/2}^{b-1/2}.$

The monomial corresponding to  $[u]_j^\vee$ is a $vv$ chain of width $1,$ where we have a bullet at $(a,b)+ (-\frac{1}{2}, j+\frac{1}{2})$ and a cross at $(a,b) + (\frac{1}{2}, j+\frac{1}{2})$. In particular, this $vv$ chain contains the rays $x=a-\frac{1}{2}$, $y\gg 0$ and $x=a+\frac{1}{2}$, $y\ll 0$. The $q$-character of $\check W(u)$ is the sum of all monomials in $\mathcal V^{\varnothing,(a+1/2),\varnothing}_{\varnothing,(a-1/2),\varnothing}$. We often denote the module $\check W(u)$ by $\check W_{a-1/2}^{a+1/2}$.

We consider tensor products of vector  representations and tensor products of covector representations using the coproduct (\ref{coproduct}).

\begin{prop} Let $u\in \mathbb C^\times, k\in \mathbb Z_{> 0}.$ \label{submodule tensor vector}
\begin{enumerate}
\item  The subspace of $W(u)\otimes W(uq_2)\otimes \cdots \otimes W(uq_2^{k-1})$ with basis  $$\big \{[u]_{r_0}\otimes [ uq_2]_{r_1}\otimes \cdots \otimes [uq_2^{k-1}]_{r_{k-1}}:\  r_0>r_1>\cdots >r_{k-1}\big \}$$  has a structure of an irreducible thin locally finite $\mathcal E_n$-module of level $1$. 
\item The subspace of $\check W(u)\otimes \check W(uq_2)\otimes \cdots \otimes \check W(uq_2^{k-1})$ with basis  $$\big \{ [u]_{r_0}^\vee \otimes [uq_2]_{r_1}^\vee \otimes \cdots \otimes [uq_2^{k-1}]_{r_{k-1}}^\vee: \  r_0>r_1>\cdots>r_{k-1}\big \}$$ has a structure of an irreducible thin locally  finite $\mathcal E_n$-module of level $1$. 
\end{enumerate} 
\end{prop}
\begin{proof}
 The proof of the proposition is standard, see Lemma 3.2 in \cite{FJMM1}, or Lemma 3.5 and Lemma 3.6 in \cite{BM}. 

 One checks that the action given by the coproduct on the subspace is well-defined and preserves the subspace by an explicit computation.  
\end{proof}

Let $u\in \mathbb C^\times$ be as in  \eqref{uinteger}. Then, each basis vector in the module in part $(i)$ of Proposition  \ref{submodule tensor vector} corresponds to an $hh$ chain of width $k$. Moreover, the $q$-character of the module is the sum of all monomials in $\mathcal H^{\varnothing,(b-1/2),\varnothing}_{\varnothing,(b+k-1/2),\varnothing}$. We denote this module by $W_{b+k-1/2}^{b-1/2}.$

Similarly, each basis vector in the module of part $(ii)$ of Proposition \ref{submodule tensor vector} corresponds to a $vv$ chain of width $k$. Moreover,  the $q$-character of the module is the sum of all monomials in $\mathcal V^{\varnothing,(a+k-1/2),\varnothing}_{\varnothing,(a-1/2),\varnothing}.$ We denote this module by $\check W_{a-1/2}^{a+k-1/2}.$

\subsection{The cases of a single \texorpdfstring{$vh$}{} or \texorpdfstring{$hv$ chain}{}} \label{vh hv chains}

We first recall the Fock modules.

Fix $s\in\hat I$ and $u\in \mathbb C^\times$. Consider the following infinite tensor product of vector representations: $$W(u)\otimes W(uq_2)\otimes W(uq_2^2)\otimes  \cdots = \bigotimes_{k\, =\, 0}^\infty  W(uq_2^k) .$$
 Let $\mathcal F_s(u)$ be the subspace spanned by all vectors of the form
\begin{equation*} |\bm r\rangle   = [u]_{r_0}\otimes [uq_2]_{r_1}\otimes [uq_2^2]_{r_2} \otimes \cdots, \end{equation*}
where $\bm r = (r_0,r_1,r_2,\dots)$ is a strictly decreasing sequence of integers which differs from the sequence $(s-1,s-2,s-3,\dots)$ only in finitely many places. Let  $\ell(\bm r)\in \mathbb Z_{\geq 0}$ be the smallest integer such that \begin{equation} \label{almost vacuum sequence} r_j = s-j-1, \quad  j\geq \ell(\bm r).\end{equation}

The vector space $\mathcal{F}_s(u)$ has a structure of irreducible thin locally finite $\mathcal{E}_n$-module of level $q$ with the action defined as follows.

Let $\bm r$ be a strictly decreasing sequence of integers satisfying   (\ref{almost vacuum sequence}). Then, the matrix coefficients (in bra-ket notation) of the action are given by 
\begin{align*}
\langle \bm r +\bm 1_j | F_i(z)  |\bm r \rangle  = &\begin{cases}  c_{j,\bm r}\cdot \delta(q_1^{-r_j-1}q_2^ju/z), \quad &i \equiv r_j+1, \\
0, &\mbox{otherwise},
\end{cases} \\\\
\langle \bm r |E_i(z) | \bm r+\bm 1_j\rangle = &\begin{cases} c_{j,\bm r}' \cdot \delta(q_1^{-r_j-1}q_2^{j}u/z), \quad &i\equiv r_j+1,\\
0, &\mbox{otherwise},\\
\end{cases}
\\ \\
  \langle \bm r|K^{\pm}_i(z)  |\bm r \rangle = &\prod_{\substack{k=0 \\ i\equiv r_k+1}}^{\ell(\bm r)} Y(q^{-1}q_1^{-r_k-1}q_2^{k}u/z) \prod_{\substack{k=0 \\ i \equiv r_k}}^{\ell(\bm r)-1} Y(qq_1^{-r_k}q_2^{k}u/z)^{-1},
\end{align*}
where $i\in \hat I$, $j\in\Z_{\geq 0}$, $\bm r + \bm 1_j = (r_0,\dots, r_j+1,\dots)$,  and
$c_{j,\bm r},$ $c'_{j,\bm r}$ are constants which are non-zero if and only if the sequence $\bm r + \bm 1_j$ is strictly decreasing. 

We do not use the explicit form  of the constants $c_{j,\bm r},$ $c'_{j,\bm r}$ in this section. We  discuss these constants later in Section \ref{constants section}.

All the other matrix coefficients  are set to be zero.  

The $\mathcal E_n$-module $\mathcal{F}_s(u)$ is called the {\it Fock module of color $s$} (in parallel realization). The Fock module is a highest $\ell$-weight module with highest $\ell$-weight vector $\ket{(s-1,s-2,s-3,...)}$ and highest $\ell$-weight ${\bs s}_{q^{-1}q_1^{-s}u}\in \hat{\mathcal P}$.

\medskip

We now define the dual Fock modules.

We keep  $s\in \hat I$ and $u\in \mathbb C^\times$ fixed. Consider the infinite tensor product of covector representations:
$$ \cdots \otimes \check{W}(uq_2^{-2})\otimes \check{W}(uq_2^{-1})\otimes \check{W}(u) =\bigotimes_{ k\,  = \, -\infty}^0 \check W(uq_2^k). $$

 Let $\check{\mathcal F}_s(u)$ be the subspace spanned by all vectors of the form
 \begin{equation*}|\bm r \rangle^\vee  = \cdots \otimes [uq_2^{-2}]_{r_2}^\vee \otimes [uq_2^{-1}]_{r_1}^\vee \otimes [u]_{r_0}^\vee, \end{equation*}
where $\bm r = (r_0,r_1,r_2,\dots)$ is a strictly increasing sequence of integers which differs from the sequence $(-s,1-s,2-s,\dots)$ only in finitely many places. Let $\ell(\bm r)\in \mathbb Z_{\geq 0}$ be the smallest integer such that \begin{equation}r_j = j-s, \quad j\geq \ell(\bm r).\label{vacuum sequence dual}\end{equation}

The vector space $\check{\mathcal{F}}_s(u)$ has a structure of irreducible thin locally finite $\mathcal{E}_n$-module of level $q^{-1}$ with the action defined as follows.

Let $\bm r$ be a strictly increasing sequence of integers satisfying   (\ref{vacuum sequence dual}). Then, the matrix coefficients  of the action are given by 
\begin{align*}
^\vee\langle \bm r -\bm 1_j | E_i(z)  |\bm r \rangle^\vee  = &\begin{cases}  \check c_{j,\bm r}\cdot \delta(q_2^{-j}q_3^{-r_j}u/z), \quad &i \equiv -r_j, \\
0, &\mbox{otherwise},
\end{cases} \\\\
^\vee\langle \bm r |F_i(z) | \bm r-\bm 1_j\rangle^\vee = &\begin{cases} \check{c}^{\,'}_{j,\bm r} \cdot \delta(q_2^{-j}q_3^{-r_j}u/z), \quad &i\equiv -r_j,\\
0, &\mbox{otherwise},\\
\end{cases}
\\ \\
^\vee\langle \bm r | K^{\pm}_i(z)  |\bm r \rangle^\vee = &\prod_{\substack{k=0 \\ i \equiv -r_k}}^{\ell(\bm r)} Y(qq_2^{-k}q_3^{-r_k}u/z)^{-1} \prod_{\substack{k=0 \\ i \equiv -r_k-1}}^{\ell(\bm r)-1} Y(q^{-1}q_2^{-k}q_3^{\,-r_k-1}u/z),
\end{align*}
where $i\in \hat I$, $j\in \Z_{\geq 0}$, $\bm r - \bm 1_j = (r_0,\dots, r_j-1,\dots)$, and the constants 
$\check c_{j,\bm r}, \check c^{\,'}_{j,\bm r} $ (see Section \ref{constants section}) 
are non-zero if and only if the sequence 
$\bm r - \bm 1_j$ is strictly increasing.  

All other matrix coefficients are set to be zero.

The $\mathcal E_n$-module $\check{\mathcal{F}}_s(u)$ is called the {\it dual Fock module of color s} (in parallel realization). The dual Fock module is a lowest $\ell$-weight module with lowest $\ell$-weight vector $\ket{(-s,1-s,2-s,...)}^\vee$ and lowest $\ell$-weight $\bs s_{qq_3^su}^{-1}\in \hat{\mathcal P}$.

The $q$-characters of $\mathcal{F}_s(u)$ and $\check{\mathcal{F}}_s(u)$  consist of monomials of integer type if and only if $u$ is of the form \eqref{uinteger}.

For $u\in \mathbb C^\times$ as in (\ref{uinteger}), the monomial corresponding to the vector $|\bm r\rangle \in \mathcal{F}_s(u)$ is a $vh$ chain with rays $x = a+s-\frac{1}{2}, y\gg 0$ and $y = b-\frac{1}{2}, x\gg 0.$ Hence, the $q$-character of $\mathcal F_s(u)$ is the sum of all monomials in $\mathcal V_{(a+s-1/2),\varnothing, \varnothing }^{(b-1/2),\varnothing, \varnothing }$. We  denote the module $\mathcal F_s(u)$ by ${\mathcal F}_{a+s-1/2}^{\, b-1/2}.$

Similarly, for $u\in \mathbb C^\times$ as in (\ref{uinteger}), the monomial corresponding to the vector $|\bm r\rangle^\vee \in \check{\mathcal F}_s(u)$ is an $hv$ chain with rays $y = b-s+\frac{1}{2}, x\ll 0$ and $x = a+\frac{1}{2}, y\ll 0.$ Hence, the $q$-character of $\check{\mathcal F}_s(u)$ is the sum of all monomials in $\mathcal V_{\varnothing, \varnothing, (b-s+1/2) }^{\varnothing, \varnothing, (a+1/2) }$. We often denote the module $\check{\mathcal F}_s(u)$ by $\check{\mathcal F}_{b-s+1/2}^{\, a+1/2}.$

\subsection{The case of multiple chains}

In this section, we consider tensor products of $\mathcal E_n$-modules from Sections \ref{hh vv chains}, \ref{vh hv chains}. That allows us to construct representations 
whose basis vectors are parametrized by  unions of non-intersecting chains.

For $\lambda, \mu, \nu, \tilde \lambda, \tilde \mu, \tilde \nu$ as in part $(i)$ of Lemma \ref{h and v sets}, let
\begin{equation} \label{v tensor}  V_{\lambda, \mu, \nu}^{\tilde \lambda, \tilde \mu, \tilde \nu} \subset  \bigg (\mathcal F_{\lambda_1}^{\tilde \lambda_1}\otimes \cdots \otimes \mathcal F_{\lambda_l}^{\tilde \lambda_l}\bigg)\otimes \bigg( \check{W}_{\mu_1}^{\tilde \mu_1}\otimes \cdots \otimes \check{W}_{\mu_m}^{\tilde \mu_m}\bigg) \otimes \bigg( \check{\mathcal F}_{\nu_1}^{\tilde \nu_1}\otimes \cdots \otimes \check{\mathcal F}_{\nu_p}^{\tilde \nu_p}\bigg)\end{equation}
be the subspace spanned by tensor products of basis vectors of factors such that the chains corresponding to monomials of these basis vectors are non-intersecting.

For $\lambda, \mu, \nu, \tilde \lambda, \tilde \mu, \tilde \nu$ as in part $(ii)$ of Lemma \ref{h and v sets} with $m>0$, let
\begin{equation} \label{h tensor} H_{\lambda, \mu, \nu}^{\tilde \lambda, \tilde \mu, \tilde \nu} \subset  \bigg (\mathcal F_{\lambda_1}^{\tilde \lambda_1}\otimes \cdots \otimes \mathcal F_{\lambda_l}^{\tilde \lambda_l}\bigg)\otimes \bigg(  {W}_{\mu_1}^{\tilde \mu_1}\otimes \cdots \otimes {W}_{\mu_m}^{\tilde \mu_m}\bigg) \otimes \bigg( \check{\mathcal F}_{\nu_1}^{\tilde \nu_1}\otimes \cdots \otimes \check{\mathcal F}_{\nu_p}^{\tilde \nu_p}\bigg)\end{equation}
be the subspace spanned by tensor products of basis vectors of factors such that the chains corresponding to monomials of these basis vectors are non-intersecting.

In general, the Drinfeld coproduct \eqref{coproduct} is not well-defined on arbitrary tensor products of modules corresponding to single chains. However, the next theorem claims that it is well-defined on the subspaces $V_{\lambda, \mu, \nu}^{\tilde \lambda, \tilde \mu, \tilde \nu}$ and $H_{\lambda, \mu, \nu}^{\tilde \lambda, \tilde \mu, \tilde \nu}$. Moreover,  these subspaces are submodules.

\begin{thm} \label{H and V modules} Let $\lambda, \mu, \nu, \tilde \lambda, \tilde \mu, \tilde \nu$ be finite strictly decreasing sequences of half-integers.

\begin{enumerate}
\item For $\lambda, \mu, \nu, \tilde \lambda, \tilde \mu, \tilde \nu$ as in part $(i)$ of Lemma \ref{h and v sets}, the Drinfeld coproduct \eqref{coproduct} defines on the 
subspace $V_{\lambda, \mu, \nu}^{\tilde \lambda, \tilde \mu, \tilde \nu}$ a structure of irreducible thin locally finite $\mathcal E_n$-module of level $q^{l-p}$. Moreover, the $q$-character of $V_{\lambda, \mu, \nu}^{\tilde \lambda, \tilde \mu, \tilde \nu}$ is  the sum of all monomials in $\mathcal V_{\lambda, \mu, \nu}^{\tilde \lambda, \tilde \mu, \tilde \nu}$.

\item For $\lambda, \mu, \nu, \tilde \lambda, \tilde \mu, \tilde \nu$ as in part $(ii)$ of Lemma \ref{h and v sets}, the Drinfeld coproduct \eqref{coproduct} defines on the
subspace $H_{\lambda, \mu, \nu}^{\tilde \lambda, \tilde \mu, \tilde \nu}$ a structure of irreducible thin locally finite $\mathcal E_n$-module of level $q^{l-p}$. Moreover, the $q$-character of $H_{\lambda, \mu, \nu}^{\tilde \lambda, \tilde \mu, \tilde \nu}$ is  the sum of all monomials in $\mathcal H_{\lambda, \mu, \nu}^{\tilde \lambda, \tilde \mu, \tilde \nu}$.
\end{enumerate}
\end{thm}
\begin{proof}
Let $\bm m$ be a single chain monomial. For any $i\in\hat I$, the eigenvalue of $K_i^\pm(z)$ on the corresponding vector $v_{\bm m}$ in the single chain module is determined by the positions of bullets and crosses in $\bm m$ and is explicitly given by
\begin{align*}
\prod_{\substack{(x-\frac{1}{2},\,y-\frac{1}{2}),\  x-y\,\equiv \,i, \\   \text{bullets}}} q\ \frac{z-q_1^{-x+1}q_3^{-y+1}}{z-q_1^{-x}q_3^{-y}}\prod_{\substack{(x+\frac{1}{2},\,y+\frac{1}{2}),\  x-y\,\equiv \,i, \\   \text{crosses}}} q^{-1}\ \frac{z-q_1^{-x-1}q_3^{-y-1}}{z-q_1^{-x}q_3^{-y}}.
\end{align*}
The action of $F_i(z)$ on $v_{\bm m}$ in the
single chain module corresponds to the expansion of $\bm m $ at bullets of color $i$ as in the left picture of Figure \ref{path expansion}. The matrix coefficient of expanding a bullet at $(x-\frac{1}{2},y-\frac{1}{2})$ contains the delta function $\delta(q_1^{-x}q_3^{-y}/z)$. 

Consider the tensor product as in \eqref{v tensor} or \eqref{h tensor}, and a tensor product of vectors which contains $v_{\bs m}$ as a tensor factor.
Consider the action of $F_i(z)$ given by the Drinfeld coproduct \eqref{coproduct}. It is a sum of $m+l+p$ terms. The summand where $F_i(z)$ acts on $v_{\bs m}$
is not well-defined only if we have a pole of $K_i^\pm(z)$ coming from a single  chain corresponding to a module  to the left of $v_{\bs m}$ in the tensor product. Moreover, the corresponding chain should have a cross at $(x+\frac{1}{2},y+\frac{1}{2})$.
Clearly, that is impossible if the chains are non-intersecting.

In addition, if after expansion of the chain $\bs m$ at a bullet $(x-\frac{1}{2}, y-\frac{1}{2})$ we obtain a chain which intersects another chain $\bs n$, then the chain $\bs n$ corresponds to a vector in a module on the left of $v_{\bs m}$, and $\bs n$ has a bullet at $(x+\frac{1}{2},y+\frac{1}{2})$. But then the corresponding eigenvalue of $K_i^\pm(z)$ on the vector corresponding to $\bs n$ has a zero at the support of the delta function $\delta(q_1^{-x}q_3^{-y}/z),$ so  the matrix coefficient vanishes. Therefore $F_i(z)$ preserves the subspace.

The case of $E_i(z)$ is similar.

The corresponding modules are clearly irreducible, thin and locally finite.
\end{proof}

We call the modules $V_{\lambda, \mu, \nu}^{\tilde \lambda, \tilde \mu, \tilde \nu}$ and $H_{\lambda, \mu, \nu}^{\tilde \lambda, \tilde \mu, \tilde \nu}$ of Theorem \ref{H and V modules} the \textit{chain modules of type $v$ and type $h$}, respectively. We call the sequences of half-integers $\lambda, \mu, \nu, \tilde \lambda, \tilde \mu, \tilde \nu$ the \textit{ray configuration} of the corresponding chain module.

\section{Injectivity of the \texorpdfstring{$q$}{}-character map and the main results}   \label{injectivity of q character}
\subsection{Linked chain modules}   

\label{linked}
Recall that the actions of $F_i(z)$ and $E_i(z)$ always change the $\ell$-weight of a vector by an affine root. In this section we give a combinatorial description of pairs of chain modules such that their $q$-characters have monomials differing by an affine root.

Let $V$ and $W$ be two different chain modules such that there exist monomials $\bs m\in\chi_q(V)$ and  $\bs m'\in\chi_q(W)$ with $\bs m=\bs A^{\pm 1}_{i,c}\bs m'$  for some $i\in\hat I$, $c\in \C^\times$. We call such modules $V$ and $W$ \textit{linked}.

Let $\bm m \in \chi_q(V)$, where $V$ is a chain module. A box $\boxxslash = (a,b) $ is called a \textit{linking box} of $\bm m$ if at least two of its corners belong to infinite rays from two distinct chains lying on the same vertical line or on the same horizontal line.  In other words, for $\varepsilon \in \{\pm1\}$, if the two corners have coordinates $(a+\varepsilon/2,b-\frac{1}{2})$, $(a+\varepsilon/2,b+\frac{1}{2})$, then $(a+\varepsilon/2,b-\frac{1}{2})$ is a cross,  $(a+\varepsilon/2,b+\frac{1}{2})$ is a bullet and there are no other bullets or crosses with first coordinate $a+\varepsilon/2$.
If the two corners have coordinates $(a-\frac{1}{2},b+\varepsilon/2)$, $(a+\frac{1}{2},b+\varepsilon/2)$, then $(a-\frac{1}{2},b+\varepsilon/2)$ is a cross,  $(a+\frac{1}{2},b+\varepsilon/2)$ is a bullet and there are no other bullets or crosses with the second coordinate $b+\varepsilon/2$.

The next lemma shows that the linking of two monomials always corresponds to a linking box.

\begin{lem} \label{must be linking box} Let $V$ and $W$ be linked chain modules. Assume that the number of rays in $V$ is not less than the number of rays in $W$. Then, for all
$\bs m\in\chi_q(V)$ and  $\bs m'\in\chi_q(W)$ such that $\bs m=\bs A_{i,c}\bs m'$ or $\bs m=\bs A^{-1}_{i,c}\bs m'$, the color $i$ and scalar $c$ are the color and content of a linking box of $\bm m$.
\end{lem}
\begin{proof}
Let $\bm m \in \chi_q(V)$ and $\bm m' \in \chi_q(W)$ be monomials such that  $\bm m' = \bm A_{i,c}^{-1}\cdot \bm m$.  Since $\bs m$ and $\bs m'$ are integer monomials, there exists a unique box $\boxxslash = (a,b) \in \mathbb Z^2$ such that $q_1^{-a}q_3^{-b}= c$ and $a-b\equiv i$. We show that $\boxxslash$ is a linking box of $\bm m$.

Consider the number of cancellations between $\bm A^{-1}_{i,c}$ and $\bs m$. If that number is zero or one, then the number of chains of $\bs m'$ is larger, contradicting our assumption.

Let the number of cancellations be two. Then, the two canceled variables in $\bs m$ cannot be both crosses or both bullets due to the fact that $\bs m'$ is thin, see  subpictures 5 and 6 in Lemma \ref{subpaths lemma}.  If there is a  bullet at $(a-\frac{1}{2},b-\frac{1}{2})$ which is canceled, then $\bm m'$ is the result of a regular expansion of $\bm m$. As a result, we have $\bs m'\in\chi_q(V),$ which is a contradiction. If $\bs m$ has a bullet $(a+\frac{1}{2},b+\frac{1}{2})$ and a cross at $(a+\frac{1}{2},b-\frac{1}{2})$, then this cross and bullet belong to different chains. Moreover, since $\bs m'$ is thin, then there are no other crosses or bullets with the first coordinate $a+\frac{1}{2}$, see subpictures 1-4 in Lemma \ref{subpaths lemma}. Therefore, the box $(a,b)$ is a linking box. The same conclusion follows if $\bm m$ has a bullet at $(a+\frac{1}{2}, b+\frac{1}{2})$ and a cross at $(a-\frac{1}{2} ,b+\frac{1}{2}).$

The cases of number of cancellations being 3 or 4, and the case of $\bs m=\bm A_{i,c}\bs m'$ are treated similarly.

\end{proof}

Now we classify all  pairs of linked modules. 

Let $V$ be a chain module. Assume that the ray configuration of $V$ contains two rays on the same line. Moreover, assume that the chains corresponding to those two rays are neighboring, in the sense that there are no other chains in between. We introduce two types of new ray configurations emerging from the ray configuration of $V$.

The first type of configuration changes the two rays which were on the same line $x=a$ or $y=b$ to the two rays $x=a\pm 1$ or $y=b\pm 1$, while keeping all other rays the same. 

The second type of configuration removes the two rays which were on the same line to create a single chain formed by ``gluing" the remaining rays of these two chains, while keeping all other rays the same. Gluing may also remove both chains, see Figure \ref{glue move} second row.

The new configuration may or may not be a ray configuration of a chain module. If the new configuration defines a chain module $W$, we say that $W$ is obtained from $V$ by\textit{ shifting} two rays in the first case and \textit{gluing} two chains in the second case.

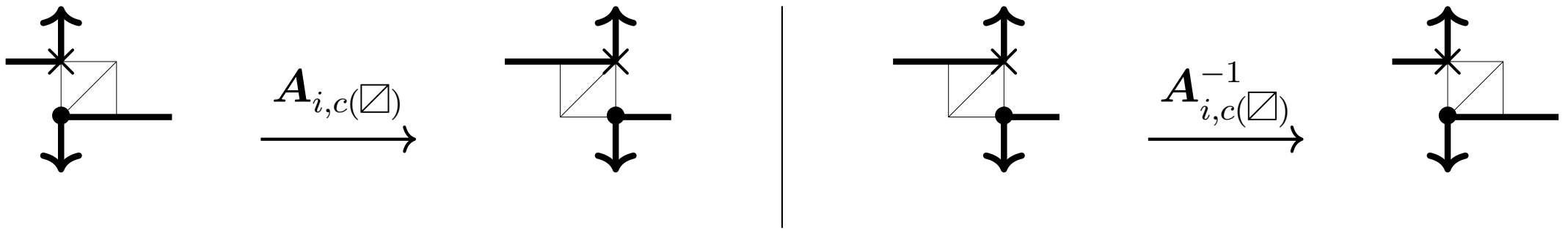
\begin{figure}[H]   
\centering
\begin{tikzpicture}

\begin{scope}[shift={(0,0)}]

\thindiagonalbox at (0,0);

\node at (0,0.5) {$\bm\times$};
\node at (0,0) {\textbullet};

\draw[->, ultra thick] (0,0)--(0,-0.5);
\draw[ultra thick] (0,0)--(1,0);
\draw[ultra thick] (0,0.5)--(-0.5,0.5);
\draw[->, ultra thick] (0,0.5)--(0,1);

\node at (2.5,0.2) {$\bm A_{i,c(\boxxslash)}$};
\draw[->,thick] (1.8,-0.2)--(3.2,-0.2);

\thindiagonalbox at (4.5,0);

\node at (5,0.5) {$\bm\times$};
\node at (5,0) {\textbullet};
\draw[ultra thick] (4,0.5)--(5,0.5);
\draw[->, ultra thick] (5,0.5)--(5,1);
\draw[ultra thick] (5,0)--(5.5,0);
\draw[->, ultra thick] (5,0)--(5,-0.5);

\draw (6.5,-1)--(6.5,1);

\thindiagonalbox at (8,0);

\node at (8.5,0.5) {$\bm\times$};
\node at (8.5,0) {\textbullet};

\draw[ultra thick] (8.5,0.5)--(7.5,0.5);
\draw[->,ultra thick] (8.5,0.5)--(8.5,1);
\draw[->,ultra thick] (8.5,0)--(8.5,-0.5);
\draw[ultra thick] (8.5,0)--(9,0);

\node at (10.5,0.2) {$\bm A^{-1}_{i,c(\boxxslash)}$};
\draw[->,thick] (9.8,-0.2)--(11.2,-0.2);

\thindiagonalbox at (12.5,0);

\node at (12.5,0) {\textbullet};
\node at (12.5,0.5) {$\bm \times$};
\draw[ultra thick] (12.5,0.5)--(12,0.5);
\draw[->,ultra thick] (12.5,0.5)--(12.5, 1);
\draw[->,ultra thick] (12.5,0)--(12.5,-0.5);
\draw[ultra thick] (12.5,0)--(13.5,0);

\end{scope}
\end{tikzpicture}
\caption{Shifting two rays}\label{shift move}
\end{figure}

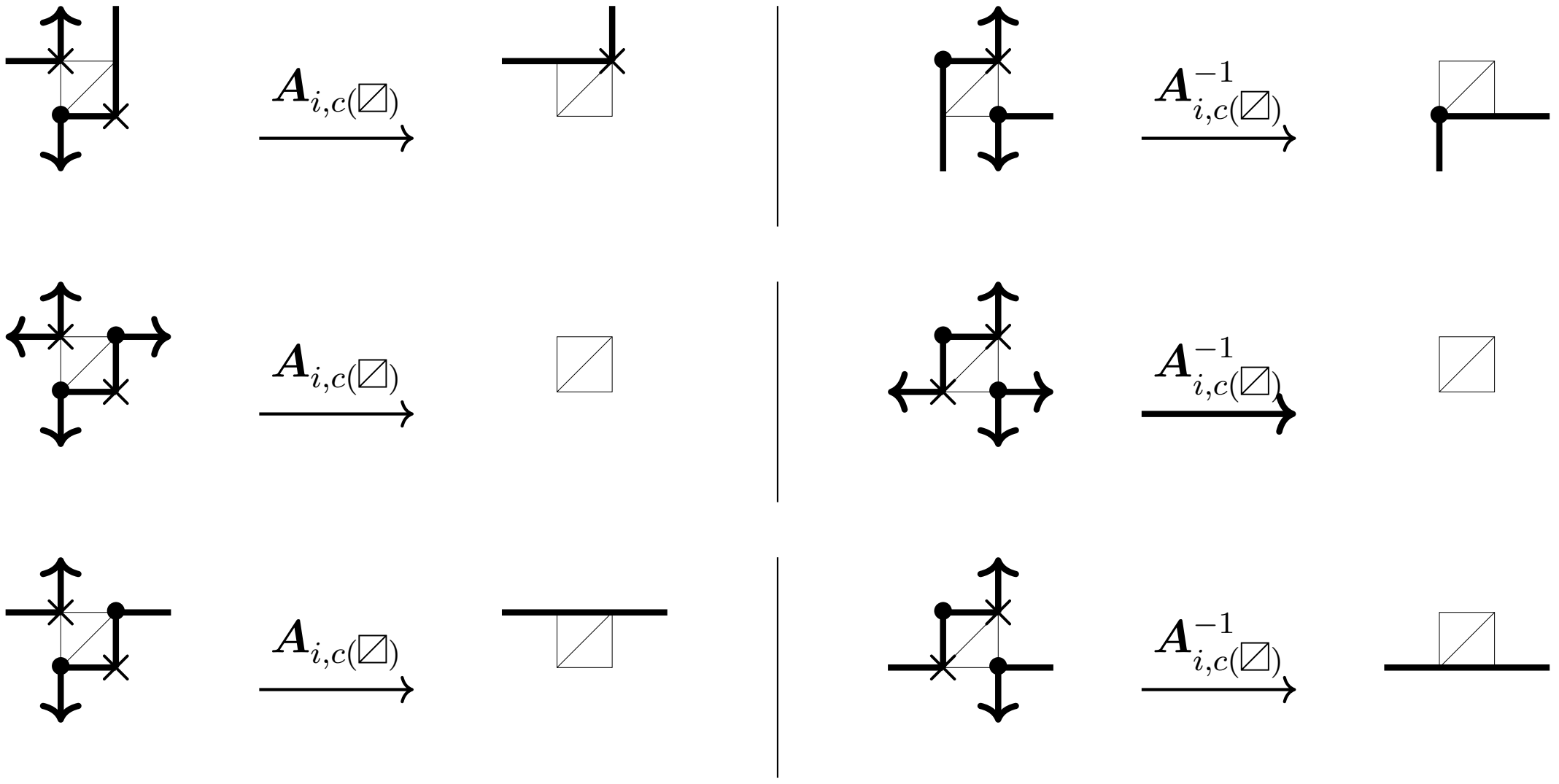
\begin{figure}[H]   
\centering
\begin{tikzpicture}

\begin{scope}[shift={(0,0)}]

\thindiagonalbox at (0,0);

\node at (0,0.5) {$\bm\times$};
\node at (0,0) {\textbullet};
\node at (0.5,0) {$\bm\times$};

\draw[->,ultra thick] (0,0)--(0,-0.5);
\draw[ultra thick] (0,0)--(0.5,0);
\draw[ultra thick] (0.5,0)--(0.5,1);
\draw[ultra thick] (0,0.5)--(-0.5,0.5);
\draw[->,ultra thick] (0,0.5)--(0,1);

\node at (2.5,0.2) {$\bm A_{i,c(\boxxslash)}$};
\draw[->,thick] (1.8,-0.2)--(3.2,-0.2);

\thindiagonalbox at (4.5,0);

\node at (5,0.5) {$\bm\times$};

\draw[ultra thick] (4,0.5)--(5,0.5);
\draw[ultra thick] (5,0.5)--(5,1);

\draw (6.5,-1)--(6.5,1);

\thindiagonalbox at (8,0);

\node at (8.5,0.5) {$\bm\times$};
\node at (8,0.5) {\textbullet};
\node at (8.5,0) {\textbullet};

\draw[ultra thick] (8.5,0.5)--(8,0.5);
\draw[ultra thick] (8,0.5)--(8,-0.5);
\draw[->, ultra thick] (8.5,0.5)--(8.5,1);
\draw[->, ultra thick] (8.5,0)--(8.5,-0.5);
\draw[ultra thick] (8.5,0)--(9,0);

\node at (10.5,0.2) {$\bm A^{-1}_{i,c(\boxxslash)}$};
\draw[->,thick] (9.8,-0.2)--(11.2,-0.2);

\thindiagonalbox at (12.5,0);

\node at (12.5,0) {\textbullet};

\draw[ultra thick] (12.5,0)--(12.5,-0.5);
\draw[ultra thick] (12.5,0)--(13.5,0);

\end{scope}

\begin{scope}[shift={(0,-2.5)}]

\thindiagonalbox at (0,0);

\node at (0,0.5) {$\bm\times$};
\node at (0,0) {\textbullet};
\node at (0.5,0) {$\bm\times$};
\node at (0.5,0.5) {\textbullet};

\draw[->, ultra thick] (0,0)--(0,-0.5);
\draw[ultra thick] (0,0)--(0.5,0);
\draw[ultra thick] (0.5,0)--(0.5,0.5);
\draw[->,ultra thick] (0.5,0.5)--(1,0.5);
\draw[->,ultra thick] (0,0.5)--(-0.5,0.5);
\draw[->, ultra thick] (0,0.5)--(0,1);

\node at (2.5,0.2) {$\bm A_{i,c(\boxxslash)}$};
\draw[->,thick] (1.8,-0.2)--(3.2,-0.2);

\thindiagonalbox at (4.5,0);

\draw (6.5,-1)--(6.5,1);

\thindiagonalbox at (8,0);

\node at (8.5,0.5) {$\bm\times$};
\node at (8,0.5) {\textbullet};
\node at (8.5,0) {\textbullet};
\node at (8,0) {$\bm \times$};

\draw[->,ultra thick] (8,0)--(7.5,0);
\draw[ultra thick] (8.5,0.5)--(8,0.5);
\draw[ ultra thick] (8,0.5)--(8,0);
\draw[->, ultra thick] (8.5,0.5)--(8.5,1);
\draw[->, ultra thick] (8.5,0)--(8.5,-0.5);
\draw[->,ultra thick] (8.5,0)--(9,0);

\node at (10.5,0.2) {$\bm A^{-1}_{i,c(\boxxslash)}$};
\draw[->, ultra thick] (9.8,-0.2)--(11.2,-0.2);

\thindiagonalbox at (12.5,0);

\end{scope}

\begin{scope}[shift={(0,-5)}]

\thindiagonalbox at (0,0);

\node at (0,0.5) {$\bm\times$};
\node at (0,0) {\textbullet};
\node at (0.5,0) {$\bm\times$};
\node at (0.5,0.5) {\textbullet};

\draw[->, ultra thick] (0,0)--(0,-0.5);
\draw[ultra thick] (0,0)--(0.5,0);
\draw[ultra thick] (0.5,0)--(0.5,0.5);
\draw[ultra thick] (0.5,0.5)--(1,0.5);
\draw[ultra thick] (0,0.5)--(-0.5,0.5);
\draw[->, ultra thick] (0,0.5)--(0,1);

\node at (2.5,0.2) {$\bm A_{i,c(\boxxslash)}$};
\draw[->,thick] (1.8,-0.2)--(3.2,-0.2);

\thindiagonalbox at (4.5,0);
\draw[ultra thick] (4,0.5)--(5.5,0.5);

\draw (6.5,-1)--(6.5,1);

\thindiagonalbox at (8,0);

\node at (8.5,0.5) {$\bm\times$};
\node at (8,0.5) {\textbullet};
\node at (8.5,0) {\textbullet};
\node at (8,0) {$\bm \times$};

\draw[ultra thick] (8,0)--(7.5,0);
\draw[ultra thick] (8.5,0.5)--(8,0.5);
\draw[ ultra thick] (8,0.5)--(8,0);
\draw[->, ultra  thick] (8.5,0.5)--(8.5,1);
\draw[->, ultra thick] (8.5,0)--(8.5,-0.5);
\draw[ultra thick] (8.5,0)--(9,0);

\node at (10.5,0.2) {$\bm A^{-1}_{i,c(\boxxslash)}$};
\draw[->,  thick] (9.8,-0.2)--(11.2,-0.2);

\thindiagonalbox at (12.5,0);

\draw[ultra thick] (12,0)--(13.5,0);

\end{scope}
\end{tikzpicture}
\caption{Gluing two chains.}\label{glue move}
\end{figure}

The operations of shifting and gluing are illustrated in Figures \ref{shift move} and \ref{glue move} in the case when the two rays are on the same vertical line. In these figures, the box $\boxxslash$  is a linking box.

\begin{lem}
Let $V$ and $W$ be two chain modules such that the number of chains in $V$ is not less than the number of chains in $W$. Then, $V$ and $W$ are linked if and only if $W$ is obtained from $V$
by either shifting two rays or by gluing two chains.
\end{lem}
\begin{proof}
The only if direction follows from Lemma \ref{must be linking box}  and Figures \ref{shift move}, \ref{glue move}.

For the if direction, it suffices to consider a module $V$ with two chains. For such a module, we clearly have linking boxes which link $V$ to $W$ both by shifting rays and gluing chains, see Figures \ref{shift move}, \ref{glue move}.
\end{proof}

\subsection{Linear independence of \texorpdfstring{$q$-characters}{}  of thin irreducible modules}

\label{linear independece q}

In this section, we show that if $U$ is a thin module whose $q$-character is a non-trivial sum of $q$-characters of chain modules, then $U$ is reducible.

It is well-known that the $q$-character is an injective map defined on the Grothendieck ring of finite-dimensional representations of a quantum affine algebra. In particular, the $q$-character of an irreducible module cannot be equal to a sum of several $q$-characters of irreducible modules. 
However, in this text, we are working with quantum toroidal algebras and modules which are not necessarily highest weight, so we have to come up with an independent proof.

Next, we study two-term extensions of chain modules.

Let $V$ and $W$ be two distinct chain modules. Let $U$ be a thin module such that $\chi_q(V)+\chi_q(W)\subset \chi_q(U)$.
Whenever we have a chain module $X$ such that $\chi_q(X)\subset \chi_q(U)$, we say that $X\subset U$ is the subspace of $\ell$-weight vectors corresponding to $\chi_q(X)$. 

Let $\bs m\in \chi_q(V)$ and $\bs m'\in \chi_q(W)$ be such that  $\bs m'=\bm A^{-1}_{i,c}\bs m$. Then, the matrix coefficients in $U$ have the form
\begin{align} \label{structure constants} \bra{\bs m'} F_i(z)\ket{\bm m} =c_{\bs m',\bs m}\delta(z/c),\qquad
\bra{\bs m} E_i(z)\ket{\bm m'} =b_{\bs m,\bs m'}\delta(z/c),
\end{align}
where  $c_{\bs m',\bs m},b_{\bs m,\bs m'}$ are constants. 

Note that if $V$ and $W$ are not linked, then all matrix coefficients between them are zero.

The next lemma says that if there exists a non-zero matrix coefficient between a vector corresponding to $V$ and a vector corresponding to $W$, then all such possible matrix coefficients are non-zero while the opposite matrix coefficients are all zero.

\begin{lem}\label{gluing lemma}
Let  $W$ be obtained from $V$ either by gluing two chains or shifting two rays. Assume that there exist $\bs m\in \chi_q(V)$ and $\bs m'\in \chi_q(W)$ such that  either
$c_{\bs m',\bs m}\neq 0$ or $b_{\bs m',\bs m}\neq 0$.

Then, for any monomials $\bs t\in \chi_q(V)$,  $\bs t'\in \chi_q(W)$ such that $\bs t'=\bm A^{-1}_{j,b}\bs t$, we have  $c_{\bs t',\bs t}\neq 0$ and $b_{\bs t,\bs t'}=0$. 
Similarly, for any monomials $\bs t\in \chi_q(V)$,  $\bs t'\in \chi_q(W)$ such that $\bs t'= \bm A_{j,b}\bs t$, we have  $b_{\bs t',\bs t}\neq 0$ and $c_{\bs t,\bs t'}=0$.
\end{lem}
\begin{proof}
We repeatedly use the following three observations.

\begin{enumerate}[leftmargin=*]
\item Let $\bs m',\bs m$ be two monomials in $\chi_q(U)$ such that $$\bs m'/\bs m\not\in\{1,\  \bm A_{i,c}^{\pm1} \bm A_{i,cq_2}^{\pm1},\ \bm  A_{i,c}^{\pm1} \bm A_{i+1,cq_1}^{\pm1}, \  \bm A_{i,c}^{\pm1} \bm A_{i+1,cq_3}^{\pm1}\}_{c\in\C^\times}.$$ 
Let $B_1(z), B_2(w)$ be  currents of the form $F_i(z)$ or $E_j(w)$. Then, from the quadratic relations, $\bra{\bs m'} B_1(z) B_2(w) \ket{\bs m} \ \neq 0$ if and only if $\bra{\bs m'} B_2(w) B_1(z) \ket{\bs m}\neq 0$.  Moreover, if $B_1(z)$ and $B_2(w)$ are not the same currents, then
there exists a unique $\bs m''$  such that  $$\bra{\bs m'} B_1(z) B_2(w) \ket{\bs m}= \bra{\bs m'} B_1(z)\ket{\bs m''}\bra{\bs m''} B_2(w) \ket{\bs m}.$$

\item If $\bs m',\bs m$ are in the $q$-character of a chain module and $\bs m'= \bm A^{-1}_{i,c}\bs m$, then $b_{\bs m,\bs m'}\ne 0$ and $c_{\bs m',\bs m}\ne 0$.

\item If $\bs m',\bs m$ are in the $q$-characters of two different chain modules
and $\bs m'= \bm A^{-1}_{i,c}\bs m$, then $b_{\bs m,\bs m'}\,c_{\bs m',\bs m}= 0$.

\end{enumerate}

\medskip

Let $W$ be obtained from $V$ by gluing two chains, both with a ray on the line $x=a+\frac{1}{2}$. We may assume that all possible linking boxes are of the form $(a,b)$ and $(a+1,b)$. The boxes $(a,b)$ correspond to the right pictures in Figure \ref{glue move}, while $(a+1,b)$ correspond to the left pictures. 

Assume that the matrix coefficient of type $c$ corresponding to a linking box $(a,b)$ is non-zero. We first show that the matrix coefficient of type $b$ corresponding to a linking box $(a+1,b)$ is non-zero, and then that the matrix coefficient of type $c$ corresponding to a box $(a,b-1)$ is non-zero.

Indeed, without loss of generality we can assume that $(a+\frac{3}{2},b+\frac{1}{2})$ is a cross. Then,
by the assumption we have 
$$
    \bra{\bm A_{i+1,q_1^{-1} c }\bs m} E_{i}(z_3) \ket{\bm A_{i+1,q_1^{-1} c}\bs m'} \    \bra{\bm A_{i+1,q_1^{-1} c }\bs m'} E_{i+1}(z_2) \ket{\bs m'}\ \bra{\bs m'} F_i(z_1)\ket{\bs m}\neq 0.
$$
Then, from the $E$-$F$ relations we have
$$
\bra{\bm A_{i+1,q_1^{-1} c }\bs m} E_{i}(z_3) \ket{\bm A_{i+1,q_1^{-1} c}\bs m'} \    \bra{\bm A_{i+1,q_1^{-1} c}\bs m'} F_i(z_1)\ket{\bs A_{i+1,q_1^{-1} c}\bs m}\ \bra{\bm A_{i+1,q_1^{-1} c }\bs m} E_{i+1}(z_2) \ket{\bs m}\neq 0,
$$
which implies that $\bra{\bm A_{i+1,q_1^{-1} c }\bs m} E_{i+1}(z_2) \ket{\bs m}\neq 0$.

Similarly, 
$$
\bra{ \bm A_{i+1,q_1^{-1} c }\bs m} E_{i+1}(z_3) \ket{\bs m }\ \bra{ \bs m} F_{i+1}(z_2)  \ket{\bm A_{i+1,q_3 c} \bs m} \ \bra{\bm A_{i+1,q_3 c }\bs m} E_{i+1}(z_1) \ket{\bs m}\neq 0
$$
implies
$$
 \bra{ \bm A_{i+1,q_1^{-1} c }\bs m} F_{i+1}(z_2)  \ket{ \bs m''} \
\bra{ \bs m''} E_{i+1}(z_3) \ket{\bm A_{i+1,q_3 c }\bs m }\
\ \bra{\bm A_{i+1,q_3 c }\bs m} E_{i+1}(z_1) \ket{\bs m}\neq 0,
$$
where $\bs m''= \bm A_{i+1,q_1^{-1}c}\bm A_{i+1,q_3c}\bs m$. Therefore, $\bra{  \bm A_{i+1,q_3 c}^{-1}\bs m''} F_{i+1}(z_2)  \ket{ \bs m''}\neq 0$.

All other cases are handled similarly.

We note that if $W$ is obtained by shifting two rays, then some intermediate vectors correspond to monomials in a module $W'$ which has one more chain than $V$.
\end{proof}

The two-term extension of the chain modules can be classified by studying the proof of Lemma \ref{gluing lemma}.  All such modules are obtained by taking a tensor product of a chain module with an appropriate single chain module.

We do not use it, so we simply formulate the result and omit further details.

\begin{prop} Let $V$ and $W$ be two non-trivial distinct chain modules. Then,  there exists an indecomposable module $U$ such that we have a short exact sequence of modules $0\to V\to U\to W\to 0$ if and only if one of the modules $V$ or $W$ is obtained from the other by gluing two chains. In this case, the module $U$ is unique up to isomorphism. \qed 
\end{prop}

Here is the main result of this section.

\begin{thm} \label{U is reducible} Let $U$ be a thin locally finite module such that there exists a chain module $V$ with $\chi_q(V)\subset \chi_q(U)$ and $\chi_q(V)\neq \chi_q(U)$. Then,  $U$ is reducible.
\end{thm}
\begin{proof} 
If $V\subset U$ is not a submodule, we have at least one other chain module $W$ such that $\chi_q(W)\subset \chi_q(U)$ and a matrix coefficient between $\bs m_1\in\chi_q(V)$ and  $\bs m_2\in\chi_q(W)$ is non-zero. Then, we claim that the submodule of $U$ generated by $\ket{\bs m_2}$ does not contain $V$ and therefore is proper.

For that, it is sufficient to show that there is no sequence of monomials $\bs m_1,\bm m_2,\dots,\bs m_{k+1}$ in $\chi_q(U)$ such that $\bs m_{k+1}=\bs m_1$ and for $i=1,\dots,k$, a matrix coefficient between $\ket{\bs m_i}$ and $\ket{\bs m_{i+1}}$ is non-zero. Let $V_i$ be the chain module containing $\bs m_i$. 
Note that by Lemma \ref{gluing lemma} we have at least three different chain modules $V_i$.
 
For each $i$, we have $\bs m_{i+1}/\bs m_i=\bm A^{\pm 1}_{j,c}$ for some $j,c$ and a choice of $\pm$. Moreover, there exists $\tilde i$ such that $(\bs m_{i+1}/\bs m_i)( \bs m_{\tilde i+1}/\bs m_{\tilde i})=1$.

Suppose we have such a sequence with $k$  the least possible.  Then $k$ is an even number and $k>2$.

 Note that we cannot have an $i$ such that  $\bs m_{i+2}=\bs m_i$. Indeed, if $V_i=V_{i+1}$, then we remove $\bs m_i$ and $\bs m_{i+1}$ and get a shorter sequence. And if $V_i\neq V_{i+1}$ then the product of matrix coefficients is zero.

For each $i$ we have $\bra{\bs m_{i+1}} B_i(z) \ket{\bs m_i}\neq 0$ where $B_i(z)$ is either $E_j(z)$ or $F_j(z)$ for some $j=j(i)$. Assume that $B_1$ is of type $E$.
Then similarly to the proof of Lemma \ref{gluing lemma}, we move all $B_i(z)$ which correspond to $E$ operators to the left of those which correspond to $F$ operators. In the process, $\bs m_i$ and $V_i$ with $i>2$ can change. The number $k$ remains the same but the number of different $V_i$ can change, but that number cannot become smaller than 3 as the number 2 is impossible. 

Then we obtain a new sequence of monomials such that $B_1,\dots, B_{k/2}$ correspond to $E$ operators  and  $B_{k/2+1},\dots, B_{k}$ correspond to $F$ operators. Then we move $B_{k/2+1}$ to the left towards the corresponding $E$ operator at the position $j$ such that $\tilde j=k/2+1$. Then we necessarily arrive at the situation when $\bs m_i=\bs m_{i+2}$ for some $i$ which is a contradiction.

If $B_1$ is of type $F$, the proof is similar by moving $F_i$ to the left of $E_j$.
\end{proof}

This result is readily extended to the tensor product case.
\begin{cor} \label{contain tensor product}Let $ (M_j)_{j=1}^m$ be a collection of $\mathcal E_n$ chain modules and let $(a_j)_{j=1}^m$ be $n$-independent over $q_1,q_2$. Let $U$ be a locally finite thin $\mathcal E_n$-module such that $\prod_{j=1}^m \chi_q(\sigma_{a_j}(M_j))\subset \chi_q(U)$  and 
 $\prod_{j=1}^m \chi_q(\sigma_{a_j}(M_j))\neq  \chi_q(U).$  Then, $U$ is reducible.
\end{cor}
\begin{proof}
The proof of Corollary \ref{contain tensor product} is parallel to the proof of Theorem \ref{U is reducible}. 
\end{proof}

\subsection{Uniqueness of a module with a given \texorpdfstring{$q$}{}-character}

\label{square complete}

In this section, we show that any two irreducible, thin and locally finite $\mathcal E_n$-modules with the same $q$-character are isomorphic.

Again, in the case of finite-dimensional representations of a quantum affine algebra, such a statement follows from the injectivity of $q$-character map.

We deduce this statement from the properties of the graphs of modules.

Let $V$ be an irreducible thin locally finite $\E_n$-module. Let $G_V$ be the graph whose vertices are the monomials in $\chi_q(V),$ with two vertices $\bs m$, $\bs m'$  connected by an edge if and only if  $\bm m\in e_i(\bs m')\cup f_i(\bs m')$ for some $i\in \hat I.$

Let $G$ be a graph with edge set $E(G)$. We call a subset $S\subset E(G)$ \textit{square complete} if there are no squares (4-cycles) in $G$ with exactly three edges in $S$.

%We call a graph $G$ \textit{square generated} if there exists a spanning tree $T$ of $G$ such that the only square complete set of edges $S$ containing the edge set of $T$ is the set of all edges of $G$.

For a subset $S_0$ of edges in $G$, we recursively define subsets $S_i$, $i\in\Z_{> 0},$ of edges of $G$ as follows. Let $S_i$ be the set of all edges of $S_{i-1}$ together with all edges which belong to a square in $G$ with other three edges in $S_{i-1}$. We have $S_i\subset S_{i+1}$. Then, the set $\langle S_0\rangle=\cup_{i=0}^\infty S_i$ is square complete and we call it the set square generated by $S_0$. 

A graph $G$ is called \textit{square generated} if and only if there exists a spanning tree $T$ in $G$ such that $\langle E(T)\rangle=E(G)$.

\begin{prop}  \label{cartesian product} Let $G = G_1\times G_2$ be the Cartesian product (box product) of the graphs $G_1$ and $G_2$. If $G_1$ and $G_2$ are square generated, then $G$ is also square generated.
\end{prop}

\begin{proof}

 For an edge $\beta = (v,v') \in E(G_2)$ and $u\in G_1$, we denote by $u\times \beta \in E(G)$ the edge between $(u,v)$ and $(u,v').$ Similarly, for an edge $\alpha = (u,u')\in E(G_1)$ and $v\in G_2,$ we denote by $\alpha\times v \in E(G)$ the edge between $(u,v)$ and $(u',v).$
By definition, all edges in $G$ are either of form $u\times \beta$ or $\alpha\times v$ as above.

Let $T_1\subset G_1$ and $T_2\subset G_2$ be the corresponding spanning trees such that $\langle E(T_1) \rangle = E(G_1)$ and $\langle E(T_2)\rangle = E(G_2).$ We define a spanning tree $T\subset G$ as follows. Fix $v_0\in G_2$. The edges of $T$ are all edges of the form $u\times \beta$, where $u\in G_1$, $\beta \in E(T_2)$, and all edges of the form $\alpha \times v_0$, where $\alpha \in E(T_1).$

The graph $G$ has squares of three types: type $G_1$-squares where the second coordinates of all 4 vertices are the same; 
type $G_2$-squares where  the first coordinates of all 4 vertices are the same; and of mixed type with vertices of the form $(u_1,v_1)$,$ (u_2,v_2)$, $(u_1,v_2)$, $(u_2,v_1),$ where $u_i\in G_1$,  $v_i\in G_2$.

Since $G_2$ is square generated, from type $G_2$-squares we get that $u\times E(G_2)\subset \langle E(T)\rangle $ for every $u\in G_1.$ Similarly, since $G_1$ is square generated, from type $G_1$-squares we get that $E(G_1)\times v_0\subset \langle E(T)\rangle .$ Lastly, from mixed squares we get that $E(G_1)\times v \subset \langle E(T)\rangle $ for every $v\neq v_0 $ in $G_2$.

It follows that $E(G) \subset \langle E(T)\rangle.$
\end{proof}

\begin{prop} \label{graph square iso} Let $V_1$ and $V_2$ be  irreducible, thin and locally finite $\E_n$-modules such that $\chi_q(V_1)=\chi_q(V_2)$. If their graph $G_{V_1}=G_{V_2}$ is square generated, then $V_1\cong V_2$ as $\mathcal E_n$-modules.
\end{prop}
\begin{proof}
Let $V$ be a module such that $\chi_q(V)=\chi_q(V_1)=\chi_q(V_2)$. 
Choose an $\ell$-weighted basis $\{v_{\bs{m}}\}$ of $V$ labeled by the vertices $\bm m$ of $G_V$. We have an edge between two vertices $\bm m, \bm m'$ if and only if the matrix coefficients of $E_i(z)$ or $F_i(z)$ are non-zero.
Suppose that $\bm m\in e_i(\bs m')$.
Then, the $(\bs m',\bs m)$ matrix coefficient of $F_i(z)$ has the form $a_{\bs m,\bs m'}\delta(z/c),$ where $a_{\bm m, \bm m'}\in \mathbb C^\times$ and $c$ is the pole of the eigenfunction of $K_i(z)$. At the same time, 
the $(\bs m,\bs m')$ matrix coefficient of $E_i(z)$ has the form $b_{\bs m',\bs m}\delta(z/c)$ for some $b_{\bm m', \bm m}\in \mathbb C^\times$. Therefore, $V$ is uniquely determined by the constants $a_{\bs m,\bs m'}$ and $b_{\bs m',\bs m}$. Such constants were called supplements in \cite{FJM}. We show that the supplement is unique up to renormalization of basis vectors.

Note that for thin modules, by relation \eqref{EFK}, we have $$a_{\bs m,\bs m'}b_{\bs m',\bs m}=(q-q^{-1})^{-1}\Res_{z=c} m_i^+(z)\frac{dz}{z}.$$ Therefore, if $a_{\bs m,\bs m'}$ is determined, then $b_{\bm m', \bm m}$ is also determined and vice-versa. 

Let $T\subset G_V$ be the spanning tree such that $\langle E(T) \rangle =  E(G).$ By renormalizing the basis vectors, we can make $a_{\bs m,\bs m'}=1$ whenever $T$ contains the edge $(\bs m,\bs m')$ and $\bs m\in e_i(\bs m')$. 

By the quadratic relations \eqref{relations 2.5}, \eqref{relations 2.6},
\eqref{relations 2.7} of $\mathcal E_n$, if the constants $a_{\bs m, \bs m'}$ and $b_{\bs m', \bs m}$ are determined for three edges of a square, then they uniquely determine the constants for the remaining edge.

Let $S_0=E(T)$ be the set of all edges of $T$. Then, we recursively determine the constants $a_{\bs m, \bs m'}$ and $b_{\bs m', \bs m}$ for edges in sets $S_i$.

It follows that all constants $a_{\bs m, \bs m'}$ and $b_{\bs m', \bs m}$ in $V$ are uniquely determined.
\end{proof}

Now, we show that the graphs of chain modules are square generated. We work with the case of $H_{\lambda, \mu, \nu}^{\tilde \lambda, \tilde \mu, \tilde \nu}$ only. The case of $V_{\lambda, \mu, \nu}^{\tilde \lambda, \tilde \mu, \tilde \nu}$ is similar. We start by defining a spanning tree. For brevity, denote by $\mathcal H = \mathcal H_{\lambda, \mu, \nu}^{\tilde \lambda, \tilde \mu, \tilde \nu}$ the set of monomials in $\chi_q(H_{\lambda, \mu, \nu}^{\tilde \lambda, \tilde \mu, \tilde \nu})$. Denote by $G$ the graph of $H_{\lambda, \mu, \nu}^{\tilde \lambda, \tilde \mu, \tilde \nu}$.

Fix $\bm m_0\in \mathcal H$. Our construction of the spanning tree $T$ depends on the choice of $\bs m_0,$  but we do not carry this dependence on $\bs m_0$ in our notation.

Let $ C,  C_0$ be single chains with the same rays.  We say that a bullet of $ C$ is \textit{good} for $ C_0$ if it lies above $ C_0$. Similarly, we say that a cross of $C$ is \textit{good} for $ C_0$ if it lies below $C_0$.

Let $\bm m \in \mathcal H$. We say that a bullet (resp. cross) on a chain $ C$ of $\bs m$ is \textit{good} if it is good for the chain $ C_0$ of $\bs m_0$ with the same rays as  $C$.

Informally speaking, expanding $\bs m$ at a good cross or a good bullet produces a monomial that is ``closer" to our fixed monomial $\bs m_0$.

Let $\mathcal A_{\bm m}$ be the set of all crosses and bullets of $\bm m$ that are good and expandable.

\begin{lem} \label{good corner} Let $\bs m \in\mathcal H, \bm m\neq \bm m_0.$  Then, the set $\mathcal A_{\bm m}$ is nonempty.
\end{lem}
\begin{proof}
The lemma is clearly true for the case when $\bs m$ and $\bs m_0$ are single chain monomials.

Note that if $(u,v)$ is a good bullet and $(u+1,v+1)$ is a bullet, then $(u+1,v+1)$ is also good by the non-intersecting property of chains. Similarly, if $(u,v)$ is a good cross and $(u-1,v-1)$ is a cross, then $(u-1,v-1)$ is also good. The general case follows.
\end{proof}

We order the chains of a monomial in $\mathcal H$ so that larger chains are below smaller ones. We order points $(u,v)\in(\Z+\frac{1}{2})^2$ lexicographically.

We define a total order in $\mathcal A_{\bm m}$ by setting $c\prec c'$ if either $c$ is in a smaller chain than $c',$ or $c,c'$ are in the same chain and $c$ precedes $c'$ in the lexicographic order of $(\Z+\frac{1}{2})^2$.

Let $T$ be the graph with the same vertex set $\mathcal H$ of $G$, with an edge between $\bm m, \bm m'$ if $\bm m'$ is obtained from $\bm m$ by expanding the maximal element of $\mathcal A_{\bm m}.$

\begin{prop}  \label{square generated}
The graph $T\subset G$ is a spanning tree and  $G$ is square generated.
\end{prop}
\begin{proof}
For $\bs m\in\mathcal H$, define depth $d(\bs m)\in\Z_{\geq 0}$ as follows. If $\bs m$ and $\bs m_0$ are single chains with the same rays, then $d(\bm  m)$ is the area of all finite open connected regions enclosed by $\bs m$ and $\bs m_0$. In the general case of multiple chains, we define $d(\bs m)$ as the sum of depths of all chains of $\bs m$ in relation to the chains of $\bs m_0$ with the same rays.

We have $d(\bs m)=0$ if and only if $\bs m=\bs m_0$. Moreover, expansion of a good corner decreases depth by one. Therefore, by Lemma \ref{good corner}, the graph $T$ is connected.

For any edge $(\bs m, \bs m')$ in $G$, we have  $|d(\bs m)-d(\bs m')|=1$.
We orient an edge $(\bs m, \bs m')$ from $\bm m$ to $\bm m'$ if $d(\bs m)=d(\bs m')-1$.

Then, we clearly have no oriented cycles in $G$.
Moreover, by definition for any $\bs m\in\mathcal H$, $\bs m\neq \bs m_0$, there exists a unique incoming edge in $T$. It implies that $T$ has no cycles, that is $T$ is a tree. 

Finally, we show that $G$ has enough squares.

Let $S_0=E(T)$ be the set of edges of $T$. Then,
by induction on depth we show that any edge connecting two vertices of depths $d$ and $d+1$ is in $S_d$. 

If $d=0$, then every such edge is in $E(T) = S_0.$ 

Let $(\bs m,\bs m')$ be an edge in $G$ with $d=d(\bm m')-1=d(\bm m)>0$. Then, $\bs m$ is obtained from $\bs m'$ by expanding a good corner $c$. Assume that $c$ is not maximal in $\mc A_{\bs m'}$. Then, the maximal corners of $\mc A_{\bs m'}$ and $\mc A_{\bs m}$ are the same. Indeed, expanding, say, a good bullet can remove a cross but it cannot remove a good cross. Denote the maximal corner of $\mc A_{\bs m}$ by $c_0$. Let $\bs m_1, \bm m_2$ be the monomials obtained from $\bs m$ and $\bs m'$ by expanding $c_0$. Then, $d(\bs m_1)+1=d(\bs m_2)=d$ and expanding $\bs m_2$ at $c$ is $\bs m_1$. It follows that $\bs m, \bs m',\bs m_1,\bs m_2$ is a square, where the edges $(\bs m,\bs m_1)$, $(\bs m',\bs m_2)$
are in $E(T)=S_0\subset S_{d-1}$ and the edge $(\bs m_1,\bs m_2)$ is in $S_{d-1}$ by induction hypothesis. Thus, the edge $(\bs m,\bs m')$ is in $S_d$. 

It follows that $\langle S_0\rangle =E(G)$.
\end{proof}

\begin{cor}
Let $V$ be an $\mathcal E_n$-module. 
\begin{enumerate}
\item If $\chi_q(V) = \chi_q(V_{\lambda, \mu, \nu}^{\tilde \lambda, \tilde \mu, \tilde \nu})$, then $V$ is isomorphic to $V_{\lambda, \mu, \nu}^{\tilde \lambda, \tilde \mu, \tilde \nu}$ as $\mathcal E_n$-modules.

\item If $\chi_q(V) = \chi_q (H_{\lambda, \mu, \nu}^{\tilde \lambda, \tilde \mu, \tilde \nu})$, then $V$ is isomorphic to $H_{\lambda, \mu, \nu}^{\tilde \lambda, \tilde \mu, \tilde \nu}$ as $\mathcal E_n$-modules.

\end{enumerate}
\end{cor}

\begin{proof}
The corollary follows from Propositions \ref{square generated} and \ref{graph square iso}.
\end{proof}

\begin{cor} \label{tensor chain}Let $ (M_j)_{j=1}^m$ be a collection of  chain modules and let $(a_j)_{j=1}^m$ be $n$-independent over $q_1,q_2$. Let $V$ be an $\mathcal E_n$-module such that $\chi_q(V) = \prod_{j=1}^m \chi_q(\sigma_{a_j}(M_j))$. Then, $$V\cong \sigma_{a_1}(M_1)\otimes \cdots\otimes \sigma_{a_m}(M_m).$$
\end{cor}
\begin{proof}

Let $G_M$ be the graph associated to the $\mathcal E_n$-module $M = \sigma_{a_1}(M_1)\otimes \cdots\otimes \sigma_{a_m}(M_m)$. From the coproduct formula, we see that $G_M = G_1\times \cdots \times G_m,$ where $G_i$ is the graph associated to the module $\sigma_{a_i}(M_i).$ Since each $G_i$ is square generated, the claim follows from Propositions \ref{cartesian product}, \ref{graph square iso}.
\end{proof}
In particular, Corollary \ref{tensor chain} says that the tensor products of $\sigma_{a_i}(M_i)$ in different orders are isomorphic.

\subsection{The main results}  \label{main results}

In this section, we collect all the main results.

\begin{thm}   \label{irreducible is tensor product}
Let $V$ be a non-trivial irreducible, thin and locally finite $\mathcal E_n$-module. Then, $$V\cong \sigma_{a_1}(M_1)\otimes \cdots\otimes \sigma_{a_m}(M_m),$$ where $(M_i)_{i=1}^m$ are non-trivial chain modules and $(a_i)_{i=1}^m$ are $n$-independent over $q_1$, $q_2$. Any two such tensor products are isomorphic if and only if they have the same factors up to a permutation.

%in which each tensor factor is an $\mathcal E_n$-module of the form  $s_a\big(V_{\lambda, \mu, \nu}^{\tilde \lambda, \tilde \mu, \tilde \nu}\big)$ or $s_a\big(H_{\lambda, \mu, \nu}^{\tilde \lambda, \tilde \mu, \tilde \nu}\big)$, where $V_{\lambda, \mu, \nu}^{\tilde \lambda, \tilde \mu, \tilde \nu}, H_{\lambda, \mu, \nu}^{\tilde \lambda, \tilde \mu, \tilde \nu}$ are in Theorem \ref{H and V modules} and $a\in \mathbb C^\times.$
 \end{thm}
\begin{proof}
Let $\bs m\in\chi_q(V)$. By Theorem \ref{expansion is contained}, we have $\text{FM}(\bm m)\subseteq \chi_q(V).$ 

Write $\bm m = \prod_{j=1}^m\bm m_j$ as in Lemma \ref{split lemma}. Since $a_i, a_j \in \mathbb C^\times$ are $n$-independent over $q_1,q_2$ for $i\neq j$, each $\bm m_j$ is  also a seed monomial. Moreover, by Proposition \ref{independent monomials} each $\sigma_{a_j^{-1}}(\bm m_j)$ is a seed monomial of integer type and $\text{FM}(\bm m) = \prod_{j=1}^m\text{FM}(\bm m_j)$. By Theorem \ref{expansions of seed },  
$\text{FM}(\bm m_j)=\chi_q(\sigma_{a_j}(M_j))$, where $M_j$ is a chain module.

Since $V$ is irreducible, $\chi_q(V)=\prod_{j=1}^m\chi_q(\sigma_{a_j}(M_j))$ follows by Corollary \ref{contain tensor product}. Then, $V=\otimes_{j=1}^m \sigma_{a_j}(M_j)$ by Corollary \ref{tensor chain}.
\end{proof}

We proceed to prime modules.
 
\begin{thm}\label{prime thm}
Let $V$ be an irreducible, thin, locally finite, and prime $\mathcal E_n$-module. 
Then, either $V=\sigma_a(V_{\lambda,\mu,\nu}^{\tilde \lambda, \tilde \mu, \tilde \nu}  )$ and 
$$\mu_i\leq \tilde \mu_{i+1}, \quad i=1,\dots,m-1, \qquad  \lambda_l\leq \tilde \mu_1\ (l,m\neq 0), \qquad  \mu_m\leq \tilde \nu_1\ (p,m\neq 0),$$
or $V=\sigma_a(H_{\lambda, \mu, \nu}^{\tilde \lambda, \tilde \mu, \tilde \nu})$ and 
$$\tilde \mu_i\leq \mu_{i+1}, \quad i=1,\dots,m-1, \qquad  \tilde  \lambda_l\leq \mu_1\ (l,m\neq 0), \qquad  \tilde \mu_m \leq \nu_1\ (p,m\neq 0).$$
\end{thm}
\begin{proof}
By Theorem \ref{irreducible is tensor product} any irreducible, thin, locally finite and prime $\mathcal E_n$-module is a shift of a single chain module.

A chain module is prime if and only if for any two neighboring chains, the condition of non-intersection is non-trivial. That corresponds to the stated inequalities.
\end{proof}

Now we classify bounded modules.

\begin{thm}
Let $V$ be an irreducible, thin, locally finite and bounded $\mathcal E_n$-module. 
Then, $$V\cong \sigma_{a_1}(M_1)\otimes \cdots\otimes \sigma_{a_m}(M_m),$$ where $(M_i)_{i=1}^m$ are chain modules of the form  $V_{\varnothing, \mu, \varnothing}^{\varnothing, \tilde \mu, \varnothing}$ or $H_{\varnothing, \mu, \varnothing}^{\varnothing, \tilde \mu, \varnothing}$ and 
$(a_i)_{i=1}^m$ are $n$-independent over $q_1$, $q_2$. Any two such tensor products are isomorphic if and only if they have the same factors up to a permutation.
\end{thm}
\begin{proof}
The theorem follows from Theorem \ref{irreducible is tensor product} and
Theorem \ref{bounded seed}.
\end{proof}

\medskip

Finally, we discuss the restriction of the bounded modules to the subalgebra $U_q\hat{\mathfrak{sl}}_n\subset \mc E_n$ generated by $ \hat{\mc U}_i$, $i=1,\dots, n-1$. 

Let $\mathcal P\subset \hat{\mathcal P}$ be the abelian subgroup generated by  $\bm i_a, i\in  I, a\in \mathbb C^\times.$
We have a group homomorphism $\rho:  \hat{\mathcal P} \to   \mathcal P$ sending 
$$
\bs{0}_a \to 1 \qquad \bs i_a\to \bs i_{a(q q_1)^{i-1}}.$$

Clearly, the restriction of an irreducible, thin, locally finite, bounded $\mc E_n$-module $V$ is a sum of thin irreducible finite dimensional $U_q\,\hat{\mathfrak{sl}}_n$-modules. The irreducible summands of the restriction of $V$ are in  bijection with dominant $U_q\,\hat{\mathfrak{sl}}_n$ monomials counted with multiplicity or, equivalently, with monomials in $\chi_q(V)$ represented by unions of chains such that all crosses have color $0$.

We show that every irreducible finite-dimensional thin $U_q\,\hat{\mathfrak{sl}}_n$-module is found inside  a thin, locally finite, and bounded $\mc E_n$-module. We restrict ourselves to the integral case. In that case, the irreducible finite-dimensional thin $U_q\,\hat{\mathfrak{sl}}_n$-modules correspond to skew Young diagrams, see \cite{NT}. Such modules were called snake modules in \cite{MY}.

Let $\eta\subset \gamma$ be partitions such that the columns of the skew Young diagram $\gamma/\eta $ are of length at most $n-1$. Then, the corresponding irreducible thin 
$U_q\,\hat{\mathfrak{sl}}_n$-module $L_{\gamma/\eta}$ has dominant monomial 
$$\bs m_{\gamma/\eta}^+=\prod_{j=1, \ \gamma_j'-\eta_j'\neq 0}^{\gamma_1} \big(\bs{\gamma_j'-\eta_j'})_{\tilde a(j)}\in \mathcal P, \qquad 
 \tilde a(j) = a q^{2(j-1)-\gamma_j'-\eta_j'},
$$
where $\gamma', \eta'$  denote the conjugate partitions and $a\in \mathbb C^\times$ is the evaluation parameter.

For the next theorem, we adopt the convention of canceling $\mu_i$ with $\tilde \mu_i$ whenever $\mu_i=\tilde \mu_i$. Namely, for $\mu_i=\tilde \mu_i$, we set
$H_{\varnothing, (\mu_1,\dots,\mu_m), \varnothing}^{\varnothing, (\tilde \mu_1,\dots,\tilde \mu_m),
\varnothing}=H_{\varnothing, (\mu_1,\dots,\mu_{i-1},\mu_{i+1},\dots,\mu_m), \varnothing}^{\varnothing, (\tilde \mu_1,\dots, \tilde \mu_{i-1},\tilde \mu_{i+1},\dots,\tilde \mu_m),
\varnothing}$. Note that the horizontal line $y=\mu_i$ can be viewed as a chain of width zero corresponding to a trivial module.

\begin{thm}
The thin irreducible $U_q\, \hat{\mathfrak{sl}}_n$-module $L_{\gamma/\eta}$  is a direct summand in the restriction of a shift of the bounded chain module $H_{\varnothing, \mu, \varnothing}^{\varnothing, \tilde \mu, \varnothing}$, where $\mu$ and $\tilde \mu$ have, respectively, the entries
$$ \gamma_1'-\eta_j'+j-\frac{3}{2} \quad  \mbox{and} \qquad \gamma_1'-\gamma_j'+j-\frac{3}{2}, \qquad j = 1,\dots, \gamma_1,$$
arranged in decreasing order.
\end{thm}
\begin{proof}
For $j=1,\dots, \gamma_1,$ set
$$ a(j) = q_2^{j-1}q_1^{\eta_j'-\gamma_1'}q_3^{\gamma_j'-\gamma_1'}, \quad b(j) = a(j)q_3^{\eta_j'-\gamma_j'}.$$
Then, the seed monomial $\bm m_{0}$ given by
\begin{equation*}
\bm m_{0} = \prod_{j=1}^{\gamma_1} \big(\bs{\gamma_j'-\eta_j'}\big)_{q^{-1}a(j)}\bs 0^{-1}_{q^{-1}b(j)} \end{equation*}
belongs to $\mathcal H_{\varnothing,\mu, \varnothing }^{\varnothing, \tilde \mu, \varnothing}$. We note that $\rho(\bs m_{0})=\bs m^+_{\gamma/\eta}$ with $a=q_1^{-1}q^{2(\gamma_1'-1)}$.
\end{proof}

\section{The language of plane partitions}  \label{language of plane partitions}
\subsection{Generalized plane partitions}\label{sec: plane partition}
In this section, we present another way to visualize  the thin monomials in the $q$-character of the chain modules $ V_{\la,\mu,\nu}^{\tilde \la,\tilde \mu,\tilde \nu}$  and $ H_{\la,\mu,\nu}^{\tilde \la,\tilde \mu,\tilde \nu}$ of Theorem \ref{H and V modules} as generalized plane partitions. The language of plane partitions is used by many authors to describe bases of thin modules, see \cite{JM}, \cite{BM}, \cite{FJMM1}, \cite{KN}.

We consider collections of unit cubes $\mbox{\mancube} = (x,y,z)$ centered at integer points in $\mathbb Z^3.$ A unit cube $\mbox{\mancube} =  (x,y,z)$ has color $i\in \hat I$ given by $i\equiv x-y,$ so that colors are constant along the $z$-direction.  A unit cube $\mbox{\mancube} = (x,y,z)$ has content $c(\mbox{\mancube}) = q_1^{-x}q_3^{-y}q_2^{-z}$.

We begin by associating single chain monomials to one layer generalized plane partitions, then proceed to the case of monomials with multiple chains. We call one-layer generalized plane partitions the \textit{generalized Young diagrams}.

To a single chain monomial we associate  the generalized Young diagram given by the collection of all unit cubes at $z=0$ whose centers are on the top of the boxes to the left and upward from the chain, see Figures \ref{vh chain plane partition}, \ref{hh plane partition}, \ref{vv chain plane partition} and \ref{hv chain plane partition}.  Then, the chain itself traces the shape of the picture as shown by the thick line.

We picture some cubes in white indicating that these cubes are present for all chains in the same module. The white cubes are present for some chain modules and absent for others.

\begin{figure}[H]
\begin{tikzpicture}

\begin{scope}[shift={(0,0)}]

\node at (4.8, -0.4) {$\cdots$};
\node at (-1.8, -0.8) {$\cdots$};
\node at (2, 0.5) {\reflectbox{$\ddots$}};

% white top first row 

\wbox at (0,0);
\wbox at (0.6,0);
\wbox at (1.2,0);
\wbox at (1.8,0);
\wbox at (2.4,0);
\wbox at (3, 0);
\wbox at (3.6, 0);

% white top first column tops

\wbox at (-0.4,-0.3);
\wbox at (-0.8,-0.6);
\wbox at (-1.2,-0.9);
\wbox at (-1.6, -1.2);
\wbox at (-2,-1.5);

% white first row front side 

\fwboxh at (2, -0.9)
\fwboxh at (2.6,-0.9);
\fwboxh at (3.2,-0.9);
\rwbox at (3.8, -0.9);

% white first column front and side

\rwboxh at (-1.4, -2.1);
\rwboxh at (-1.8, -2.4);
\fwbox at (-2.4, -2.4);

% black first row top

\bbox at (0.2,-0.3);
\bbox at (0.8,-0.3);
\bbox at (1.4,-0.3);
%\bbox at (2,-0.3);

% black first row face and side

\fbboxh at (1,-1.2);
%\fbboxh at (1.6,-1.2);
\rbboxh at (1.6, -1.2);

%black second row top

\bbox at (-0.2,-0.6);
\bbox at (0.4,-0.6);

%black  second row face and side

\fbboxh at (0,-1.5);
\rbboxh at (0.6,-1.5);

% black third row top

\bbox at (-0.6,-0.9);

% black third row face and side

\fbboxh at (-1,-1.8);
\rbboxh at (-0.4,-1.8);

\node at (-2.6, -2.5) {\reflectbox{$\ddots$}};

\end{scope}

\begin{scope}[shift={(7,-2.7)}]
\draw[step=0.5cm, gray, thin] (0,0) grid (3.5, 3);

\node at (1,1) {\textbullet};
\node at (1.5,1) {$\bm \times$};
\node at (1.5,1.5) {\textbullet};
\node at (2,1.5) {$\bm \times$};
\node at (2,2) {\textbullet};
\node at (2.5,2) {$\bm \times$};
\node at (2.5,2.5) {\textbullet};

\draw[line width=2] (1,1)--(1,0);
\draw[line width=2] (1,1)--(1.5,1);
\draw[line width=2] (1.5,1)--(1.5,1.5);
\draw[line width=2] (1.5,1.5)--(2,1.5);
\draw[line width=2] (2,1.5)--(2,2);
\draw[line width=2] (2,2)--(2.5,2);
\draw[line width=2] (2.5,2)--(2.5,2.5);
\draw[line width=2] (2.5,2.5)--(3.5,2.5);

\end{scope}

\end{tikzpicture}
\caption{The generalized Young diagram corresponding to a $vh$ chain.}
\label{vh chain plane partition}
\end{figure}
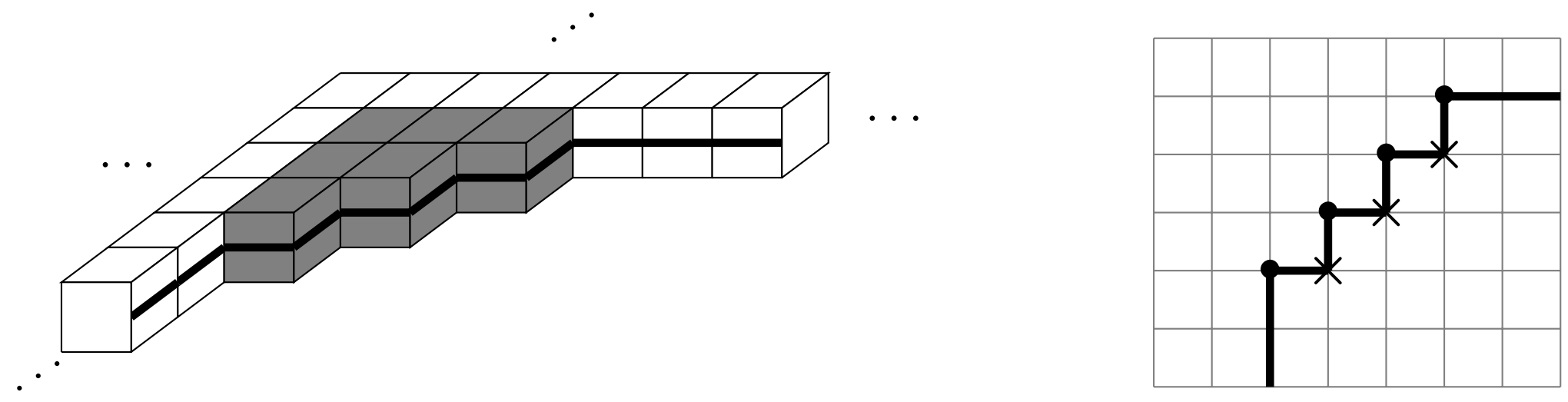

\begin{figure}[H]
\begin{tikzpicture}

\begin{scope}[shift={(0,0)}]

\node at (4, -0.4) {$\cdots$};
\node at (2, 0.5) {\reflectbox{$\ddots$}};
\node at (-1.5, -0) {$\cdots$};
\node at (-3, -1.5) {$\cdots$};
\node at (-2.3, -0.6) {$\cdots$};

% white top of first row
\wbox at (-0.6,0);
\wbox at (0,0);
\wbox at (0.6,0);
\wbox at (1.2,0);
\wbox at (1.8,0);
\wbox at (2.4,0);
\wbox at (3,0);

% white front of first row
\fwboxh at (2,-0.9);
\fwboxh at (2.6,-0.9);

% white side of first row

\rwbox at (3.2,-0.9);

% black top of first row
\bbox at (-1, -0.3);
\bbox at (-0.4,-0.3);
\bbox at (0.2,-0.3);
\bbox at (0.8,-0.3);
\bbox at (1.4,-0.3);

%black front of first row

\fbboxh at (0.4,-1.2);
\fbboxh at (1,-1.2);

% black side of first row

\rbboxh at (1.6,-1.2);

% black top of second row

\bbox at (-1.4,-0.6);
\bbox at (-0.8,-0.6);
\bbox at (-0.2,-0.6);

%black front of second row

\fbboxh at (-0.6,-1.5);

%black side of second row

\rbboxh at (-0,-1.5);

% black top of third row
\bbox at (-1.2,-0.9);
\bbox at (-1.8,-0.9);

% black sid of third row

\rbboxh at (-1,-1.8);

%black front of third row

\fbboxh at (-1.6,-1.8);
\fbboxh at (-2.2,-1.8);

\end{scope}

\begin{scope}[shift={(7,-2.7)}]
 \draw[step=0.5cm, gray, thin] (0,0) grid (3.5,3);

\node at (1,0.5) {$\bm \times$};
\node at (1,1) {\textbullet};
\node at (1.5,1) {$\bm \times$};
\node at (1.5,1.5) {\textbullet};
\node at (2.5,1.5) {$\bm \times$};
\node at (2.5,2) {\textbullet};

\draw[line width = 2] (1,0.5)--(0,0.5);
\draw[line width = 2] (1,0.5)--(1,1);
\draw[line width = 2] (1,1)--(1.5,1);
\draw[line width = 2] (1.5,1)--(1.5,1.5);
\draw[line width = 2] (1.5,1.5)--(2.5,1.5);
\draw[line width = 2] (2.5,1.5)--(2.5,2);
\draw[line width = 2] (2.5,2)--(3.5,2);

 \end{scope}

\end{tikzpicture}
\caption{The generalized Young diagram corresponding to an $hh$ chain.}
\label{hh plane partition}
\end{figure}
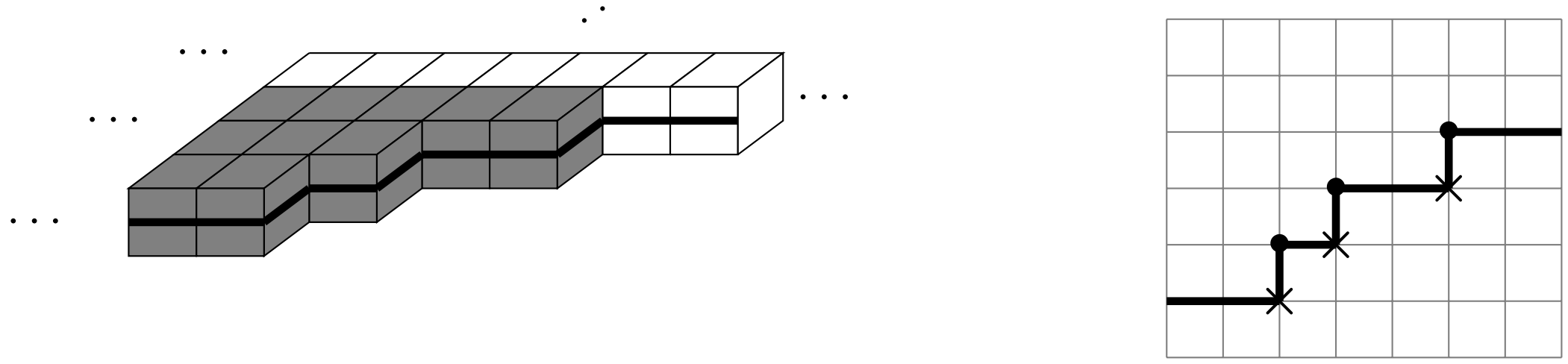

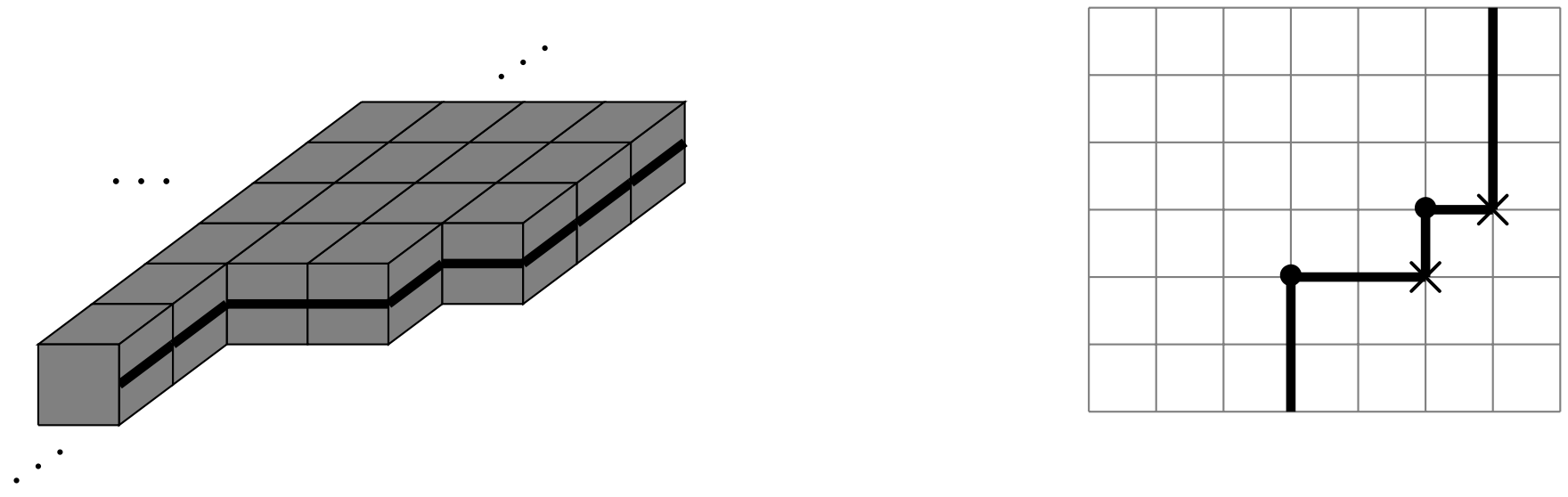
\begin{figure}[H]

\begin{tikzpicture}

\begin{scope}[shift={(2,-0.5)}]

\node at (-2, 1.5) {$\cdots$};
\node at (0.8, 2.5) {\reflectbox{$\ddots$}};

\node at (-2.8, -.5) {\reflectbox{$\ddots$}};

% black fourth row top

\bbox at (-0.4,1.2);
\bbox at (-1,1.2);
\bbox at (-1.6,1.2);

% black fourth row front

\fbboxh at (-1.4,0.3);
\fbboxh at (-0.8,0.3);

%black fourth row side

\rbboxh at (-0.2,0.3);

% black third row top

\bbox at (0.6,1.5);
\bbox at (0,1.5);
\bbox at (-0.6,1.5);
\bbox at (-1.2,1.5);

% black third row front

\fbboxh at (0.2,0.6);

% black first column  side

\rbboxh at (0.8,0.6);
\rbboxh at (1.2,0.9);
\rbboxh at (1.6,1.2);

%black fifth  row top

\bbox at (-2,0.9);

% black sixth row top

\bbox at (-2.4,0.6);

% black sixth row side

\rbboxh at (-2.2,-0.3);

% black sixth row front

\fbbox at (-2.8,-0.3);

% black fifth row side

\rbboxh at (-1.8,0);

% black second row top

\bbox at (1,1.8);
\bbox at (0.4,1.8);
\bbox at (-0.2,1.8);
\bbox at (-0.8,1.8);

% black first row top

\bbox at (1.4,2.1);
\bbox at (0.8,2.1);
\bbox at (0.2,2.1);
\bbox at (-0.4,2.1);

\end{scope}

\begin{scope}[shift={(7,-0.7)}]
 \draw[step=0.5cm, gray, thin] (0,0) grid (3.5,3);

\node at (1.5,1) {\textbullet};
\node at (2.5,1) {$\bm \times$};
\node at (2.5,1.5) {\textbullet};
\node at (3,1.5) {$\bm \times$};

\draw[line width = 2] (1.5,-0)--(1.5,1);
\draw[line width = 2] (1.5,1)--(2.5,1);
\draw[line width = 2] (2.5,1)--(2.5,1.5);
\draw[line width = 2] (2.5,1.5)--(3,1.5);
\draw[line width = 2] (3,1.5)--(3,3);

 \end{scope}

\end{tikzpicture}
\caption{The generalized Young diagram corresponding to a $vv$ chain.}
\label{vv chain plane partition}
\end{figure}

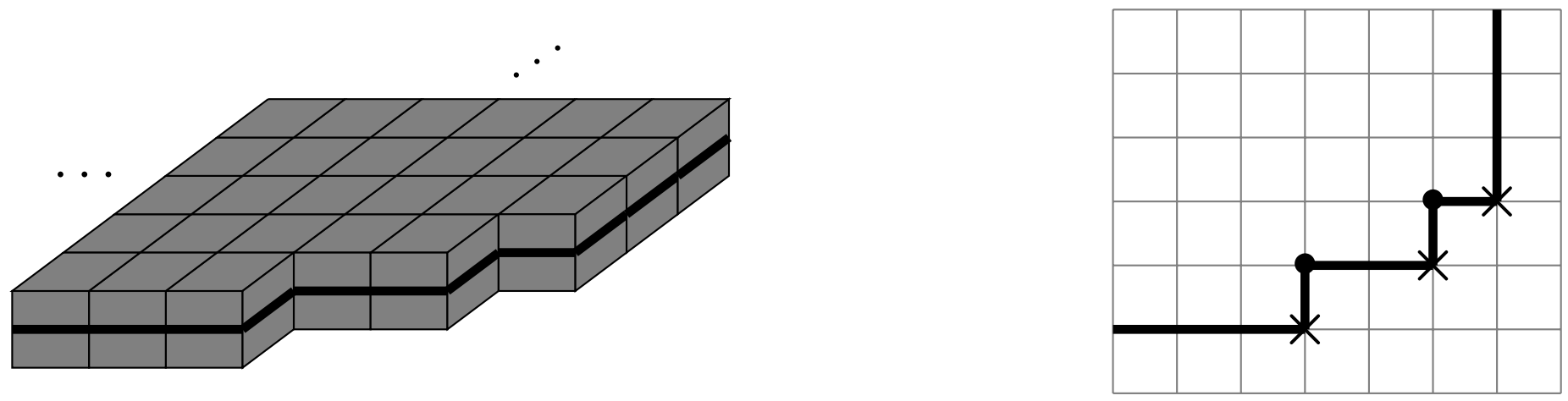
\begin{figure}[H]

\begin{tikzpicture}

\begin{scope}[shift={(2,-0.5)}]

\node at (-3, 1.5) {$\cdots$};
\node at (0.5, 2.5) {\reflectbox{$\ddots$}};

% black fourth row top

\bbox at (-0.4,1.2);
\bbox at (-1,1.2);
\bbox at (-1.6,1.2);
\bbox at (-2.2,1.2);
\bbox at (-2.8,1.2);

% black fourth row front

\fbboxh at (-1.4,0.3);
\fbboxh at (-0.8,0.3);

%black fourth row side

\rbboxh at (-0.2,0.3);

% black third row top

\bbox at (0.6,1.5);
\bbox at (0,1.5);
\bbox at (-0.6,1.5);
\bbox at (-1.2,1.5);
\bbox at (-1.8,1.5);
\bbox at (-2.4,1.5);

% black third row front

\fbboxh at (0.2,0.6);

% black first column  side

\rbboxh at (0.8,0.6);
\rbboxh at (1.2,0.9);
\rbboxh at (1.6,1.2);

%black fifth row front

\fbboxh at (-2.4,0);
\fbboxh at (-3,0);
\fbboxh at (-3.6,0);

%black fifth  row top

\bbox at (-2,0.9);
\bbox at (-2.6,0.9);
\bbox at (-3.2,0.9);

% black fifth row side

\rbboxh at (-1.8,0);

% black second row top

\bbox at (1,1.8);
\bbox at (0.4,1.8);
\bbox at (-0.2,1.8);
\bbox at (-0.8,1.8);
\bbox at (-1.4,1.8);
\bbox at (-2,1.8);

% black first row top

\bbox at (1.4,2.1);
\bbox at (0.8,2.1);
\bbox at (0.2,2.1);
\bbox at (-0.4,2.1);
\bbox at (-1,2.1);
\bbox at (-1.6,2.1);

\end{scope}

\begin{scope}[shift={(7,-0.7)}]
 \draw[step=0.5cm, gray, thin] (0,0) grid (3.5,3);

\node at (1.5,0.5) {$\bm \times$};
\node at (1.5,1) {\textbullet};
\node at (2.5,1) {$\bm \times$};
\node at (2.5,1.5) {\textbullet};
\node at (3,1.5) {$\bm \times$};

\draw[line width = 2] (0,0.5)--(1.5,0.5);
\draw[line width = 2] (1.5,0.5)--(1.5,1);
\draw[line width = 2] (1.5,1)--(2.5,1);
\draw[line width = 2] (2.5,1)--(2.5,1.5);
\draw[line width = 2] (2.5,1.5)--(3,1.5);
\draw[line width = 2] (3,1.5)--(3,3);

 \end{scope}

\end{tikzpicture}
\caption{The generalized Young diagram corresponding to an $hv$ chain.}
\label{hv chain plane partition}
\end{figure}

A generalized partition is a two-sided  stabilizing ordered sequence of integers. More precisely, a generalized partition $\bs\ga$ is a sequence $(\ga_j)_{j\in\Z}$ with $\gamma_j \in \mathbb Z\cup \{-\infty, +\infty\}, \gamma_j\geq \gamma_{j+1},$ such that there exists $\gamma^{\pm}$ for which  $\gamma_{-N} = \gamma^+$ and $\gamma_N = \gamma^-$ for large enough $N$.

Then, the usual partitions are identified with  generalized partitions such that $\ga^+ = +\infty$, $\ga^-=0$. For example, $\bs \gamma = (\cdots, +\infty, 3,2,1,0,\cdots)$ is identified with the partition $(3,2,1).$

A $vh$-chain corresponds to a generalized partition with $\ga^{-} >-\infty$, $\ga^+ = +\infty$.  An $hh$-chain corresponds to $\ga^{-} =-\infty$, $\ga^+ =+\infty$.  A $vv$-chain corresponds to $\ga^{-} >-\infty$, $\ga^+ <+\infty$. 
 An $hv$-chain corresponds to $\ga^{-} =-\infty$, $\ga^+ <+\infty$. Thus, every chain is a generalized partition. Our pictures can be thought of as Young diagrams of generalized partitions.

\medskip

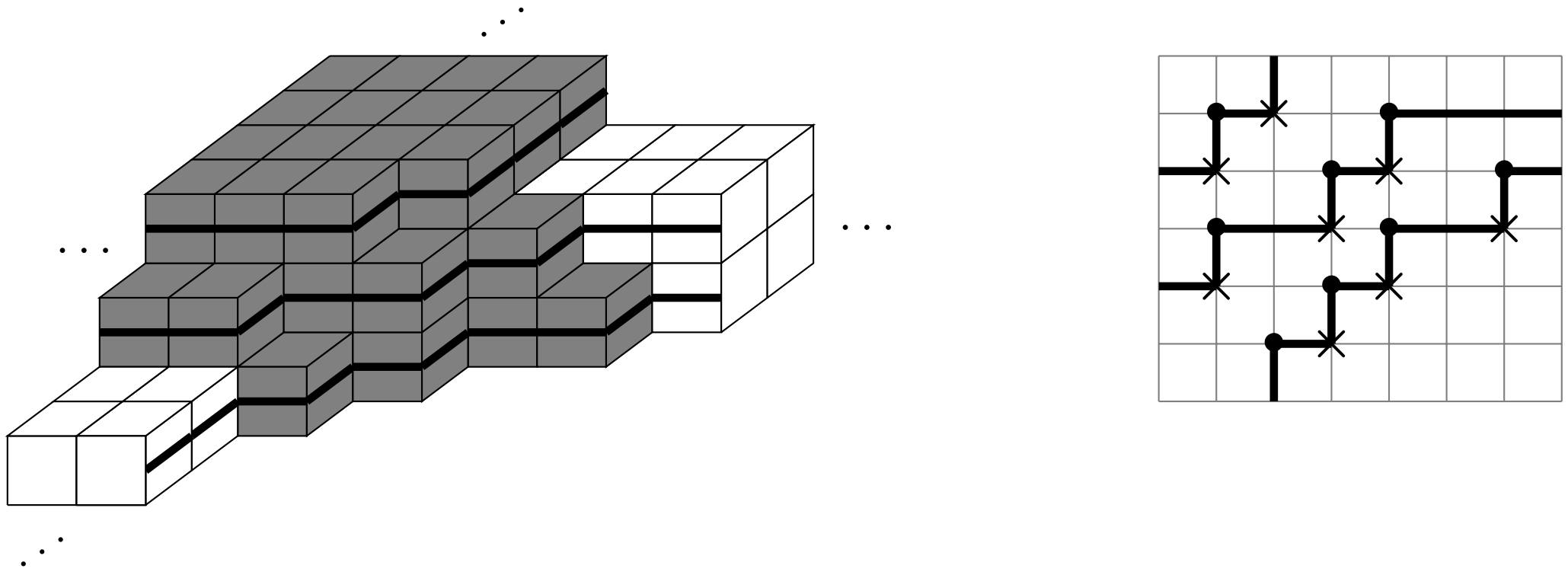
\begin{figure}[H]
\begin{tikzpicture}

\begin{scope}[shift={(0,-1)}]

\node at (1.3, 1.9) {\reflectbox{$\ddots$}};
\node at (4.5, 0) {$\cdots$};
\node at (-2.3, -0.2) {$\cdots$};
\node at (-2.7, -2.7) {\reflectbox{$\ddots$}};

\bbox at (2,-0.3);

% black first row face and side

\fbboxh at (1,-1.2);
\fbboxh at (1.6,-1.2);
\rbboxh at (2.2, -1.2);

% black third row top

\bbox at (-0.6,-0.9);

%third row face and side

\fbboxh at (-1,-1.8);
\rbboxh at (-0.4,-1.8);

% 2nd black top third row
\bbox at (-1.8,-0.3);
\bbox at (-1.2,-0.3);

% 2nd black front  third row

\fbboxh at (-2.2,-1.2);
\fbboxh at (-1.6,-1.2);

% 2nd black side third row

\rbboxh at (-1,-1.2);

% 2nd black top second row

\bbox at (0.4,0);
\bbox at (-0.2,0);
\bbox at (-0.8,0);
\bbox at (-1.4,0);

% 2nd black front second row
\fbboxh at (0,-1.5); % black first floor
\rbboxh at (0.6,-1.5); % black first floor
\fbboxh at (-0.6,-0.9);
\fbboxh at (0,-0.9);

% 2nd black side second row

\rbboxh at (0.6,-0.9);

% 2nd black last top first row 

\bbox at (1.4,0.3);

% 2nd white top first row

\wbox at (3,0.6);
\wbox at (2.4,0.6);
\wbox at (1.8,0.6);
\wbox at (1.2,0.6);
\wbox at (0.6,0.6);
\wbox at (0,0.6);

% 2nd white front first row

\fwboxh at (2,-0.3);
\fwboxh at (2.6,-0.3);

% 2nd white side first row

\rwbox at (3.2,-0.3);

%  white front third row

\fwbox at (-3,-2.4);
\fwbox at (-2.4,-2.4);

% 1st white side below

\rwboxh at (-1.4, -2.1);
\rwboxh at (-1.8, -2.4);
\fwboxh at (-2.4, -2.4);

% 1st white front right

\fwboxh at (2.6,-0.9);
\rwbox at (3.2, -0.9);

% 1st white top

\wbox at (-2,-1.5);
\wbox at (-2.6,-1.5);
\wbox at (-1.6,-1.2);
\wbox at (-2.2,-1.2);

% 2nd black top first row

\bbox at (0.8,0.3);
\bbox at (0.2,0.3);
\bbox at (-0.4,0.3);
\bbox at (-1,0.3);

%3rd black top first row
\bbox at (0.8,0.9);
\bbox at (0.2,0.9);
\bbox at (-0.4,0.9);
\bbox at (-1,0.9);

%3rd black top -1 row

\bbox at (1.2,1.2);
\bbox at (0.6,1.2);
\bbox at (0,1.2);
\bbox at (-0.6,1.2);

%3rd black top -2 row

\bbox at (1.6,1.5);
\bbox at (1,1.5);
\bbox at (0.4,1.5);
\bbox at (-0.2,1.5);

% 3rd black top 2 row
\bbox at (-0.2,0.6);
\bbox at (-0.8,0.6);
\bbox at (-1.4,0.6);

%3rd black front 2 row

\fbboxh at (-0.6,-0.3);
\fbboxh at (-1.2,-0.3);
\fbboxh at (-1.8,-0.3);

%3rd black side

\rbboxh at (0,-0.3);

%3rd white top

\wbox at (3.4,0.9);
\wbox at (2.8,0.9);
\wbox at (2.2,0.9);

%3rd white side to bottom

\rwbox at (3.6,-0);
\rwbox at (3.6,-0.6);

%3rd black front first row

\fbboxh at (0.4,0);

%3rd black first side first column

\rbboxh at (1,0);
\rbboxh at (1.4, 0.3);
\rbboxh at (1.8, 0.6);

% 2nd black front first row

\fbboxh at (1, -0.6);

% 2nd black side first row

\rbboxh at (1.6,-0.6);

\fwbox at (-2.4,-2.4);

\end{scope}

\begin{scope}[shift={(7,-2.5)}]
 \draw[step=0.5cm, gray, thin] (0,0) grid (3.5,3);

\node at (0.5,2) {$\bm \times$};
\node at (1,2.5) {$\bm \times$};
\node at (0.5,2.5) {\textbullet};

\node at (1,0.5) {\textbullet};
\node at (1.5,0.5) {$\bm \times$};
\node at (1.5,1) {\textbullet};
\node at (2,1) {$\bm \times$};
\node at (2,1.5) {\textbullet};
\node at (3,1.5) {$\bm \times$};
\node at (3,2) {\textbullet};

\node at (0.5,1) {$\bm \times$};
\node at (0.5,1.5) {\textbullet};
\node at (1.5,1.5) {$\bm \times$};
\node at (1.5,2) {\textbullet};
\node at (2,2) {$\bm \times$};
\node at (2,2.5) {\textbullet};

\draw[line width = 2] (0.5,2)--(0.5,2.5);
\draw[line width = 2] (0.5,2.5)--(1,2.5);
\draw[line width = 2] (0.5,2)--(0,2);
\draw[line width = 2] (1,2.5)--(1,3);

\draw[line width = 2] (0.5,1)--(0,1);
\draw[line width = 2] (0.5,1)--(0.5,1.5);
\draw[line width = 2] (0.5,1.5)--(1.5,1.5);
\draw[line width = 2] (1.5,1.5)--(1.5,2);
\draw[line width = 2] (1.5,2)--(2,2);
\draw[line width = 2] (2,2)--(2,2.5);
\draw[line width = 2] (2,2.5)--(3.5,2.5);

\draw[line width = 2] (1,0.5)--(1,0);
\draw[line width = 2] (1,0.5)--(1.5,0.5);
\draw[line width = 2] (1.5,0.5)--(1.5,1);
\draw[line width = 2] (1.5,1)--(2,1);
\draw[line width = 2] (2,1)--(2,1.5);
\draw[line width = 2] (2,1.5)--(3,1.5);
\draw[line width = 2] (3,1.5)--(3,2);
\draw[line width = 2] (3,2)--(3.5,2);

\end{scope}

\end{tikzpicture}
\caption{A generalized plane partition corresponding to three chains.}
\label{plane partition stack}
\end{figure}

If we have a union of $s$ non-intersecting chains, then we think of it as a sequence of $s$ generalized partitions and picture this sequence by stacking the $s$ corresponding Young diagrams from bottom to top. For that, we shift the Young  diagram corresponding to $z$-th chain by the vector $(z-1,z-1,z-1)$. Then,  the result is a picture of height $s$. The non-intersecting condition is equivalent to the \textit{non-floating condition}: every cube of height more than 0 should be sitting on the top of another cube below, see Figure \ref{plane partition stack}.

A generalized plane partition is a two-sided stabilizing ordered  sequence of generalized partitions. More precisely, a generalized plane partition $\hat{\bs\ga}$ is a sequence 
$(\bs \ga_z)_{z\in\Z}$ such that each $\bs \ga_z$ is a generalized partition, $(\bs{\ga}_i)_j\geq (\bs\ga_{i+1})_j$ for all $i,j\in\Z$, and there exist generalized partitions $\bs \ga^{\pm}$ with $\bs \ga_N=\bs \ga^-$ and  $\bs \ga_{-N}=\bs \ga^+$ for large enough $N$.

We picture a generalized plane partition $\hat{\bs\ga}$ by placing the generalized Young diagram of $\ga_z\in \hat{\bs \gamma}$ at height $z\in \mathbb Z$, see Figure \ref{plane partition stack}. In these pictures, whenever $\gamma_z  = (+\infty)_{j\in \mathbb Z}$ or $\gamma_z = (-\infty)_{j\in \mathbb Z}$, we omit the corresponding generalized Young diagrams from the picture for better visualization.

Any monomial in a chain module corresponds to a generalized plane partition with $\bs \ga^-=(-\infty)_{j\in\Z}$  and such that $\bs \ga_0\neq \bs \ga_{-1}=\bs \ga^+=(+\infty)_{j\in\Z}$. We call such generalized plane partitions {\it finite height generalized plane partitions}. For a  finite height generalized plane partition $\hat{\bs \gamma}$, we denote the corresponding monomial by $\bs m_{\hat{\bs \gamma}}$.

For a cube $\mbox{\mancube} = (x,y,z)$ and a generalized plane partition $\hat{\bm \gamma}$, we write  
$(x,y,z)\in \hat{\bs \gamma}$ if the cube $(x,y,z)$ belongs to  the picture of  $\hat{\bs \gamma}$, that is, if $(x,y)\in \ga_z$.

A cube $\mbox{\mancube} = (x,y,z)\in \hat{\bs \gamma}$ is called \textit{convex}  if the removal of this cube produces a generalized plane partition, which we denote by $\hat{\bs\gamma}\setminus \mbox{\mancube}$.  We denote by $\mbox{CV}_i(\hat{\bs \gamma})$ the set of convex cubes of color $i\in \hat I$. Similarly, a cube $\mbox{\mancube} = (x,y,z)\not  \in \hat{\bs \gamma}$ is called \textit{concave}  if the addition of this cube  produces a generalized plane partition, which we denote by $\hat{\bs\gamma}\sqcup \mbox{\mancube}$. We denote by $\mbox{CC}_i(\hat{\bs \gamma})$ the set of concave cubes of color $i\in \hat I$.

Let $\hat{\bm \gamma}$ be a finite height generalized plane partition.  Let $S(\hat{\bm\gamma})$ denote the set of all generalized plane partitions that are obtained from $\hat{\bm\gamma}$ by a finite sequence of adding or removing cubes. It follows that $S(\hat{\bs \gamma})$ is in bijection with the set of monomials $S(\bm m_{\hat{\bs\gamma}})$ obtained from $\bm m_{\hat{\bs\gamma}}$ by positive and negative expansions, cf. Definition \ref{good monomial}.
In particular, the corresponding chain module has a basis labeled by elements of $S(\hat{\bm\gamma})$.

\subsection{The structure constants} \label{constants section}

Recall that the matrix coefficients of the action of $E_i(z), F_i(z)$ on a thin $\mathcal E_n$-module are given by delta functions multiplied by appropriate constants, cf. \eqref{structure constants}. Following \cite{FJM}, we call the set of these constants the supplement.
In this section, we describe the supplement in a combinatorial way following the ideas in Section 3.3 of \cite{JM}. In particular,  we describe the structure  constants $c_{j,\bm r},$ $c'_{j,\bm r}$, $\check c_{j,\bm r}, \check c^{\,'}_{j,\bm r}$ of the Fock and dual  Fock spaces in Section \ref{vh hv chains}.

%Consider the generalized plane partition $\hat{\bs \gamma}$ associated to a monomial $\bs m = \bm m_{\hat{\bs \gamma}}$ formed by $s$ non-intersecting chains. We define the \textit{support} of $\hat{\bs \gamma}$, $\mbox{supp}(\hat{\bm \gamma})$, as the union of all \textcolor{blue}{grey cubes coming from the nontrivial layers at height $z = 0,1,\dots,s-1.$ See Figure \ref{plane partition stack}. }

Below, we consider infinite products of $\ell$-roots. These infinite products are understood as ``layered". Namely, we first multiply and simplify factors corresponding to each fixed $z$-coordinate, and then  multiply the results. After that, the final result is a well-defined monomial in $\hat{\mathcal P}$.

We recover the thin monomial $\bm m_{\hat{\bs \gamma}}$ associated to a finite height generalized plane partition $\hat{\bs \gamma}$ by taking the product over all $\ell$-roots centered at cubes in $\hat{\bs\gamma}$ as follows.

\begin{lem} \label{product over gamma} Let $\hat{\bs \gamma}$ be the finite height generalized plane partition associated with a monomial  $\bs m_{\hat{\bs \gamma}}$ formed by non-intersecting chains. Then,
$$\bs m_{\hat{\bs\gamma}} = \prod_{(x,y,z) \in \hat{\bs \gamma}} \bs A_{x-y,\, c((x,y,z))}^{-1}.$$
\end{lem}

\begin{proof}
It suffices to prove for a single chain monomial. In this case, the cubes in $\hat{\bs \gamma}$ are identified with the boxes to the left and upwards from the chain, see Figures \ref{vh chain plane partition} - \ref{hv chain plane partition}.

If the corner $(u,v)\in (\mathbb Z + \frac{1}{2})^2$ happens to have the neighboring boxes $(u\pm \frac{1}{2}, v\pm \frac{1}{2})$ in $\hat{\bs\gamma}$, then the product of $\bs A_{i,c}^{-1}$ over these boxes is so that the corner $(u,v)$ has no bullet or cross after cancellation. The same conclusion follows if only the neighboring boxes $(u-\frac{1}{2}, v\pm \frac{1}{2})$ or $(u\pm \frac{1}{2}, v-\frac{1}{2})$ are in $\hat{\bs\gamma}$.

Now, if the corner $(u,v)$ is so that only  the neighboring boxes $(u-\frac{1}{2}, v\pm \frac{1}{2}), (u+\frac{1}{2},v-\frac{1}{2})$ are in  $\hat{\bs\gamma}$, then $(u,v)$ is a bullet. Similarly,  if only the neighboring box $(u-\frac{1}{2}, v-\frac{1}{2})$ is in  $\hat{\bs\gamma}$, then $(u,v)$ is a cross. As a result, the product of  $\bs A_{i,c}^{-1}$ over all the boxes in  $\hat{\bs\gamma}$ equals the monomial formed by the chain.
\end{proof}
Note, in particular, that 
$$\prod_{x\in\Z}\bs A^{-1}_{x,a q_1^{-x}}= \prod_{y\in\Z}\bs A^{-1}_{y,a q_3^y}=1.$$

We order the unit cubes in $\mathbb Z^3$ lexicographically by setting $\mbox{\mancube}_1 = (x,y,z)>\mbox{\mancube}_2 = (x',y',z')$ if and only if $z>z'$ or if $z=z', \ y>y'$ or if $z=z', y = y', \ x>x'.$ Whenever $\mbox{\mancube}_1>\mbox{\mancube}_2,$ we say that $\mbox{\mancube}_2$ \textit{precedes} $\mbox{\mancube}_1$ and that $\mbox{\mancube}_1$ \textit{succeeds} $\mbox{\mancube}_2$.

Given a cube $\mbox{\mancube},$ we split $\bm m_{\hat{\bs \gamma}}$ into  factors by multiplying cubes preceding and succeeding $\mbox{\mancube}$ as follows. Let
\begin{equation} \label{smaller and bigger cube} \bm m_{\hat{\bs \gamma}, <\, \mbox{\mancube}} =  \prod_{\substack{ (x,y,z) \in \hat{\bs \gamma}\\ (x,y,z)<\,\mbox{\mancube} }} \bs A_{x-y,\, c((x,y,z))}^{-1}, \qquad \bm m_{\hat{\bs \gamma}, >\, \mbox{\mancube}} =  \prod_{\substack{ (x,y,z) \in \hat{\bs \gamma}\\ (x,y,z)>\,\mbox{\mancube} }} \bs A_{x-y,\, c((x,y,z))}^{-1}.\end{equation}
By Lemma \ref{product over gamma}, if $\mbox{\mancube}\in \mathbb Z^3$ is of color $i\in\hat I,$ we have
\begin{align}   \label{monomial break}
\bm m_{\hat{\bs\gamma}}  = \begin{cases}  \bm m_{\hat{\bs\gamma}, <\, \mbox{\mancube}} \, \bm m_{\hat{\bs\gamma}, >\, \mbox{\mancube}}\, \bs A_{i,c(\mbox{\mancube})}^{-1}, & \quad \mbox{if} \quad \mbox{\mancube}\in \hat{\bs \gamma},\\ \bm m_{\hat{\bs\gamma}, <\, \mbox{\mancube}} \, \bm m_{\hat{\bs\gamma}, >\, \mbox{\mancube}}, & \quad \mbox{if} \quad  \mbox{\mancube} \not\in  \hat{\bs \gamma}.
\end{cases}
\end{align}

For monomials $\bm m_{\hat{\bs\gamma}}, \bm m_{\hat{\bs\gamma}, <\, \mbox{\mancube}}$ and  $\bs m_{\hat{\bs\gamma}, >\, \mbox{\mancube}},$ we denote by $ m_{i,\hat{\bs\gamma}}(z),  m_{i,\hat{\bs\gamma}, <\, \mbox{\mancube}}(z)$ and $m_{i,\hat{\bs\gamma}, >\, \mbox{\mancube}}(z)$ their $i$-th rational functions, $i\in \hat I$. We write  $m_{i,\hat{\bs\gamma}}^{\pm}(z),  m_{i,\hat{\bs\gamma}, <\, \mbox{\mancube}}^{\pm}(z), etc$ to denote the expansion at $z=\infty$ and $z=0,$ respectively.

Let $V(\hat{\bs\gamma}) = \bigoplus_{\hat{\bs \nu}\in S(\hat{\bs \gamma})} \mathbb C \ket{\hat{\bs{\nu}}}$ be a vector space with  basis labeled by generalized plane partitions obtained from $\hat{\bs\gamma}$ by adding and removing cubes.

Let $K_i^\pm(z)$, $F_i(z)$, $E_i(z)$,  $i\in \hat I,$ act on $V(\hat{\bs\gamma})$ by the following formulas
\begin{align}  \label{formulas action}
\begin{split}
& K_i^{\pm}(z) \ket{\hat{\bs \nu}} =  m_{i,\hat{\bs\nu}}^{\pm}(z) \ket{\hat{\bs \nu}},  \\[10pt] 
& F_i(z) \ket{\hat{\bs\nu}} = \sum_{\mbox{\mancube} \, \in CC_i(\hat{\bs \nu})} c_{\hat{\bs\nu}, \, \mbox{\mancube}} \, \delta(z/c(\mbox{\mancube})) \, \ket{{\hat{\bs \nu}\sqcup \,\mbox{\mancube}}}, \\
& E_i(z) \ket{\hat{\bs \nu}} = \sum_{\mbox{\mancube} \, \in CV_i(\hat{\bs \nu})} b_{\hat{\bs\nu}, \, \mbox{\mancube}} \, \delta(z/c(\mbox{\mancube})) \, \ket{{\hat{\bs \nu}\setminus \mbox{\mancube}}},
\end{split}
\end{align}
where the supplement is given by
\begin{align}
& c_{\hat{\bs \nu}, \, \mbox{\mancube}} =  m_{i, \hat{\bs\nu}, >\, \mbox{\mancube}}(c(\mbox{\mancube})),
\\
& b_{{\hat{\bs \nu}}, \, \mbox{\mancube}} = \frac{1}{q-q^{-1}}\mathop{\mbox{res}}_{z = c(\mbox{\mancube})} m_{i, \hat{\bs\nu}, <\, \mbox{\mancube}}(z)  \frac{dz}{z}.
\end{align}

One can check explicitly that formulas \eqref{formulas action} are consistent. Namely, we have the following theorem.

\begin{thm} \label{plane partition module} Let $\hat{\bs \gamma}$ be a finite height generalized plane partition. Then formulas \eqref{formulas action} define a structure of an irreducible, locally finite, thin $\mathcal E_n$-module on $V(\hat{\bs\ga})$.

The module $V(\hat{\bs\gamma})$ is isomorphic to the chain module associated with the ray configuration of  $\bm m_{\hat{\bs\gamma}}.$
 \end{thm}
\begin{proof}
The proof is similar to the proof of Theorem 3.5 in \cite{JM}.
\end{proof}

We note that Theorem \ref{plane partition module} provides another proof of the existence of the chain modules $V_{\lambda, \mu, \nu}^{\tilde \lambda, \tilde \mu, \tilde \nu}$ and $H_{\lambda, \mu, \nu}^{\tilde \lambda, \tilde \mu, \tilde \nu}$, cf. Theorem  \ref{H and V modules}.

\section{Representations of shifted quantum toroidal algebras}\label{sec: shifted}

\subsection{Shifted quantum toroidal algebras}

\label{shifted algebras}

We now recall the  shifted quantum toroidal algebras, see \cite{N}, cf. \cite{FT}.

A shift is a  $\mathbb Z$-valued $\hat I$-tuple $\textbf{s} = (s_0,\dots, s_{n-1})$. For a shift $\bm s$, the \textit{\textbf{s}-shifted quantum toroidal $\mathfrak{gl}_n$ algebra } $\mathcal E_{n,\textbf{s}} =\mathcal E_{n,\textbf{s}}(q_1,q_2,q_3) $ is the unital associative algebra with generators  $E_{i,k}, F_{i,k}, H_{i,l}$ and invertible generators $K_i^{\pm}$, $i\in \hat I, k\in \mathbb Z, l \in \mathbb Z\setminus\{0\},$ satisfying the same relations \eqref{relations 2.2}-\eqref{relations 2.9} as before, except that the currents $K_i^{\pm}(z)$  have the form
\begin{align*} & K_i^{-}(z) = K_i^-\, \exp \big(-(q-q^{-1})\sum_{l=1}^{\infty}H_{i,- l}z^l\big), \\
&K_i^+(z) = K_i^+ z^{s_i} \exp \big((q-q^{-1})\sum_{l=1}^{\infty}H_{i, l}z^{-l}\big).
\end{align*}

Notice that the quantum toroidal algebra $\mathcal E_n$ is a quotient of the $\bm 0$-shifted quantum toroidal algebra $\mathcal E_{n,\bs 0}$ by the relations $K_i^{-}K_i^+ = 1,$  $i\in\hat I$.

Let $\bm s_1,\bm s_2$ be shifts. The coproduct formula \eqref{coproduct} defines a coproduct map $$\Delta_{\bm s_1, \bm s_2}: \mathcal E_{n,\bm s_1+\bm s_2} \to \mathcal E_{n, \bm s_1} \widehat{\otimes}\,  \mathcal E_{n,\bm s_2}.$$
Here, the hat indicates that the map $\Delta_{\bm s_1, \bm s_2}$ applied to an element of $\mathcal E_{n,\bm s_1+\bm s_2}$  may produce an infinite sum. As with the coproduct  \eqref{coproduct}, we use   $\Delta_{\bm s_1, \bm s_2}$ only if all these sums are convergent.

The shifted algebra  $\mathcal E_{n,\bs 0}$ has the following family of one-dimensional modules.
For $\bs a=(a_0,\dots,a_{n-1})$, $a_i\in\C^\times$, let $L_{\bs a}=\C v_{\bs a}$ denote the $1$-dimensional $\mathcal E_{n,\bm 0}$-module such that  $$E_i(z)\cdot v_{\bs a} = F_i(z)\cdot v_{\bs a} =0, \quad K_i^{\pm }(z)\cdot v_{\bs a}   = a_iv_{\bs a}.$$

Let $V$ be an irreducible thin  $\mathcal E_{n,\bs 0}$-module. Let $\bs \phi= (\phi_i(z))$ be an $\hat I$-tuple of rational functions such that $\phi(0)$ and $\phi(\infty)$ are well-defined and non-zero. Suppose $(\phi_i^\pm(z))$ is an $\ell$-weight in $V$. Then, 
$$
V\otimes L_{\bs a}, \qquad a_i^{-2}=\phi_i(0)\phi_i(\infty),
$$
has balanced $\ell$-weights and therefore factors to an $\mathcal E_n$-module. 

More generally, let $\bs a(z)=(a_0(z),\dots,a_{n-1}(z))$, $a_i(z)\in\C[z]$, $a_i(0)\neq 0$, be an $\hat I$-tuple of non-zero polynomials. Let $\bs s=(\deg a_i(z))_{i\in\hat I}$. Then, we have a one-dimensional $\E_{n,\bs s}$-module $L_{\bs a(z)}=\C v_{\bs a(z)}$ such that 
$$E_i(z)\cdot v_{\bs a(z)} = F_i(z)\cdot v_{\bs a(z)} =0, \quad K_i^{\pm }(z)\cdot v_{\bs a(z)}   = a_i(z)v_{\bs a(z)}.$$

\subsection{Infinite chain configurations and truncations} \label{infinite chain configurations}
We study thin $\E_{n,\bs s}$-modules by the same methods used for the locally finite thin $\E_n$-modules. This leads to the concept of truncated chain configurations, which we now develop. These configurations allow unions of infinitely many non-intersecting chains and include truncation boxes that restrict their expansions. In particular, in this section we call the finite unions of non-intersecting chains used in the study of locally finite modules the finite chain configurations.

Consider a union of infinitely many non-intersecting chains. As before, the chains are ordered from bottom-right  to top-left. Denote the chains by $C_s$ and let $C=\cup_s C_s$. We have three cases: $C_s$ is a chain for all $s\in\Z$, $s\in\Z_{>0}$, or $s\in\Z_{<0}$. In the first case, we have no smallest and no largest chain. In the second case, we have a smallest chain. In the third case, we have a largest chain. We set $C_s=\varnothing$ for $s\leq 0$ in the second case, and $C_s=\varnothing$ for $s\geq 0$ in the third case. Then,  we simply write $C=\cup_{s\in\Z} C_s$ in all three cases.

In all cases, we define the expansion of $C$ as before. Namely, we allow the expansion of any chain $C_s$ at a bullet or cross as long as the result is a union of non-intersecting chains. 

First, we rule out the cases where $C$ can be expanded in infinitely many places, as we need the number of possible expansions of $C$ to be finite in order to define an action of the toroidal algebra, see Section \ref{infinite chain modules}. 

We say that a chain $C$ is a shift of a chain $C'$ if $C$ is obtained from $C'$ by a shift by $q_2^{-1}$. In other words, $(u,v)$ is a bullet (or cross) in $C$ if and only if $(u+1,v+1)$ is a bullet (or cross) in $C'$. 

Note that a shift of an empty chain is an empty chain.

If $C$ is a shift of $C'$, then $C$ and $C'$ are non-intersecting. Moreover, if both $C$ and $C'$ are present in a union of non-intersecting chains, then the crosses of $C'$ and the bullets of $C$ cannot be expanded.

\begin{lem}\label{finitely many poles}  Let $C=\cup_{i\in\Z} C_i$ be a union of infinitely many non-intersecting chains. Then, $C$ can be expanded only in  finitely many 
places if and only if $C_{s+1}$ is a shift of $C_{s}$ for all $s\gg 0$ and $s\ll 0$.
\end{lem}
 \begin{proof}
If $C_{i}$ is a shift of $C_{i-1}$ and $C_{i+1}$ is a shift of $C_{i}$, then neither bullets nor crosses of $C_i$ can be expanded. The if part of the lemma follows.
 
 If $C_{s+1}$ is not  a shift of $C_{s}$ and both are non-empty and of the same type, then there is either a bullet of $C_{s+1}$ or a cross of $C_{s}$ which can be expanded. That proves the only if part.
 \end{proof}

\begin{figure}[H]
\begin{tikzpicture}

\begin{scope}[shift={(0,-1)}]
 
 \draw[step=0.5cm, gray, very thin] (-2,0.5) grid (1.5,3.5);

\node at (-1.75, 3.4) {$\ddots$};

\node at (-1.5,2.5) {$\bm \times$};
\node at (-1,2) {$\bm \times$};
\node at (-0.5,1.5) {$\bm \times$};
\node at (0,1) {$\bm \times$};
\node at (-1.5,3) {\textbullet};
\node at (-1,2.5) {\textbullet};
\node at (-0.5,2) {\textbullet};
\node at (0,1.5) {\textbullet};
\node at (-1,3) {$\bm \times$};
\node at (-0.5,2.5) {$\bm \times$};
\node at (0,2) {$\bm \times$};
\node at (0.5,1.5) {$\bm \times$};

\draw[line width = 2] (-2,2.5)--(-1.5,2.5);
\draw[line width = 2] (-1.5,2.5)--(-1.5, 3);
\draw[line width = 2] (-1.5,3)--(-1,3);
\draw[line width = 2] (-1,3)--(-1,3.5);
\draw[line width = 2] (-1,2)--(-2,2);
\draw[line width = 2] (-1,2)--(-1,2.5);
\draw[line width = 2] (-1,2.5)--(-0.5,2.5);
\draw[line width = 2] (-0.5,2.5)--(-0.5,3.5);
\draw[line width = 2] (-0.5,1.5)--(-2,1.5);
\draw[line width = 2] (-0.5,1.5)--(-0.5,2);
\draw[line width = 2] (-0.5,2)--(0,2);
\draw[line width = 2] (0,2)--(0,3.5);
\draw[line width = 2] (0,1)--(-2,1);
\draw[line width = 2] (0,1)--(0,1.5);
\draw[line width = 2] (0,1.5)--(0.5,1.5);
\draw[line width = 2] (0.5,1.5)--(0.5,3.5);

\draw[line width = 2] (0.5,1.5)--(0.5,3.5);

\end{scope}

\begin{scope}[shift={(5,-1)}]

\draw[step=0.5cm, gray, very thin] (-2,0.5) grid (1.5,3.5);

\node at (-1,2.5) {\textbullet};
\node at (-0.5,2.5) {$\bm \times$};
\node at (-0.5,3) {\textbullet};
\node at (0,3) {$\bm \times$};

\draw[line width = 2] (-1,2.5)--(-1,0.5);
\draw[line width = 2] (-1,2.5)--(-0.5,2.5);
\draw[line width = 2] (-0.5,2.5)--(-0.5,3);
\draw[line width = 2] (-0.5,3)--(0,3);
\draw[line width = 2] (0,3)--(0,3.5);

\node at (-0.5,2) {\textbullet};
\node at (-0,2) {$\bm \times$};
\node at (-0,2.5) {\textbullet};
\node at (0.5,2.5) {$\bm \times$};

\draw[line width = 2] (-0.5,2)--(-0.5,0.5);
\draw[line width = 2] (-0.5,2)--(-0,2);
\draw[line width = 2] (0,2)--(0,2.5);
\draw[line width = 2] (0,2.5)--(0.5,2.5);
\draw[line width = 2] (0.5,2.5)--(0.5,3.5);

\node at (0,1.5) {\textbullet};
\node at (0.5,1.5) {$\bm \times$};
\node at (0.5,2) {\textbullet};
\node at (1,2) {$\bm \times$};

\draw[line width = 2] (0,1.5)--(0,0.5);
\draw[line width = 2] (0,1.5)--(0.5,1.5);
\draw[line width = 2] (0.5,1.5)--(0.5,2);
\draw[line width = 2] (0.5,2)--(1,2);
\draw[line width = 2] (1,2)--(1,3.5);

\node at (0.5,1) {\textbullet};
\node at (1,1) {$\bm \times$};
\node at (1,1.5) {\textbullet};
\node at (1.5,1.5) {$\bm \times$};

\draw[line width = 2]  (0.5,0.5)--(0.5,1);
\draw[line width = 2]  (0.5,1)--(1,1);
\draw[line width = 2] (1,1)--(1,1.5);
\draw[line width = 2] (1,1.5)--(1.5,1.5);
\draw[line width = 2] (1.5,3.5)--(1.5,1.5);

\node at (1.3, 0.75) {$\ddots$};

\end{scope}

\begin{scope}[shift={(10,-1)}]

\draw[step=0.5cm, gray, very thin] (-2,0.5) grid (1.5,3.5);

\node at (-1.75, 3.4) {$\ddots$};

\node at (-1.5,2.5) {$\bm \times$};
\node at (-1,2) {$\bm \times$};
\node at (-1.5,3) {\textbullet};
\node at (-1,2.5) {\textbullet};
\node at (-1,3) {$\bm \times$};
\node at (-0.5,2.5) {$\bm \times$};

\node at (0,1.5) {$\bm \times$};
\node at (0,2) {\textbullet};
\node at (0.5,1) {$\bm \times$};
\node at (0.5,1.5) {\textbullet};
\node at (1,0.5) {$\bm \times$};
\node at (1,1) {\textbullet};
\draw[line width = 2] (0,1.5)--(-2,1.5);
\draw[line width = 2] (0,1.5)--(0,2);
\draw[line width = 2] (0,2)--(1.5,2);
\draw[line width = 2] (0.5,1)--(-2,1);
\draw[line width = 2] (0.5,1)--(0.5,1.5);
\draw[line width = 2] (0.5,1.5)--(1.5,1.5);

\node at (1.3, 0.75) {$\ddots$};

\draw[line width = 2] (-2,2.5)--(-1.5,2.5);
\draw[line width = 2] (-1.5,2.5)--(-1.5, 3);
\draw[line width = 2] (-1.5,3)--(-1,3);
\draw[line width = 2] (-1,3)--(-1,3.5);
\draw[line width = 2] (-1,2)--(-2,2);
\draw[line width = 2] (-1,2)--(-1,2.5);
\draw[line width = 2] (-1,2.5)--(-0.5,2.5);
\draw[line width = 2] (-0.5,2.5)--(-0.5,3.5);
\draw[line width = 2] (1,0.5)--(-2,0.5);
\draw[line width = 2] (1,0.5)--(1,1);
\draw[line width = 2] (1,1)--(1.5,1);

\end{scope}

\end{tikzpicture}
\caption{Examples of infinite chain configurations.}
\label{all possible infinite hv chains}
\end{figure}
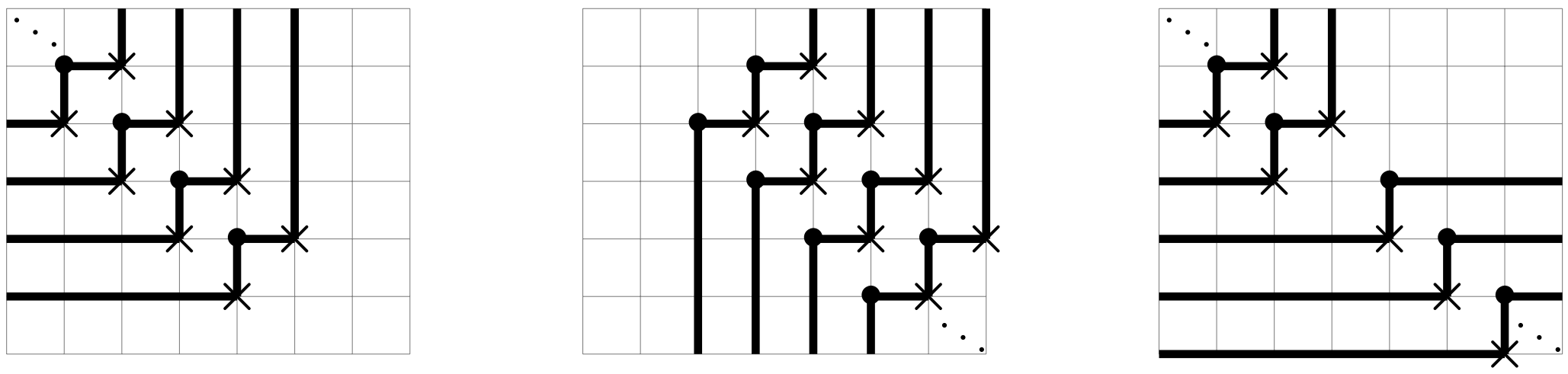

We call a union of infinitely many chains $C$ an infinite chain configuration if $C$ satisfies the conditions of Lemma \ref{finitely many poles}. For an infinite chain configuration $C$, we call the chains $C_s$ with $s\gg 0$ and $s\ll 0$ the positive and negative tails of $C$. Note that in our definition the tail is defined only up to a shift. Examples of infinite chain configurations are given in Figure \ref{all possible infinite hv chains}.

We denote by $S(C)$ the set of all chain configurations obtained by finitely many expansions of $C$.

For an infinite chain configuration $C = \cup_{z\in \mathbb Z} C_{z},$ we associate a generalized plane partition $\hat{\bs\gamma}_C$ as before. Namely, for $z\in \mathbb Z,$  we shift the generalized Young diagram representing the chain $C_{z}$ by the vector $(z-1,z-1,z-1),$ and we stack the Young diagram of
$C_{z+1}$ on the top of $C_z$. The collection of cubes thus obtained has infinite height and grows indefinitely in the positive or negative $z$-direction. We call \textit{infinite height generalized plane partitions} the generalized plane partitions of this form. See Figure \ref{vertical fock 321} for an example.

Given a finite or infinite chain configuration $C$, we introduce truncation boxes which restrict the expansions of $C$. The set of truncation boxes $B$ is a finite collection of boxes $\square_k=(b_k,b_k')\in \mathbb Z^2$. We call the set $\mathcal C =(C,B)$ a truncated chain configuration. The usual finite or infinite chain configurations are the ones in which the set of truncation boxes is empty, $B = \varnothing$.  

We define the expansion of a truncated chain configuration  $\mathcal C = (C,B)$ as before, except that we prohibit the expansion of any box in $B$. That means that we do not allow to expand  bullets at $(b_k-\frac{1}{2},b_k'-\frac{1}{2})$, and we do not allow to expand crosses at $(b_k+\frac{1}{2},b_k'+\frac{1}{2})$ for $(b_k,b_k')\in B$. 

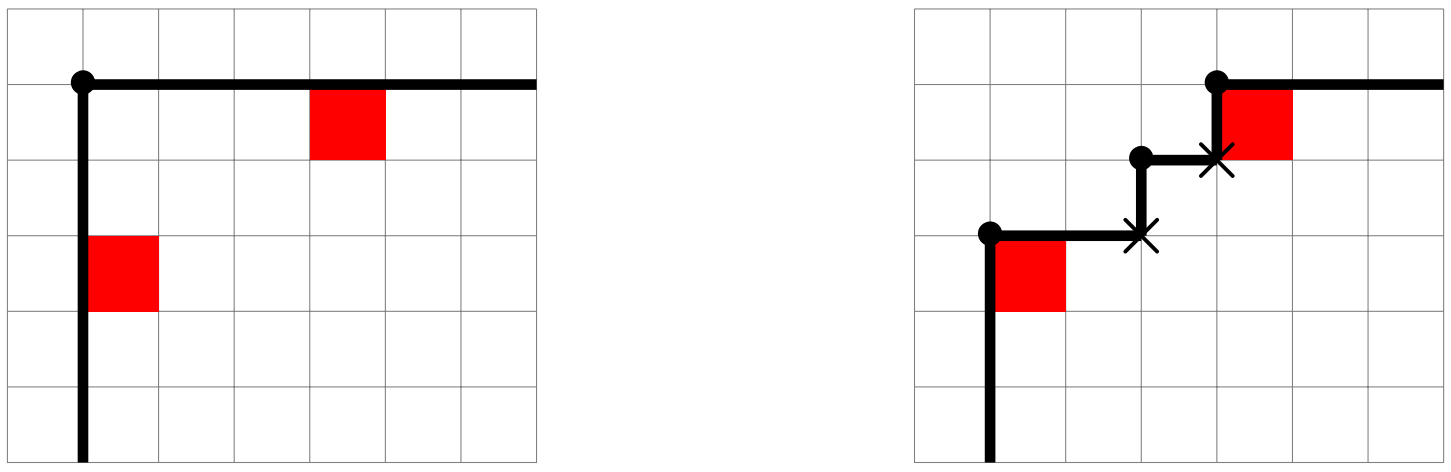
\begin{figure}[H]   
\begin{tikzpicture}

\begin{scope}[shift={(0,0)}]

\draw[step=0.5cm, gray, very thin] (-2,0.5) grid (1.5,3.5);

\node at (-1.5,3) {\textbullet};

\redbox at (0,2.5);

\redbox at (-1.5,1.5);

\draw[line width = 2] (-1.5,3)--(-1.5,0.5);
\draw[line width = 2] (-1.5,3)--(1.5,3);

\end{scope}

\begin{scope}[shift={(6,0)}]

\draw[step=0.5cm, gray, very thin] (-2,0.5) grid (1.5,3.5);

\redbox at (0,2.5);

\redbox at (-1.5,1.5);

\node at (0,3) {\textbullet};
\node at (0,2.5) {$\bm \times$};
\node at (-0.5,2.5) {\textbullet};
\node at  (-0.5, 2) {$\bm \times$};
\node at (-1.5,2) {\textbullet};

\draw[line width = 2] (-1.5,2)--(-1.5,0.5);
\draw[line width = 2] (-1.5,2)--(-0.5,2);
\draw[line width = 2] (-0.5,2)--(-0.5,2.5);
\draw[line width = 2] (-0.5,2.5)--(0,2.5);
\draw[line width = 2] (0,2.5)--(0,3);
\draw[line width = 2] (0,3)--(1.5,3);

\end{scope}

\end{tikzpicture}
\caption{Expansion of a $vh$ chain with two truncation boxes in red.}
\label{vh truncated chain}
\end{figure}

In the pictures we mark the prohibited boxes in red, see Figure \ref{vh truncated chain}. Recall that every chain $C_i$ divides the plane in two regions. Combinatorially, if a truncation box belongs to one  of those two regions, then this truncation box stays in the corresponding region after any expansion. In the example of Figure \ref{vh truncated chain}, the expansions produce a total of ten different states.

A truncation box $(b,b')$ of color $i=b-b'$ without any chains corresponds to a one-dimensional module $L_{a(z)}$, where $a(z)=(1,\dots,1,z-q_1^{-b}q_3^{-b'},1\dots,1)$ with the non-trivial coordinate in the $i$-th position. Adding truncations means multiplying by the corresponding one-dimensional modules and taking the irreducible quotient. Such a procedure was used in Lemma 4.12 of \cite{FFJMM4}, and in Proposition 6.7 of \cite{FJMM3} to construct two-step modules used in the Bethe ansatz. The truncations of the MacMahon module are discussed already in Section 4.6 of \cite{FJMM1}. The terminology is taken from \cite{FJM}, see Section 2.4, where such procedure was used in the context of $qq$-characters.

\subsection{The \texorpdfstring{$\ell$}{}-weights of  truncated chain configurations}       \label{l weight chain config}

To a truncated chain configuration  we associate an $\hat{I}$-tuple of rational functions as follows.

For $i\in\hat I, a\in\C^\times$, let $ X_{i,a}$ be the $\hat I$-tuple of rational functions $$X_{i,a} = (1,\dots,1,z-a,1,\dots,1),$$ where $z-a$ is at the $i$-th position. Clearly,  we have  $\bm i_a =(1,\dots, 1, q,1\dots,1)\cdot    X_{i,aq^{-1}}X_{i,aq}^{-1}$, where $q$ lies in the $i$-th position.

For a truncated chain configuration $\mc C=(C,B)$, let $b_i(C)\subset (\mathbb Z+\frac{1}{2})^2$ and $r_i(C)\subset (\mathbb Z+\frac{1}{2})^2$ be the sets of bullets and crosses of $C$ of color $i$, respectively. Similarly, let $T_i(B)$ be the set of truncation boxes of color $i$. Recall that the content of the half-integer point $w = (x+\frac{1}{2},y+\frac{1}{2})\in (\mathbb Z+\frac{1}{2})^2$ is the scalar $c(w)=qq_1^{-x}q_3^{-y}$ and the content of a box $\square = (b,b')\in \mathbb Z^2$ is the scalar $c(\square) = q_1^{-b}q_3^{-b'}.$

Define the $\hat I$-tuple $\bs m_{\mathcal C}$ of rational functions corresponding to $\mathcal C = (C,B)$ by 
$$
\bs{m}_{\mathcal C}=\prod_{i\in\hat I}\bigg(\prod_{w\in r_i(C)} X_{i,q^{-1}c(w)}^{-1}X_{i,qc(w)} \prod_{w\in b_i(C)}  X_{i,q^{-1}c(w)}X_{i,qc(w)}^{-1}  \prod_{\square \in  T_i(B)} X_{i,c(\square)}\bigg).
$$

While the product above is infinite in the case of an infinite chain configuration, after all cancellations  the result is a finite product. Indeed, let $C=\cup_{s\in\Z}C_s$ be such that $C_{s+1}$ is a shift of $C_s$ for all $s\geq N$  and $s\leq -N.$ Consider $C^{(N)}=\cup_{s\geq N} C_s$ and $C_{(N)}=\cup_{s\leq-N} C_s$. Then,
\begin{align}   \label{tails}
\begin{split}\bs m_{(C^{(N)},\varnothing)}&=\prod_{i\in\hat I}\Big(\prod_{w\in r_i(C_N)}X_{i,qc(w)} \prod_{w\in b_i(C_N)}X_{i,qc(w)}^{-1}\Big),\\
 \bs m_{(C_{(N)},\varnothing)}&=\prod_{i\in\hat I}\Big(\prod_{w\in r_i(C_{-N})}X_{i,q^{-1}c(w)}^{-1} \prod_{w\in b_i(C_{-N})}X_{i,q^{-1}c(w)}\Big).
 \end{split}
\end{align}
are both finite products. 

Clearly, if $C=C'\sqcup C''$ and $B=B'\sqcup B''$
is a union of disjoint chain configurations with no overlapping truncation boxes, then $\bs m_{\mc C}=\bs m_{\mc C'}\cdot \bs m_{\mc C''}$, where  $\mc C'=(C',B')$, $\mc C''=(C'',B'')$.

Note also that if $C$ is a finite truncated chain configuration with $B=\varnothing$, then $\bs m_{\mathcal C}$ differs from the monomial defined in Section \ref{monomial pic} by a factor of the form $(q^{d_0},q^{d_1},\dots, q^{d_{n-1}})$, where $q^{-2d_i}$ is the value of $i$-th rational function in the $\hat I$-tuple $\bs m_{\mathcal C}$ at $z=0$. 

\medskip

The expansions of chains in $\mc C=(C,B)$ do not change $B$. So, now we fix $B$ and we omit $B$ from the notation of monomials. For example, we write simply $\bs m_C$ for $\bs m_{\mc C}$.

Let $(C^+,B)$ be a truncated chain configuration. Let $C$ be obtained from $C^+$ by finitely many expansions of bullets or crosses. For $i\in \hat I$, let $\omega_i(C|C^+)$ be the difference between the number of expansions of bullets of color $i$ and the number of expansions of crosses of color $i$ in $C$. Let $\al_i=(1,\dots, 1, q^{-1},q^2,q^{-1},1,\dots,1)$, where $q^2$ is in the $i$-th position. We set
$$   
\bs m_{C^+\!\to\,C}=\omega_{C^+\!\to \,C}\cdot \bs m_C, \qquad \mbox{where} \quad \omega_{C^+ \! \to\, C}=\prod_{i\in\hat I}\al_i^{-\omega_{i}(C|C^+)}.
$$
We call the $\hat I$-tuple $\omega_{C^+\to C}$ the relative weight at $C$ with respect to $C^+$, and $\bm m_{C^+ \! \to \, C}$ the monomial of $C$ relative to $C^+$.

Then, if $C'$ is obtained from $C$ by an expansion of  a bullet of color $i$ and content $c$, then $\bs m_{C^+\!\to\, C'}=\bs m_{C^+\!\to \, C} \cdot \bm A^{-1}_{i,qc}$. Similarly, if $C'$ is obtained from $C$ by an expansion of a cross of color $i$ and content $c$, then $\bs m_{C^+\! \to \, C'}=\bs m_{C^+\!\to \, C} \cdot \bm A_{i,q^{-1}c}$. In particular, $\bs m_{C^+ \! \to \, C}$ differs from $\bs m_{C^+}$ by a finite product in the variables $\bm A^{\pm 1}_{i,c}, i\in \hat I, c\in \mathbb C^\times.$

We write $\bs m_{C^+\! \to\, C}=\bs m_{C^+}\prod_{s=1}^k\bs A_{i_s,c_s}\prod_{s=1}^l\bs A_{j_s,b_s}^{-1}$. In this product, the factors correspond to all boxes which have to be added or removed to obtain $C$ from $C^+$. In particular, some $(i_s,c_s)$ and some $(j_s,b_s)$ could be the same (but still correspond to different boxes in the sense of generalized plane partitions). We write as a product over boxes
$$
\bs m_{C^+\!\to\, C}=\bs m_{C^+}\prod_{ (x,y,z)  \in \hat{\bs \gamma}_{C^+}\setminus \hat{\bs\gamma}_C}\bs A_{x-y,c((x,y,z))}\prod_{{   (x,y,z) \in \hat{\bs \gamma}_{C}\setminus \hat{\bs\gamma}_C^+}}\bs A_{x-y,c((x,y,z))}^{-1}.
$$

Then, for any  cube $\mbox{\mancube} = (x,y,z)$, we split $\bs m_{C^+\!\to \, C}$ into two parts $\bs m_{C^+\!\to \, C, >\,\mbox{\mancube}}$ and $\bs m_{C^+\!\to \, C, <\,\mbox{\mancube}},$ according to the order of boxes:
\begin{align*}
&\bs m_{C^+\!\to\, C,<\,\mbox{\mancube}} = \bs m_{C^+}\prod_{ \substack{(x,y,z)  \in \hat{\bs \gamma}_{C^+}\setminus \hat{\bs\gamma}_C \\ (x,y,z)<\,\mbox{\mancube}}}\bs A_{x-y,c((x,y,z))}\prod_{\substack{ (x,y,z) \in \hat{\bs \gamma}_{C}\setminus \hat{\bs\gamma}_C^+\\ (x,y,z)<\,\mbox{\mancube}}}\bs A_{x-y,c((x,y,z))}^{-1}.
\\
&\bs m_{C^+\!\to \, C,>\,\mbox{\mancube}}=\prod_{ \substack{(x,y,z)  \in \hat{\bs \gamma}_{C^+}\setminus \hat{\bs\gamma}_C \\ (x,y,z)>\,\mbox{\mancube}}}\bs A_{x-y,c((x,y,z))}\prod_{\substack{ (x,y,z) \in \hat{\bs \gamma}_{C}\setminus \hat{\bs\gamma}_C^+\\ (x,y,z)>\,\mbox{\mancube}}}\bs A_{x-y,c((x,y,z))}^{-1}.
\end{align*}

 This split is similar to \eqref{smaller and bigger cube}, but differs by a factor independent of $C$.

As before, for $i\in \hat I,$ we denote by $ m_{i,C^+\!\to\, C}(z),   m_{i, C^+\!\to\, C,<\,\mbox{\mancube}}(z),$ etc, the $i$-th rational function of $\bs m_{C^+\!\to\, C}, \bs m_{ C^+\!\to\, C,<\,\mbox{\mancube}}$, etc. We write  $m_{i,C^+\!\to\, C}^{\pm}(z),  m^{\pm}_{i,C^+\!\to\, C, <\, \mbox{\mancube}}(z)$, etc, to denote the expansion at $z=\infty$ and $z=0,$ respectively.

\subsection{Thin representations of shifted quantum toroidal algebras}  \label{infinite chain modules}

Let $\mathcal C^+ = (C^+,B)$ be a truncated chain configuration.
Let $V(\mathcal C^+) = \bigoplus_{C\in S(C^+)} \mathbb C \ket{C}$ be a vector space with a basis labeled by the monomials of truncated chain configurations $C$ obtained from $C^+$ by finitely many expansions. 

Recall that $\bs m_{C^+\!\to C}$ is the $\hat I$-tuple defined in Section \ref{l weight chain config}.
For $C\in S(C^+)$ and $\mbox{\mancube}\in CC_i(\hat \gamma_C)$,  we denote $C\cup \mbox{\mancube}$ the truncated chain configuration obtained from $C$ by expanding at the box $\mbox{\mancube}$. In particular, we have $\bs m_{C^+\!\to \, C\cup \, \mbox{\mancube}} =\bs m_{C^+\!\to \, C} \bs A^{-1}_{i,c(\mbox{\mancube})}$. 

Similarly, for $C\in S(C^+)$ and $\mbox{\mancube}\in CV_i(\hat \gamma_C)$,  we denote $C\setminus \mbox{\mancube}$ the truncated chain configuration obtained from $C$ by expanding at the box $\mbox{\mancube}$. Let $\bm B$ be the set of all cubes which have contents equal to the contents of the truncation boxes in $B$.
 
Let $K_i^\pm(z)$, $F_i(z)$, $E_i(z)$,  $i\in \hat I,$ act in $V(\mc C^+)$ by the following formulas
\vspace{0.1cm}
\begin{align}   \label{shift formulas}
\begin{split}
& K_i^{\pm}(z) \ket{ C} =  m_{i, C^+\!\to \, C}^{\pm}(z) \ket{ C},  \\[10pt] 
& F_i(z) \ket{ C} = \sum_{\mbox{\mancube} \, \in CC_i(\hat{\bs \gamma}_C)\setminus{\bm B}} c_{{C^+ \!\to\, C}, \, \mbox{\mancube}} \, \delta(z/c(\mbox{\mancube})) \, \ket{ C \sqcup \mbox{\mancube}}, \\
& E_i(z) \ket{C} = \sum_{\mbox{\mancube} \, \in CV_i(\hat{\bs \gamma}_C)\setminus{\bm B}} b_{C^+ \!\to \, C, \, \mbox{\mancube}} \, \delta(z/c(\mbox{\mancube})) \, \ket{C \setminus \mbox{\mancube} },
\end{split}
\end{align}
where the supplement is given by
\begin{align}
& c_{{C^+\!\to\,C}, \, \mbox{\mancube}} =  m_{i, C^+\!\to \, C, >\, \mbox{\mancube}}(c(\mbox{\mancube})),
\\
& b_{C^+\!\to\, C, \, \mbox{\mancube}} = \frac{1}{q-q^{-1}}\mathop{\mbox{res}}_{z = c(\mbox{\mancube})} m_{i, C^+\!\to C\,, <\, \mbox{\mancube}}(z)  \frac{dz}{z}.
\end{align}

Let $\bm s(\mc C^+) = (s_0,\dots, s_{n-1})$ be the shift where $s_i$ denotes the degree of the $i$-th rational function of the $\hat I$-tuple $\bm m_{C^+}=\bs m_{\mc C^+}.$ Note that the shift  depends on the truncation boxes $B$.

\begin{thm}  \label{truncated chain theorem} Let $\mathcal C^+ = (C^+,B)$ be a truncated chain configuration. Then, formulas \eqref{shift formulas} define a structure of an irreducible, thin $\mathcal E_{n,\bs s(\mc C^+)}$-module on $V(\mathcal C^+)$.
\end{thm} 
\begin{proof}
The proof is similar to the proof of Theorem 3.5 in \cite{JM}, cf. Theorem \ref{plane partition module}.
\end{proof}
Note that if $C$ is obtained from $C^+$ by finitely many expansions, then the modules  $V(\mathcal C^+)$  and $V(\mathcal C)\otimes L_{\bs a}$ are isomorphic, where $\mathcal C = (C,B), \mathcal C^+ = (C^+, B)$ and   $L_{\bs a}$ is the one-dimensional module with $\bs a=\omega_{C^+\to C}$.

We call a thin $\mathcal E_{n,\bs s}$-module $V$ 
\textit{totally prime} if $V$ is not isomorphic to a tensor product of two modules $W_1$ and $W_2$ of $\mathcal E_{n,\bs s_1}$ and $\mathcal E_{n,\bs s_2}$, where neither $W_1$ nor $W_2$ have the form $L_{\bs a}$, $\bs a\in(\C^\times)^n$. 

For example, all modules corresponding to the configurations in Figures \ref{all possible infinite hv chains}-\ref{vh truncated chain} are totally prime. However, the  MacMahon $\E_{n,\bs 0}$-module, see \cite{FJMM1},  is thin and prime but not totally prime, cf. Figure \ref{MacMahon} below.

We propose the following conjecture on irreducible, thin, totally prime modules over the shifted quantum toroidal algebras.

\begin{conj}\label{main conjecture} Let $V$ be an irreducible, thin, and totally prime $\mc E_{n,\bs s}$-module. Then, $V$ is isomorphic to a module obtained from a truncated chain configuration as in Theorem \ref{truncated chain theorem} by a shift $\sigma_a$, $a\in\C^\times$.
\end{conj}

\subsection{Examples}
We give a few explicit examples.

\begin{example}

Let $C$ be an infinite chain configuration such that the negative tail is empty and positive tail is a chain of $vh$-type. We call such a $C$ the MacMahon configuration.  

Clearly, all chains of a MacMahon configuration are of $vh$-type.
The corresponding module has a basis which can be labeled by plane partitions with three boundary conditions, see \cite{FJMM1}. 

The MacMahon module with empty boundary conditions corresponds to the infinite chain configuration given in Figure \ref{MacMahon}.

\begin{figure}[H]   
\begin{tikzpicture}

\begin{scope}[shift={(0,-1)}]
 
\draw[step=0.5cm, gray, very thin] (-2,0.5) grid (1.5,3.5);

\node at (-1.5,3) {\textbullet};
\node at (-1, 2.5) {\textbullet};
\node at (-0.5,2) {\textbullet};
\node at (0, 1.5) {\textbullet};

\draw[line width = 2] (-1.5,3)--(-1.5,0.5);
\draw[line width = 2] (-1.5,3)--(1.5,3);
\draw[line width = 2] (-1,2.5)--(-1,0.5);
\draw[line width = 2] (-1,2.5)--(1.5,2.5);
\draw[line width = 2] (-0.5,2)--(-0.5,0.5);
\draw[line width = 2] (-0.5,2)--(1.5,2);
\draw[line width = 2] (0,1.5)--(0,0.5);
\draw[line width = 2] (0,1.5)--(1.5,1.5);

\draw (0.03,1)--(0.505,1.5);

\node at (-1.8, 3.4) {$\ddots$};

\end{scope}

\end{tikzpicture}
\caption{MacMahon infinite chain configuration.}
\label{MacMahon}
\end{figure}
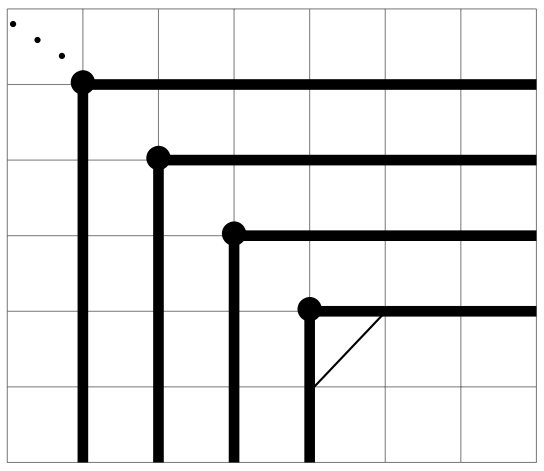

Note that the corresponding $\hat I$-tuple is $\bs m_{C^+} =  (1,\dots,(z-a)^{-1},\dots, 1)$, where $a = c(\boxxslash)$. In particular, the
shift is $\bs s=(0,\dots,0,-1,0,\dots,0)$, so it is  a module for the shifted algebra $\E_{n,\bs s}$. The shifted algebra $\E_{n,\bs s}$ with antidominant shift $\bs s$, $s_i\leq 0,$ contains the Borel subalgebra of $\mathcal E_n$. Our module restricted to the Borel subalgebra is the
negative prefundamental module in terminology of \cite{HJ}. The $\E_n$ MacMahon module of \cite{FJMM1} is obtained by tensoring with a one-dimensional $\E_{n,-s}$ module $L_{(1,\dots,1,(zK-K^{-1}a),1,\dots,1)}$, where $K\in\C^\times$ is a generic complex number.

We note that a MacMahon configuration can be modified to a representation of $\E_n$ by tensoring with a one-dimensional module if and only if the shift $\bs s$ is antidominant. In the language of \cite{FJMM1}, that corresponds to the colorless vertical boundary partition.

Similarly, let $C$ be an infinite chain configuration such that the positive tail is empty while the negative tail is of $hv$-type. We call such a $C$ the dual MacMahon configuration. The MacMahon configurations give highest $\ell$-weight modules while the dual MacMahon configuration corresponds to lowest $\ell$-weight modules.

\end{example}

\begin{example} 
Recall that $\E_n$ has two vector modules, see  Section \ref{hh vv chains},  which
are given by a single $hh$ or $vv$ chain. In the language of plane partitions, the basis of the vector modules is given by plane partitions of height one over a horizontal or vertical half-line. The algebra $\E_1$  has a third vector representation, see \cite{FJMM5}. This vector representation has a basis of generalized plane partitions over a single cube. In other words, the basis is labeled by infinite towers of boxes placed in three-dimensional space vertically on the top of each other.

The algebra $\E_n$ does not have such a module. However, the shifted algebra $\E_{n, \bs s},$ where $\bs s=(0,\dots,0,1,-2,1,0,\dots,0)$,
does. We call these modules the vertical vector representations (here ``vertical" refers to the height of the plane partition, rather than the vertical direction in the plane).

Figure \ref{Infinite same monomial} gives two infinite chain configurations for a vertical vector representation. The left one has positive and negative \(vh\)-tails and the right one has positive and negative \(hv\)-tails. Note, in particular, that in contrast to the finite chain configurations, two different infinite chain configurations may produce the same module.

\begin{figure}[H]
\begin{tikzpicture}

\begin{scope}[shift={(0,-1)}]
 
\draw[step=0.5cm, gray, very thin] (-2,0.5) grid (1.5,3.5);

\node at (-1.5,3) {\textbullet};
\node at (-1, 2.5) {\textbullet};

\draw[line width = 2] (-1.5,3)--(-1.5,0.5);
\draw[line width = 2] (-1.5,3)--(1.5,3);
\draw[line width = 2] (-1,2.5)--(-1,0.5);
\draw[line width = 2] (-1,2.5)--(1.5,2.5);

\node at (-0.5, 1.5) {\textbullet};
\node at (0,1.5) {$\bm \times$};
\node at (0,2) {\textbullet};

\node at (0,1) {\textbullet};
\node at (0.5,1) {$\bm \times$};
\node at (0.5,1.5) {\textbullet};

\draw[line width = 2] (-0.5,1.5)--(-0.5,0.5);
\draw[line width = 2] (-0.5,1.5)--(0,1.5);
\draw[line width = 2] (0,1.5)--(0,2);
\draw[line width = 2] (0,2)--(1.5,2);
\draw[line width = 2] (0,1)--(0,0.5);
\draw[line width = 2] (0,1)--(0.5,1);
\draw[line width = 2] (0.5,1)--(0.5,1.5);
\draw[line width = 2] (0.5,1.5)--(1.5,1.5);

\node at (0.9, 0.85) {$\ddots$};
\node at (-1.8, 3.4) {$\ddots$};

\end{scope}

\begin{scope}[shift={(6,-1)}]
\draw[step=0.5cm, gray, very thin] (-2,0.5) grid (1.5,3.5);

\node at (-1,2) {$\bm \times$};
\node at (-1,2.5) {\textbullet};
\node at (-0.5,2.5) {$\bm \times$};

\node at (-1.5,2.5) {$\bm \times$};
\node at (-1.5,3) {\textbullet};
\node at (-1,3) {$\bm \times$};

\draw[line width = 2] (-1.5,2.5)--(-2,2.5);
\draw[line width = 2] (-1.5,2.5)--(-1.5,3);
\draw[line width = 2] (-1.5,3)--(-1,3);
\draw[line width = 2] (-1,3)--(-1,3.5);

\draw[line width = 2] (-1,2)--(-2,2);
\draw[line width = 2] (-1,2)--(-1,2.5);
\draw[line width = 2] (-1,2.5)--(-0.5,2.5);
\draw[line width = 2] (-0.5,2.5)--(-0.5,3.5);

\node at (0,1.5) {$\bm \times$};
\node at (0.5,1) {$\bm \times$};

\draw[line width = 2] (0,1.5)--(-2,1.5);
\draw[line width = 2] (0,1.5)--(0,3.5);

\draw[line width = 2] (0.5,1)--(-2,1);
\draw[line width = 2] (0.5,1)--(0.5,3.5);

\node at (0.9, 0.85) {$\ddots$};
\node at (-1.8, 3.4) {$\ddots$};

\end{scope}

\end{tikzpicture}
\caption{Two configurations for a vertical vector representation.}
\label{Infinite same monomial}
\end{figure}
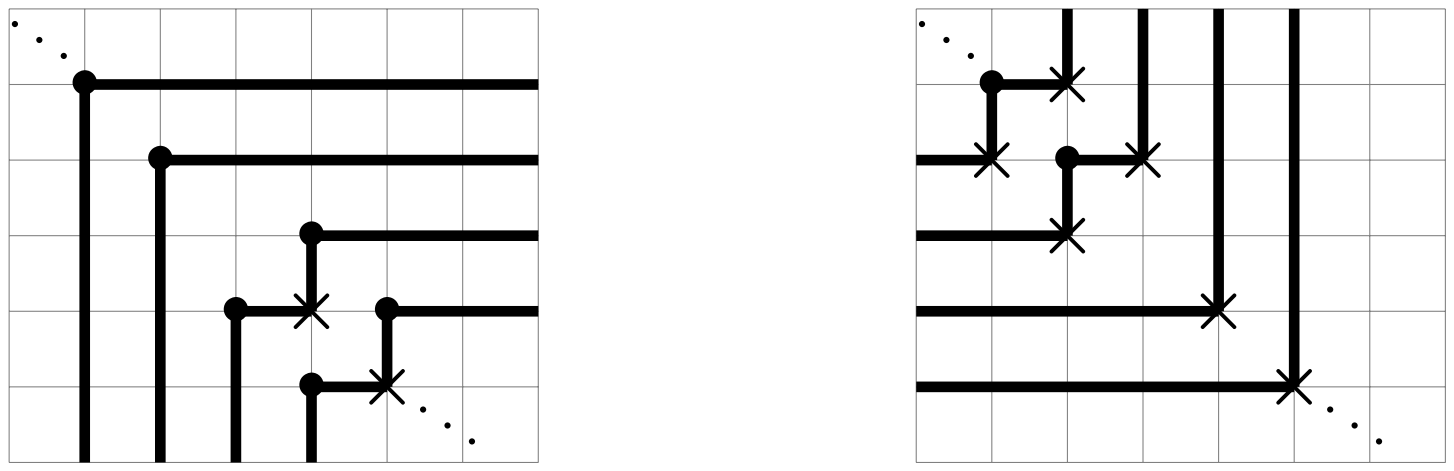

The expansions correspond to adding a cube to or removing a cube from the infinite tower of boxes over the same box. In particular, all expansions in the vertical vector representation have the same color.

\end{example}

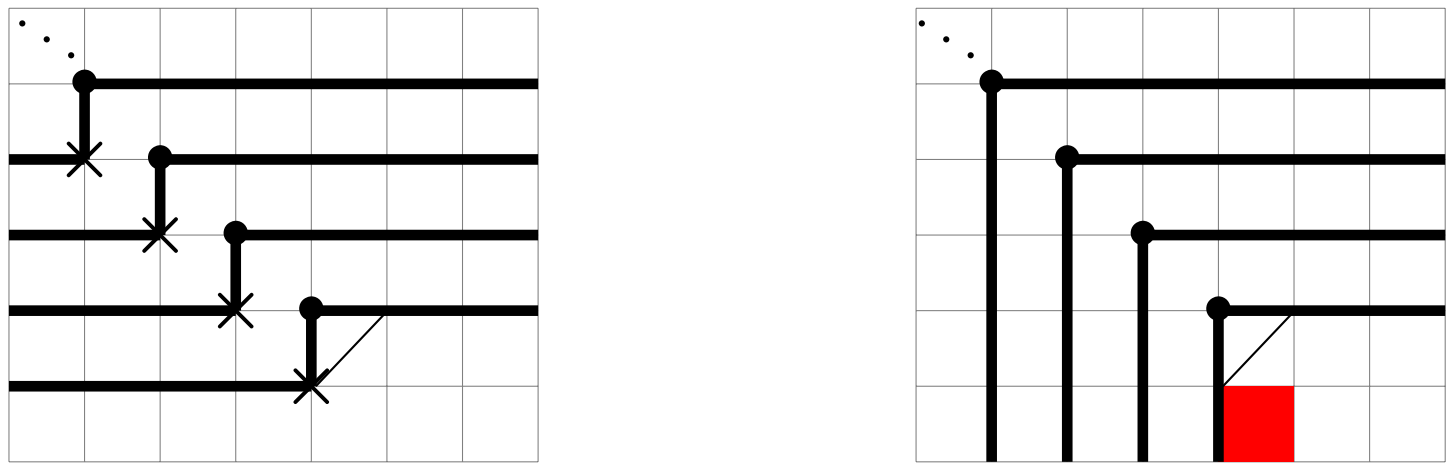
\begin{figure}[H]   
\begin{tikzpicture}

\begin{scope}[shift={(0,0)}]
\draw[step=0.5cm, gray, very thin] (-2,0.5) grid (1.5,3.5);

\node at (-1.75, 3.4) {$\ddots$};

\node at (-1.5,2.5) {$\bm \times$};
\node at (-1,2) {$\bm \times$};
\node at (-0.5,1.5) {$\bm \times$};
\node at (0,1) {$\bm \times$};
\node at (-1.5,3) {\textbullet};
\node at (-1,2.5) {\textbullet};
\node at (-0.5,2) {\textbullet};
\node at (0,1.5) {\textbullet};
\draw (0.03,1)--(0.505,1.5);

\draw[line width = 2] (-2,2.5)--(-1.5,2.5);
\draw[line width = 2] (-1.5,2.5)--(-1.5, 3);
\draw[line width = 2] (-1.5,3)--(1.5,3);
\draw[line width = 2] (-1,2)--(-2,2);
\draw[line width = 2] (-1,2)--(-1,2.5);
\draw[line width = 2] (-1,2.5)--(1.5,2.5);
\draw[line width = 2] (-0.5,1.5)--(-2,1.5);
\draw[line width = 2] (-0.5,1.5)--(-0.5,2);
\draw[line width = 2] (-0.5,2)--(1.5,2);
\draw[line width = 2] (0,1)--(-2,1);
\draw[line width = 2] (0,1)--(0,1.5);
\draw[line width = 2] (0,1.5)--(1.5,1.5);

\end{scope}

\begin{scope}[shift={(6,0)}]
 
\draw[step=0.5cm, gray, very thin] (-2,0.5) grid (1.5,3.5);

\node at (-1.5,3) {\textbullet};
\node at (-1, 2.5) {\textbullet};
\node at (-0.5,2) {\textbullet};
\node at (0, 1.5) {\textbullet};
\redbox at (0,0.5);

\draw[line width = 2] (-1.5,3)--(-1.5,0.5);
\draw[line width = 2] (-1.5,3)--(1.5,3);
\draw[line width = 2] (-1,2.5)--(-1,0.5);
\draw[line width = 2] (-1,2.5)--(1.5,2.5);
\draw[line width = 2] (-0.5,2)--(-0.5,0.5);
\draw[line width = 2] (-0.5,2)--(1.5,2);
\draw[line width = 2] (0,1.5)--(0,0.5);
\draw[line width = 2] (0,1.5)--(1.5,1.5);

\draw (0.03,1)--(0.505,1.5);

\node at (-1.8, 3.4) {$\ddots$};

\end{scope}

\end{tikzpicture}
\caption{Two configurations for a vertical Fock module.}
\label{vertical fock}
\end{figure}

\begin{example}
Recall that a Fock module, see  Section \ref{vh hv chains},
is given by a single $vh$ chain. In the language of plane partitions, the basis of the Fock module is given by plane partitions of height one. The algebra $\E_1$  has two other Fock modules, see \cite{FJMM5}. These Fock modules have a basis of plane partitions over a horizontal or vertical half-line. In other words, the basis is labeled by usual partitions placed in the three-dimensional space vertically over the half-line. 

The algebra $\E_n$ does not have such a module. However, for $\bs s=(0,\dots,0,1,-1,0,\dots,0)$, the shifted algebra $\E_{n, \bs s}$ does, see  Figure 29 and Appendix E.2.2 in \cite{NW}. We call such modules the vertical Fock modules (not to be confused with vertical lines).

Figure \ref{vertical fock} shows two configurations yielding the same vertical Fock module. The first has a positive \(hh\)-tail and an empty negative tail; the second is obtained by truncating the MacMahon configuration. 

\end{example}

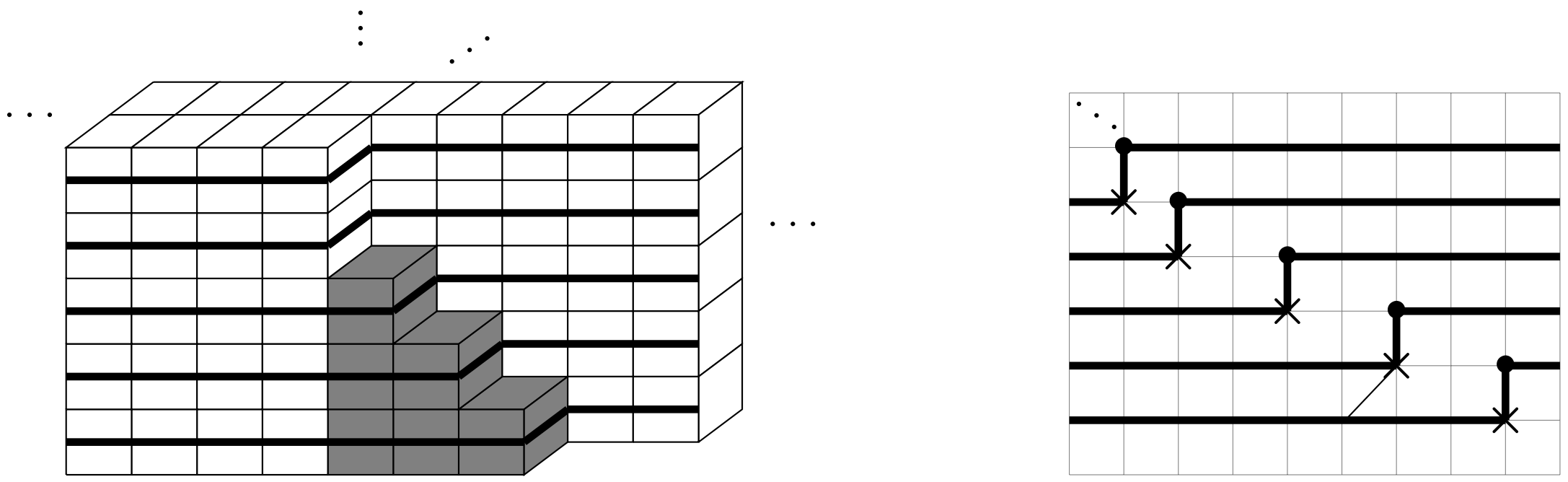
\begin{figure}[H]

\begin{tikzpicture}

\begin{scope}[shift={(0,0)}]

% first layer first row front

\fbboxh at (-0.6,0);
\fbboxh at (-1.2,0);
\fbboxh at (-1.8,0);
\fwboxh at (-2.4,0);
\fwboxh at (-3,0);
\fwboxh at (-3.6,0);
\fwboxh at (-4.2,0);

%first layer first row  row top
\bbox at (-0.2,0.9);
%\bbox at (-0.8, 0.9);

% first layer first row side

\rbboxh at (0,0);

% first layer second row front

\fwboxh at (-4.2,0.6);
\fwboxh at (-3.6,0.6);
\fwboxh at (-3,0.6);
\fwboxh at (-2.4,0.6);
\fbboxh at (-1.8,0.6);
\fbboxh at (-1.2,0.6);

% first layer second row side

\rbboxh at (-0.6,0.6);

% first layer second row top

%\bbox at (-1.4,1.5);
\bbox at (-0.8,1.5);

% first layer third row front

\fwboxh at (-4.2,1.2);
\fwboxh at (-3.6,1.2);
\fwboxh at (-3,1.2);
\fwboxh at (-2.4,1.2);
\fbboxh at (-1.8,1.2);

% first layer third row side

\rbboxh at (-1.2,1.2);

% first layer third row top

\bbox at (-1.4,2.1);

% first layer fourth row front

\fwboxh at (-4.2,1.8);
\fwboxh at (-3.6,1.8);
\fwboxh at (-3,1.8);
\fwboxh at (-2.4,1.8);

% first layer fourth row side

\rwboxh at (-1.8,1.8);

% first layer fifth row front 

\fwboxh at (-4.2,2.4);
\fwboxh at (-3.6,2.4);
\fwboxh at (-3,2.4);
\fwboxh at (-2.4,2.4);

% first layer fifth row side

\rwboxh at (-1.8,2.4);

% first layer fifth row top
\wbox at (-3.8,3.3);
\wbox at (-3.2,3.3);
\wbox at (-2.6,3.3);
\wbox at (-2,3.3);

% second layer fifth row top

\wbox at (-3.4, 3.6);
\wbox at (-2.8, 3.6);
\wbox at (-2.2, 3.6);
\wbox at (-1.6, 3.6);
\wbox at (-1,3.6);
\wbox at (-0.4, 3.6);
\wbox at (0.2, 3.6);
\wbox at (0.8, 3.6);
\wbox at (1.4, 3.6);

%second layer fifth row side

\rwbox at (1.6,2.7)

% second layer fifth row front

\fwboxh at (-1.4, 2.7);
\fwboxh at (-0.8, 2.7);
\fwboxh at (-0.2, 2.7);
\fwboxh at (0.4, 2.7);
\fwboxh at (1, 2.7);

%second layer fourth row side
\rwbox at (1.6,2.1)

% second layer fourth row front
\fwboxh at (-1.4, 2.1);
\fwboxh at (-0.8, 2.1);
\fwboxh at (-0.2, 2.1);
\fwboxh at (0.4, 2.1);
\fwboxh at (1, 2.1);

%second layer third row side

\rwbox at (1.6,1.5)

% second layer third row front

%\fwboxh at (-1.4, 1.5);
\fwboxh at (-0.8, 1.5);
\fwboxh at (-0.2, 1.5);
\fwboxh at (0.4, 1.5);
\fwboxh at (1, 1.5);

% second layer second row side

\rwbox at (1.6,0.9)

%second layer second row front

%\fwboxh at (-0.8, 0.9);
\fwboxh at (-0.2, 0.9);
\fwboxh at (0.4, 0.9);
\fwboxh at (1, 0.9);

% second layer first row side

\rwbox at (1.6,0.3)

% second layer first row front

\fwboxh at (0.4, 0.3);
\fwboxh at (1, 0.3);

\node at (-0.5, 4) {\reflectbox{$\ddots$}};

\node at (-4.5, 3.3) {$\dots$};

\node at (2.5, 2.3) {$\dots$};

\node at (-1.5, 4.2) {$\vdots$};

\end{scope}

\begin{scope}[shift={(7,0)}]
\draw[step=0.5cm, gray, very thin] (-2,0) grid (2.5,3.5);

\node at (-1.75, 3.4) {$\ddots$};
\node at (-1.5,2.5) {$\bm \times$};
\node at (-1,2) {$\bm \times$};
\node at (0,1.5) {$\bm \times$};
\node at (1,1) {$\bm \times$};
\node at (-1.5,3) {\textbullet};
\node at (-1,2.5) {\textbullet};
\node at (0,2) {\textbullet};
\node at (1,1.5) {\textbullet};
\node at (2,0.5) {$\bm \times$};
\node at (2,1) {\textbullet};
\draw (0.03+0.5,1-0.5)--(0.505+0.5,1.5-0.5);

\draw[line width = 2] (-2,2.5)--(-1.5,2.5);
\draw[line width = 2] (-1.5,2.5)--(-1.5, 3);
\draw[line width = 2] (-1.5,3)--(2.5,3);
\draw[line width = 2] (-1,2)--(-2,2);
\draw[line width = 2] (-1,2)--(-1,2.5);
\draw[line width = 2] (-1,2.5)--(2.5,2.5);
\draw[line width = 2] (0,1.5)--(-2,1.5);
\draw[line width = 2] (-0,1.5)--(-0,2);
\draw[line width = 2] (0,2)--(2.5,2);
\draw[line width = 2] (1,1)--(-2,1);
\draw[line width = 2] (1,1)--(1,1.5);
\draw[line width = 2] (1,1.5)--(2.5,1.5);
\draw[line width = 2] (2,1)--(2.5,1);
\draw[line width = 2] (2,0.5)--(-2,0.5);
\draw[line width = 2] (2,0.5)--(2,1);

\end{scope}

\end{tikzpicture}
\caption{Vertical Fock space (3,2,1).}
\label{vertical fock 321}
\end{figure}

We illustrate the name ``vertical" by giving the state corresponding to a partition $(3,2,1)$ as a generalized plane partition in Figure \ref{vertical fock 321}. The state $(3,2,1)$ in the horizontal Fock module is shown in Figure \ref{vh chain plane partition}.

\begin{example}
The slanted $\E_2$ module shown in Figure 8  in \cite{JM} is a tensor product of one vertical vector representation and two vertical Fock modules. In particular, the sum of all three shifts is $(-2,2)+(1,-1)+(1,-1)=(0,0)$.

The other $\E_2$  slanted modules, see Figure 9 in \cite{JM}, can be obtained as tensor products of two vertical Fock modules with a module corresponding to an infinite chain configuration with appropriate positive and negative $vh$ tails.

\end{example}

{\bf Acknowledgments.\ }
The authors are partially supported by Simons Foundation grant number \#709444.

The paper was written without AI assistance. ChatGPT 6 Ultra was used to proofread the final draft.

\end{document}